\documentclass[12pt]{article}
\usepackage{graphicx} 
\usepackage{amsmath, amsthm, amssymb}
\usepackage[shortlabels]{enumitem}
\usepackage{mathtools}
\usepackage{mathrsfs}
\usepackage{tikz, tikz-cd}
\usepackage[colorlinks]{hyperref}
\usepackage[noabbrev,capitalize]{cleveref}
\usepackage{float}

\hypersetup{
    colorlinks=true,
    linkcolor=purple,
    citecolor=orange,
    filecolor=magenta,      
    urlcolor=cyan,
    pdftitle={A new tautological relation},
    }
    
\usepackage{fullpage}
\usepackage{caption}

\newcommand{\ghost}{\mathsf{Ghost}}
\newcommand{\Mbar}{\overline{\mathcal{M}}}

\newcommand{\M}{\mathcal{M}}
\newcommand{\CC}{\mathcal{C}}
\newcommand{\SSS}{\Sigma}
\newcommand{\Mb}{\Mbar}
\newcommand{\R}{\mathsf{R}}
\newcommand{\RH}{\mathsf{RH}}
\newcommand{\CH}{\mathsf{CH}}
\newcommand{\HH}{\mathsf{H}}
\newcommand{\QQQ}{\mathsf{Q}}
\newcommand{\CCC}{\mathsf{C}}
\renewcommand{\L}{\mathsf{L}}
\newcommand{\pr}{\mathsf{pr}}
\newcommand{\Qinv}{\QQQ_{\mathrm{inv}}}
\newcommand{\AAA}{\mathsf{A}}
\newcommand{\qq}{\mathbb{Q}}

\newcommand{\C}{\mathbb C}
\newcommand{\Q}{\qq}
\newcommand{\gr}{\operatorname{gr}}
\newcommand{\ct}{\mathrm{ct}}
\newcommand{\Sym}{\mathrm{Sym}}
\newcommand{\Aut}{\mathrm{Aut}}

\newcommand{\BM}{\mathrm{BM}}
\newcommand{\rt}{\mathrm{rt}}
\newcommand{\ThreeSpin}{\mathsf{3Spin}}
\DeclareMathOperator{\id}{id}

\newtheorem{thm}{Theorem}[section]
\newtheorem{prop}[thm]{Proposition}
\newtheorem{prediction}[thm]{Prediction}
\newtheorem{lem}[thm]{Lemma}
\newtheorem{cor}[thm]{Corollary}

\newtheorem{conjecture}[thm]{Conjecture}

\theoremstyle{definition}
\newtheorem{definition}[thm]{Definition}

\newtheorem{rem}[thm]{Remark}
\newtheorem{question}[thm]{Question}

\usepackage{bead-diagrams}

\newcommand{\GraphBead}[3]{%
  \node[beads/bead,beads/#1] (#2) at (#3) {};
  \node[beads/#1,beads/letter] at (#2.center)
    {\csname GraphLetter#1\endcsname};
}

\tikzset{
  graphs/dots/.style={
    rectangle, draw=none, fill=none, inner xsep=1mm,
    inner ysep=0pt, outer sep=0pt, font=\normalsize
  },
  moregraphs/length/.style={beads/label,text=black,font=\small},
  moregraphs/center/.style={
    baseline={([yshift=-.5ex]current bounding box.center)}
  },
  moregraphs/vertex/.style={circle,fill=black,inner sep=1.7pt},
  moregraphs/edge/.style={thick},
  moregraphs/leg/.style={thick}
}

\newcommand{\MoreGraphLeg}[4]{%
  \draw[beads/edge] (#1.#2)
    -- ++(#2:\pgfkeysvalueof{/tikz/beads/leg length})
    node[beads/label,anchor=#3] {#4};
}

\newcommand{\GraphPath}[3]{%
  \GraphBead{P}{#1-left}{2,#2}
  \GraphBead{P}{#1-right}{5,#2}
  \node[graphs/dots] (#1-dots) at (3.5,#2) {$\cdots$};

  \ifnum#3=1\relax
    \GraphBead{A}{#1-extra}{1,{#2/2}}
    \draw[beads/edge]
      (v) -- (#1-extra)
      (#1-extra) -- (#1-left);
  \else
    \draw[beads/edge] (v) -- (#1-left);
  \fi

  \draw[beads/edge]
    (#1-left) -- (#1-dots.west)
    (#1-dots.east) -- (#1-right)
    (#1-right) -- (w);
}

\newcommand{\GraphVertexNames}{%
  \node[beads/label,anchor=south] at (v.north) {$v$};
  \node[beads/label,anchor=south] at (w.north) {$w$};
}

\newcommand{\GraphLegs}{%
  \draw[beads/edge] (v.west)
    -- ++(180:\pgfkeysvalueof{/tikz/beads/leg length})
    node[beads/label,anchor=east] {1};
  \draw[beads/edge] (w.east)
    -- ++(0:\pgfkeysvalueof{/tikz/beads/leg length})
    node[beads/label,anchor=west] {2};
}

\newcommand{\ThreePathGraph}{%
  \begin{tikzpicture}[beads/defaults,baseline=-.5ex]
    \GraphBead{O}{v}{0,0}
    \GraphBead{P}{w}{7,0}
    \GraphPath{upper}{1.5}{0}
    \GraphPath{middle}{0}{0}
    \GraphPath{lower}{-1.5}{0}
    \GraphVertexNames
  \end{tikzpicture}%
}

\newcommand{\TwoPathGraph}[4]{%
  \begin{tikzpicture}[beads/defaults,baseline=-.5ex]
    \GraphBead{#1}{v}{0,0}
    \GraphBead{#2}{w}{7,0}
    \GraphPath{upper}{1.5}{#3}
    \GraphPath{lower}{-1.5}{#4}
    \GraphVertexNames
    \GraphLegs
  \end{tikzpicture}%
}

\usepackage{array}
\usepackage[table]{xcolor}

\newcommand{\BracketedBead}[1]{\left[\BeadDiagram{#1}\right]}

\newlength{\BeadNumberWidth}
\newlength{\BeadSumWidth}
\newlength{\BeadHalfSumWidth}
\newcolumntype{N}{%
  >{\columncolor{black!10}\centering\arraybackslash\sffamily\bfseries}%
  m{\BeadNumberWidth}}

\newcommand{\BeadSum}[2][\BeadSumWidth]{%
  \begin{minipage}[c]{#1}
    \vspace{3mm}
    \(\begin{aligned}#2\end{aligned}\)
    \par\vspace{3mm}
  \end{minipage}%
}

\usetikzlibrary{arrows.meta,decorations.pathmorphing}
\newcommand{\BeadRewriteArrow}{%
  \mathrel{%
    \begin{tikzpicture}[baseline={-\the\fontdimen22\textfont2}]
      \draw[
        line width=.5pt,
        -{Stealth[length=1.5mm,width=1.2mm]},
        decorate,
        decoration={snake,amplitude=.35mm,segment length=2.4mm,
                    pre length=.4mm,post length=1.8mm}
      ] (0,0) -- (14mm,0);
    \end{tikzpicture}%
  }%
}

\newcommand{\BeadInput}[1]{%
  \left[
    \BeadDiagram[beads/labels={$\boxed{1}$,$\boxed{2}$,$\boxed{3}$}]{#1}
  \right]%
}

\newcommand{\BeadPair}[2]{%
  \left[
    \vcenter{\hbox{%
      $\displaystyle
        \BeadDiagram[
          baseline=0pt,
          beads/labels={$\boxed{1}$,$\star$,}
        ]{#1}\,\,
        \BeadDiagram[
          baseline=0pt,
          beads/labels={$\star$,$\boxed{2}$,$\boxed{3}$}
        ]{#2}
      $%
    }}%
  \right]%
}

\title{A new tautological relation on $\Mbar_{2,14}$ in cohomology}
\author{Samir Canning, Bochao Kong, Rahul Pandharipande,\\ Dan Petersen,  Aaron Pixton, Johannes Schmitt}
\date{November 2026}

\begin{document}

\maketitle

\begin{abstract}
We construct a relation $\QQQ$ in the cohomological tautological ring 
$\RH^*(\Mbar_{2,14})$
which lies outside the space of $3$-spin relations. $\QQQ$ is a simple
linear combination of pure boundary strata classes of codimension 9 supported on the divisor $\Delta_0 \subset \Mbar_{2,14}$ of irreducible nodal curves,
and transforms according to the irreducible representation of $\SSS_{14}$ associated to the partition $(2^7)$.
The result shows that the set of $3$-spin relations cannot be complete in cohomology. Motivic conjectures relating algebraic cycles to special values of $L$-functions suggest that $\QQQ$ is nonzero in Chow. More precisely, we prove that the Chow nonvanishing of $\QQQ$ 
follows from a conjecture of Harder concerning the failure of the Manin--Drinfeld principle in genus 2.

We prove that $\QQQ$ vanishes in cohomology by pairing with generators of the tautological ring. To prove $\QQQ$ is not in the set of $3$-spin relations requires
an argument 
involving the $\SSS_{14}$-representation-theoretic properties of the $3$-spin relations
and the cohomology of local systems on $\mathcal{A}_2$ and $\mathcal{M}_2$. 
Parallel constructions for other moduli spaces $\Mbar_{g,n}$ 
are discussed in the paper. There are many open questions. 
\end{abstract}

\tableofcontents

\section{Introduction}
\subsection{Tautological relations}
Let $\RH^*(\Mbar_{g,n}) \subset \HH^{2*}(\Mbar_{g,n})$ be the tautological ring in the cohomology{\footnote{All cohomology, homology, and Chow groups of $\Mbar_{g,n}$ are taken with $\mathbb{Q}$-coefficients.}} of
the moduli spaces of stable curves of genus $g$ with $n$ markings. 
The ideal of cohomological relations, by 
definition, is
the kernel of the canonical $\mathbb{Q}$-algebra homomorphism from the strata algebra to
the tautological ring:
$$\mathsf{CohRel}^*_{g,n} = \mathsf{Ker}\Big( \mathsf{S}^*_{g,n} 
\rightarrow \mathsf{RH}^*(\Mbar_{g,n})\Big)\, .$$ Similarly,
the ideal of relations in Chow is defined by:
$$\mathsf{ChowRel}^*_{g,n} = \mathsf{Ker}\Big( \mathsf{S}^*_{g,n} 
\rightarrow \mathsf{R}^*(\Mbar_{g,n})\Big)\, ,$$
where $\R^*(\Mbar_{g,n}) \subset \CH^{*}(\Mbar_{g,n})$ is the tautological ring in Chow.
Using the cycle class map, 
$$\mathsf{ChowRel}^*_{g,n} \subset \mathsf{CohRel}^*_{g,n}\, .$$

The investigation of the ideal of tautological relations goes back to Mumford \cite{DMumford}.
We refer the reader to \cite{FP3,P1} for an introduction to the study of tautological classes. In addition to governing the tautological ring,  the relations
also constrain Gromov-Witten theory and play a role in
several other directions in mathematics.

An explicit ideal, $\mathsf{3Spin}^*_{g,n} \subset \mathsf{S}_{g,n}^*$, was proposed by Pixton in \cite{Pixton2012-ConjecturalRelations} 
and was proven  to yield relations in both cohomology and  Chow \cite{Janda,PPZ}:
$$\mathsf{3Spin}^*_{g,n}\subset \mathsf{ChowRel}^*_{g,n} \subset \mathsf{CohRel}^*_{g,n}\, .$$
Pixton conjectured\footnote{OpenAI \cite{OpenAI} has proposed a counterexample to Pixton's conjecture \eqref{pixconj} in Chow for $g=10^{60}$ and \\$n=3\binom{10^{60}}{3}$ with
methods unrelated to our cohomological constructions for $g=2$
and $n=14$.}
\begin{equation} \label{pixconj}
\mathsf{3Spin}^*_{g,n}= \mathsf{ChowRel}^*_{g,n} 
\end{equation}
for all $g$ and $n$. A folklore conjecture which has generally been believed is the equality
\begin{equation} \label{cycleconj}
\mathsf{ChowRel}^*_{g,n} = \mathsf{CohRel}^*_{g,n}
\end{equation}
for all $g$ and $n$.

\subsection{The new relation \texorpdfstring{$\QQQ$}{\QQQ}} \label{newrel}

The main result here is a new tautological relation $\mathsf{Q} \in \mathsf{CohRel}^9_{2,14}$.
The formula for $\mathsf{Q}$ is a linear combination of pure boundary strata classes. The construction requires a family of stable graphs for $\Mbar_{2,14}$ 
with 8 vertices (all of genus 0) and $9$ edges.

To start, we color the markings $\{1, \ldots, 7\}$ red and the markings $\{8,\ldots, 14\}$ blue. 
We fix the pairing of
the red and blue markings given by
\begin{equation} \label{pppr}
(1,8)\, , \ (2,9)\, , \ldots, \, (7,14)\, .
\end{equation}
For positive integers $a,b,c$ with sum $9$, let $\Theta_{abc}$ be the {\em subdivided theta} graph consisting of two vertices $v,w$ connected by three parallel paths of lengths $a,b,c$, with one red leg and one blue leg (both unlabeled) on every vertex other than $v$. For example, $\Theta_{135}$ is:

\vspace{10pt}
\begin{center}
    
\begin{tikzpicture}[
    vertex/.style={circle,fill=black,inner sep=1.7pt},
    edge/.style={thick},
    leg/.style={thick},
]

\node[vertex,label=left:$v$]  (v) at (0,0) {};
\node[vertex,label=right:$w$] (w) at (6,0) {};

\draw[edge] (v) -- (w);

\node[vertex] (u1) at (2,1.6) {};
\node[vertex] (u2) at (4,1.6) {};

\draw[edge]
    (v) -- (u1) -- (u2) -- (w);

\node[vertex] (d1) at (1.2,-1.5) {};
\node[vertex] (d2) at (2.4,-2.0) {};
\node[vertex] (d3) at (3.6,-2.0) {};
\node[vertex] (d4) at (4.8,-1.5) {};

\draw[edge]
    (v) -- (d1) -- (d2) -- (d3) -- (d4) -- (w);


\foreach \x in {u1,u2}{
    \draw[leg,red]  (\x) -- ++(-0.28,0.55);
    \draw[leg,blue] (\x) -- ++( 0.28,0.55);
}

\foreach \x in {d1,d2,d3,d4}{
    \draw[leg,red]  (\x) -- ++(-0.22,-0.55);
    \draw[leg,blue] (\x) -- ++( 0.22,-0.55);
}

\draw[leg,red]  (w) -- ++(0.45,0.45);
\draw[leg,blue] (w) -- ++(0.45,-0.45);

\end{tikzpicture}
\end{center}
\vspace{10pt}
\noindent We define the class $[\Theta_{abc}]\in\mathsf{S}^9_{2,14}$ by the following procedure.
\begin{enumerate}
    \item[(i)] Place marking numbers of the corresponding 
    colors on the red and blue legs in a way which respects the pairing \eqref{pppr}. The result is a stable graph $\Gamma$ for $\Mbar_{2,14}$. 
    Let $[\Gamma]\in \mathsf{S}^9_{2,14}$ be the associated class in the strata algebra.
    \item[(ii)] Anti-symmetrize with respect to the $\SSS_7$-action
on the red markings and the $\SSS_7$-action on the blue markings:
\begin{equation}\label{fres} [\Theta_{abc}] = \frac{1}{7!\cdot 7!}
\sum_{(\sigma,\tau)\in \SSS_7\times \SSS_7\subset \SSS_{14}}
\mathsf{sgn}(\sigma) \cdot\mathsf{sgn}(\tau) \cdot (\sigma,\tau)[\Gamma]\, .
\end{equation}
\end{enumerate}
The final result  \eqref{fres} is independent of the initial placing (i).

\begin{definition} \label{q-def} Let $\mathsf{Q} \in \mathsf{S}^9_{2,14}$ be given by
$$\mathsf{Q} = [\Theta_{135}]-[\Theta_{144}]-[\Theta_{225}]+2[\Theta_{234}]-[\Theta_{333}]\, .$$
\end{definition}

\begin{thm}\label{qth_rel}
We have $\mathsf{Q} \in \mathsf{CohRel}^9_{2,14}$.
\end{thm}

\begin{thm} \label{qth_not3spin} We have $\mathsf{Q} \notin \mathsf{3Spin}^9_{2,14}$.
\end{thm}

As a consequence of Theorem~\ref{qth_rel} and Theorem~\ref{qth_not3spin}, at least one of the two conjectured equalities 
\eqref{pixconj} and \eqref{cycleconj} must not hold.\footnote{We expect \emph{both} conjectured equalities to be false in general: $\QQQ$ likely falsifies \eqref{cycleconj} and OpenAI claims to falsify \eqref{pixconj}.}
The $\mathbb{Q}$-vector spaces $$\mathsf{3Spin}^9_{2,14} \subset \mathsf{ChowRel}^9_{2,14} \subset \mathsf{CohRel}^9_{2,14} \subset \mathsf{S}^9_{2,14}$$ are naturally representations of $\SSS_{14}$.
Since $\QQQ$ is defined by
the $\SSS_7\times \SSS_7$ anti-symmetrization of an $(\SSS_2)^7$-invariant{\footnote{The subgroup $(\SSS_2)^7\subset \SSS_{14}$ is
determined by the pairing \eqref{pppr}.}} class, $\QQQ$ generates a copy of the irreducible $\SSS_{14}$-representation of type $(2^7)$ in $\mathsf S_{2,14}^9$.   Theorems~\ref{qth_rel} and \ref{qth_not3spin} therefore imply: 
\begin{cor} \label{qth2} 
The class $\QQQ$
generates an irreducible $\SSS_{14}$-subrepresentation 
$$\mathcal{Q}\subset \mathsf{CohRel}^9_{2,14}/ \mathsf{3Spin}^9_{2,14}$$ 
of type $(2^7)$.
\end{cor}

\noindent Concerning the Chow tautological relations,
there are two distinct possibilities:
\begin{enumerate}
\item[(i)] $\mathcal{Q} \subset
\mathsf{ChowRel}^9_{2,14}/\mathsf{3Spin}^9_{2,14}$, 
\item[(ii)]
$\mathcal{Q} \subset
\mathsf{CohRel}^9_{2,14}/\mathsf{ChowRel}^9_{2,14}$,
\end{enumerate}
which imply the failure of the conjectural equalities \eqref{pixconj} and \eqref{cycleconj}, respectively, in the $(2^7)$-isotypical component. To settle the dichotomy, the following question must be answered.

\begin{question}\label{Chowvan} Is $\mathsf{Q} \in \mathsf{ChowRel}^9_{2,14}$? 
\end{question}

\noindent We expect the answer to Question \ref{Chowvan} to be {\em no}, and, therefore, we expect 
possibility (ii) holds. 

Since $\QQQ$ is a linear
combination of strata classes with vertices all of genus 0, $\QQQ$ yields a universal genus 0 relation in Gromov-Witten theory by the splitting axiom. Whether new results in Gromov-Witten theory can be deduced from $\QQQ$ is an interesting direction which we do not pursue here.

In Section \ref{sec-qp},
the formula for $\QQQ$ is generalized: for each homogeneous symmetric polynomial $P$ of degree $m$ in three variables, there is an associated class $$\QQQ_P\in \mathsf{S}_{2,14+2m}^{9+m}\, ,$$ 
with $\QQQ_1=6\QQQ$. 
Which of the above properties of $\QQQ$ hold for the other classes $\QQQ_P$?
Some results can be found in Sections \ref{secvan} and \ref{sec:boundarynonvanishing}.

\subsection{Harder's conjecture and the Chow class of \texorpdfstring{$\QQQ$}{Q}}

The question whether $\QQQ$ vanishes in Chow has an arithmetic interpretation. The motive $S[18]$ associated to the unique cusp form of weight $18$ for $\mathrm{SL}_2(\mathbb Z)$ occurs in $\HH^{17}(\Mbar_{2,14})$, transforming as the $\SSS_{14}$-representation $(2^7)$. Its $L$-function has a simple zero at the central point. Beilinson's conjectures therefore predict the existence of a nonzero homologically trivial $8$-cycle on $\Mbar_{2,14}$, defined over $\Q$, transforming according to the representation $(2^7)$. Our cycle $\QQQ$ provides an explicit candidate for this cycle.

There is a more precise connection with a conjecture of Harder. Let $\mathcal A_2$ be the moduli space of principally polarized abelian surfaces. In the cohomology of the local system $\mathbb V_{7,7}$ on
$\mathcal A_2$ one finds an extension
\begin{equation}\label{harder-extension}
0 \longrightarrow S[18] \longrightarrow W_{18}\HH^3(\mathcal A_2,\mathbb V_{7,7})\longrightarrow \Q(-9) \longrightarrow 0\,.
\end{equation}
Beilinson's conjectures predict likewise that $\mathrm{Ext}^1(\Q(-9),S[18])$, the extension group in ``mixed motives over $\Q$'', is 1-dimensional. A 1993 conjecture of Harder predicts that the extension\footnote{Harder's conjecture is  more general and concerns also cusp forms of higher level and other weights.} \eqref{harder-extension} is nontrivial and hence generates the motivic extension group. It turns out that this conjecture can be expressed in terms of the cycle $\QQQ$.

\begin{thm}\label{intro-harder}
The Abel--Jacobi class of the homologically trivial cycle $\QQQ$ coincides with the extension class \eqref{harder-extension}. In particular, Harder's conjecture implies that $\QQQ$ is nonzero in Chow. 
\end{thm}

Theorem \ref{intro-harder} gives a promising avenue towards proving the nonvanishing of $\QQQ$, even without assuming a conjectural category of mixed motives over $\Q$. If the extension \eqref{harder-extension} is nontrivial in mixed Hodge structures or $\ell$-adic Galois representations then this already implies nonvanishing of $\QQQ$. 

In Section \ref{chow-quest}, we explain Beilinson's conjectures and their predictions for $\Mbar_{g,n}$ in case $g = 1$ and $g=2$. We then prove \cref{intro-harder} and explain why there are good reasons to believe in Harder's conjecture.\footnote{AI systems run by B.K.\ and Yuan Lu have independently proposed different arguments, both involving substantial computations, purporting to prove that the crystalline realization of the extension \eqref{harder-extension} is non-split at the prime $p=19$, implying in particular nonvanishing of $\QQQ$ in Chow. We have not yet verified the AI arguments and do not use them anywhere in the paper. We hope, however, to return to the topic in a sequel.}

\subsection{Path of the proof of Theorem \ref{qth_rel}}
The proof of Theorem \ref{qth_rel} relies upon the equality $$\RH^*(\Mbar_{2,14}) = \HH^{2*}(\Mbar_{2,14})$$
established by Petersen and Tommasi \cite{Petersen2016-GenusTwo,petersentommasi}.
Since $\HH^{2*}(\Mbar_{2,14})$ satisfies Poincar\'e duality, 
$\mathsf{Q} \in \mathsf{CohRel}^9_{2,14}$
if and only if
\begin{equation}\label{reqvan}
\int_{[\QQQ]} \gamma = 0 
\end{equation}
for every $\gamma \in \RH^8(\Mbar_{2,14})$.
The required vanishing \eqref{reqvan} is proven in 
Section \ref{secvan}.

\subsection{Path of the proof of Theorem \ref{qth_not3spin}}

Consider the following maps of moduli spaces
$$\Mbar_{1,14\cup\{ \bullet,\star\}}/\SSS_2\,  \stackrel{\epsilon}{\longrightarrow} \, \Delta_0 \, \stackrel{\iota_0}{\longrightarrow}\, 
\Delta\, \stackrel{\iota_\Delta}{\longrightarrow}\,  \Mbar_{2,14}\, ,$$
where $\epsilon$ is the normalization of the divisor 
$\Delta_0\subset \Mbar_{2,14}$ of curves with a non-separating node, and $\Delta$ is the full boundary $\Mbar_{2,14}\smallsetminus \M_{2,14}$. 
The $\SSS_2$-action permutes the  markings $\{ \bullet,\star\}$ which are
glued to define $\epsilon$.
Both $\iota_0$ and $\iota_\Delta$ are closed 
inclusions of substacks.

Since the stable graphs which appear in Definition \ref{q-def} of
$\QQQ\in \mathsf{S}_{2,14}^9$ 
all include non-separating edges,
$\QQQ$ canonically determines
homology classes
$$\QQQ_{\Delta_0} \in \HH_{16}(\Delta_0)\ \ \ 
\text{and} \ \ \ \QQQ_\Delta \in \HH_{16}(\Delta)$$
which satisfy push-forward equations:
\begin{equation}\label{pfe2}
\iota_{0*}(\QQQ_{\Delta_0})=\QQQ_\Delta\, , \ \ \ \iota_{\Delta*}(\QQQ_\Delta)=0\, .
\end{equation}
The second equality in \eqref{pfe2}
is a consequence of
Theorem \ref{qth_rel}.

The first step in the proof of Theorem \ref{qth_not3spin} is to show the non-vanishing
$$\QQQ_{\Delta_0} \neq 0 \in \HH_{16}(\Delta_0)\, ,$$
which is deduced from the formula for $\QQQ$ in Section \ref{sec:boundarynonvanishing}.
Next, using the boundary spectral sequence for $\Mbar_{2,14}$ and properties of
the cohomology of local systems on $\mathcal{A}_2$ \cite{localsystemsA2}, we deduce the stronger non-vanishing 
$$\QQQ_{\Delta} \neq 0 \in \HH_{16}(\Delta)\, $$
in Section \ref{locsys}. 

The third step in the proof uses $\SSS_{14}$-representation theory. We have already seen that
$\QQQ\in (\mathsf{CohRel}^9_{2,14})^{[2^7]}$.
 In Section \ref{S14rep}, the $\SSS_{14}$-action on $\mathsf{3Spin}_{2,14}^9$ is proven to 
 have the following property:
$$\QQQ\in (\mathsf{3Spin}_{2,14}^9)^{[2^7]} \implies \QQQ_{\Delta} = 0 \in \HH_{16}(\Delta)\, .$$
Therefore, by contradiction, we obtain 
$\QQQ \notin \mathsf{3Spin}^9_{2,14}$.

\subsection{Propagation}

Via pushforward and pullback by appropriate gluing and forgetful maps,  $\QQQ$ can be used to produce cohomological relations outside the span of the $3$-spin relations, for $(g,n)$ other than $(2,14)$. Similarly, if $\QQQ$ is nonzero in Chow, then pushforward and pullback produces nonzero homologically trivial tautological classes, for the same range of $g$ and $n$.

\begin{thm}\label{thm:highergenus}
    Let $g \geq 2$ and $g+n \geq 16$. Then we have: 
    \begin{enumerate}
        \item[\em{(i)}] $\ThreeSpin_{g,n}^*\subsetneq \mathsf{CohRel}_{g,n}^*$.
        \item[\em{(ii)}] If $\QQQ \neq 0 \in \R^*(\Mbar_{2,14})$, then also $\mathsf{ChowRel}_{g,n}^*\subsetneq \mathsf{CohRel}_{g,n}^*$.
    \end{enumerate}
\end{thm}

  \noindent Theorem \ref{thm:highergenus} is proven in Section \ref{sec:highergenus}. Whether these constructions are optimal is not at all clear. 

\begin{question} For $g+n<16$, do the equalities
$\ThreeSpin_{g,n}^* = \mathsf{ChowRel}_{g,n}^* = \mathsf{CohRel}_{g,n}^*$ hold?
\label{smallergn}
\end{question}

\subsection{Use of Artificial Intelligence}

Motivated by the recent papers \cite{CLS, pixton2026-nongorenstein}, B.K.\ began searching for a counterexample to Pixton’s conjecture using OpenAI Codex. 
The search attempts via intersection number calculations, higher-spin, and double-ramification constructions were inconclusive. B.K.\ then prompted the AI to prove the completeness of the 3-spin relations in $\R^*(\Mbar_{2,n})$. 
However, while monitoring the milestones reported throughout the run, B.K.\ noticed that an agent was crafting a file named \verb|potential_counterexample.md|. An instruction to adversarially check the counterexample was sent. The run concluded with the AI claiming the existence of a relation in $\RH^*(\Mbar_{2,14})$ not in the span of the 3-spin relations.

The AI prompt was constructed in a style similar to that of OpenAI's CDC conjecture prompt \cite{OpenAI2026-CDCPrompt}. One GPT-6 Astra coordinated the work, while nine GPT-5.6 Sols and seven GPT-6 Astras ran in parallel to freely discuss the questions. The agents were encouraged to send messages to each other, and to do adversarial checks once a milestone was reached.

The AI's strategy was to prove nontriviality of $\gr^W_{16}\HH_c^{17}(\M_{2,14})$ in the $(2^7)$-isotypical component, using the cohomology of local systems and the boundary weight spectral sequence, and from this deduce the existence of cohomological relations outside $3$-spin. The same nontriviality is established in our \cref{reduction}, but our proof is different and shorter than the AI-generated argument, which we were unable to fully verify.

The arguments using local systems  show the existence of a relation, but do not produce formulas. However, \cref{reduction} implies the existence of a class in $\HH_{16}(\Mbar_{1,16}/\SSS_2)$ whose pushforward to $\Delta_0$ is nonzero, but which pairs to zero with every element of $\HH^{16}(\Mbar_{2,14})$. Since all even cohomology is tautological in this range, finding such a class is a high dimensional linear algebra problem, which was investigated by J.S.\ using OpenAI Codex running GPT-6 Astra to implement the linear algebra, and Claude Code running Fable 5.1 to check the implementations and improve the presentation of the resulting written reports. 
The explicit formula of Definition \ref{q-def} for the new cohomological relation $\QQQ$ is
an outcome of these calculations. The explicit formula for the $\Delta_0$-nonvanishing detector cycle $\mathsf{D}_P$ of Section~\ref{con-dp} was also found by AI-coded computer linear algebra (run by A.P.).

Our paper, however, does {\em not} depend upon any AI-generated documents, or computer linear algebra. The proofs are
all written by us and can be verified with pen and paper.
In fact, ideas very close to the results of our paper were in the air for more than a decade. The existence of nontrivial cohomological tautological relations on $\Mbar_{2,14}$, supported on $\Delta_0$, in an $\SSS_{14}$-representation of type $(2^7)$, was already explained in a lecture at ETH Z\"urich in 2013 by D.P., see \cite[Section 5]{Danslecture}, where he writes {\em it would be very useful to be able to explicitly write down these tautological relations}. 

\subsection*{Acknowledgements}
B.K. would like to thank OpenAI's Academic Researchers program for providing computational resources. 
We thank Carel Faber, Jeremy Feusi, Hannah Larson, and Yuan Lu for several early discussions about the new relation $\QQQ$. 

S.C. was supported by the SNSF Ambizione grant PZ00P2\_223473. R.P. was
supported by SNF-200020-219369 and SwissMap.
A.P. was supported by the NSF grant DMS-2301506.

\section{Proof of Theorem \ref{qth_rel}}  \label{secvan}

\subsection{The classes \texorpdfstring{$\QQQ_{P}$}{QP}} \label{sec-qp}

In Section \ref{secvan},
we will perform computations in the partially antisymmetrized strata algebra and tautological Chow ring{\footnote{We will switch to cohomology for Corollary~\ref{cor:theta_relations} only.}}
\[
\epsilon: (\mathsf{S}^*_{g,n+a+b})^{[1^a]\boxtimes[1^b]} \rightarrow
\R^*(\Mbar_{g,n+a+b})^{[1^a]\boxtimes[1^b]}\, ,
\]
where $[1^a]\boxtimes[1^b]$ indicates that we take classes simultaneously antisymmetric with respect to permuting $a$ red markings $n+1,\ldots,n+a$ and $b$ blue markings $n+a+1,\ldots,n+a+b$. Equivalently, we will work in the image of the idempotent
\[
\AAA_{n,a,b} = \frac{1}{a!\cdot b!}\sum_{\sigma\in \Sigma_a,\tau\in \Sigma_b}\mathsf{sgn}(\sigma)\mathsf{sgn}(\tau)(\sigma,\tau)
\]
in the group algebra of $\SSS_{n+a+b}$.

Fix $k\ge 1$. For any positive integers $a,b,c$ with $a+b+c=k+2$, let $\Theta_{abc}$ be the subdivided theta graph defined in the introduction (with $k$ red legs and $k$ blue legs). Let $\Gamma$ be a stable graph obtained from $\Theta_{abc}$ by labeling the red and blue legs in any way such that markings $i,k+i$ are placed on the same vertex for each $i$. The class
\[
[\Theta_{abc}]=\AAA_{0,k,k}([\Gamma])\in (\mathsf{S}^*_{2,2k})^{[1^k]\boxtimes[1^k]} 
\]
does not depend on the initial labeling.


\begin{lem}\label{lem:sep_pullbacks}
The class 
$\epsilon([\Theta_{abc}])$ vanishes on restriction to every separating boundary divisor of $\Mbar_{2,2k}$
other than the rational tails divisors in which the rational component contains exactly one red leg and one blue leg.
\end{lem}
\begin{proof}
The boundary strata of shape $\Theta_{abc}$ are disjoint from all separating boundary divisors other than the
specified rational tails divisors: $\Theta_{abc}$ has no separating edge and all the vertices are of genus zero with at most one red leg and one blue leg.
\end{proof}

\begin{prop}\label{prop:nonsep_pullback}
Let $j:\Mbar_{1,2+2k}\to\Mbar_{2,2k}$ glue the first two markings. Then, there exist classes
\[
\CCC_a\in \R^{k+2}(\Mbar_{1,2+2k})
\]
for $1\le a \le k$ such that
\[
j^*\epsilon([\Theta_{abc}]) = \CCC_a + \CCC_b + \CCC_c
\]
for all $a,b,c$ with sum $k+2$.
\end{prop} Proposition~\ref{prop:nonsep_pullback} will be proven in Section~\ref{subsec:beads}.

\begin{lem}\label{lem:rational_tail_intersection_number}
Let $\iota_k:\Mbar_{2,k}\to\Mbar_{2,2k}$ glue a rational tail at each marking $1\le i \le k$, with new markings $i$ and $i+k$. Then,
\[
\int_{\Mbar_{2,2k}}\epsilon([\Theta_{abc}])\cdot \iota_{k*}(\psi_1) = \frac{1}{k\cdot k!}\, .
\]
\end{lem}
\begin{proof}
Let $\Gamma$ be any stable graph given by labeling the red and blue legs of $\Theta_{abc}$, and let $\Gamma_\iota$ be the stable graph corresponding to $\iota_k$. We claim that there are either 1 or 0 generic $(\Gamma,\Gamma_\iota)$-graphs (as defined in \cite{GP}). More precisely, if the legs on each vertex of $\Gamma$ are of the form $\{i,i+k\}$ (for some $i$), then there is a generic $(\Gamma,\Gamma_\iota)$-graph given by bubbling off each such pair of legs, and otherwise there is no such graph.

The intersection number of $[\Gamma]$ and $(\iota_k)_*\psi_1$ is 0 in the latter case, and is equal to the integral of $\psi$ on the vertex where legs $1,1+k$ were located in the first case, which is either 0 or 1 depending on whether that vertex was valence $2$ or $3$ in the subdivided theta graph. Averaging over all such $\Gamma$ (and noting that the nonzero contributions all show up with a positive sign) gives the result.
\end{proof}

\begin{definition} \label{Def-QP}
For each integer $m\ge 0$, let 
$\mathcal{S}_m$ be the
$\Q$-vector space of homogeneous symmetric polynomials of degree $m$ in the variables $x,y,z$. Given a polynomial $P\in\mathcal{S}_m$, define the class
\[
\QQQ_P = \sum_{a+b+c = m+9} f_{abc}\cdot[\Theta_{abc}]\, \in\,  \mathsf{S}_{2,14+2m}^{9+m}\,,
\]
where the coefficients are determined by
\[
xyz(x-y)^2(x-z)^2(y-z)^2P = \sum_{a+b+c=m+9}f_{abc}x^ay^bz^c.
\]
\end{definition}

For $P=1$, the degree 0 constant polynomial, Definition \ref{Def-QP} specializes to 
Definition \ref{q-def}:
$$\QQQ_1 = 6\QQQ \in \mathsf{S}_{2,14}^{9}\, ,$$
up to scaling.

Let $\mathsf{GorKer}_{g,n}^*\subset \mathsf{S}_{g,n}^*$ denote the ideal of classes $\gamma\in \mathsf{S}_{g,n}^*$ for which
$\epsilon(\gamma)$ lies in the Gorenstein kernel of the tautological ring:
\begin{equation} \label{pprr}
\int_{\Mbar_{g,n}} \epsilon(\gamma) \cdot \widehat{\gamma} = 0\, , \ \ \ \forall\,  \widehat{\gamma} \in \R^*(\Mbar_{g,n})\, .
\end{equation}
We have the inclusion $\mathsf{CohRel}_{g,n}^* \subset \mathsf{GorKer}_{g,n}^*$.
With these preliminaries (and assuming Proposition~\ref{prop:nonsep_pullback}), we can now prove the main result of Section \ref{secvan}.

\begin{thm}\label{thm:gorenstein_kernel}
Let $m\geq 0$, and let $P(x,y,z) \in \mathcal{S}_m$ be a homogeneous symmetric polynomial of degree $m$. Then, $$\QQQ_{P} \in \mathsf{GorKer}_{2,14+2m}^{9+m}\, .$$
\end{thm}
\begin{proof} Let $k=7+m$, and let $\widehat{P} = xyz(x-y)^2(x-z)^2(y-z)^2P$, the polynomial used in the definition of $\QQQ_P$. For the pairing \eqref{pprr}, since $\epsilon(\QQQ_{P})$ has codimension $k+2$, 
we must consider only classes  $\widehat{\gamma}\in \R^{k+1}(\Mbar_{2,2k})$. 
By \cite{FP2,GV}, such a class $\widehat{\gamma}$ can be written in terms of tautological classes supported on strata with at least $k$ rational components. 
By Proposition~\ref{prop:nonsep_pullback} (along with the fact that $\widehat{P}(x,y,y) = 0$) and Lemma~\ref{lem:sep_pullbacks}, the pairing vanishes for all strata other than the rational tails strata in which the red/blue legs bubble off in pairs. Since $\kappa_1$ can be replaced by $\psi$ and boundary classes on a genus $2$ vertex, the only remaining class to pair with is (up to symmetry) the class considered in Lemma~\ref{lem:rational_tail_intersection_number}. Since $\widehat{P}(x,x,x) = 0$, the sum of the coefficients of $\widehat{P}$ is zero and thus this intersection number vanishes as well.
\end{proof}

Since Petersen \cite{Petersen2016-GenusTwo} showed that all even cohomology is tautological in genus $2$ when there are fewer than $20$ markings, we immediately obtain the following consequence via Poincar\'e duality.
\begin{cor}\label{cor:theta_relations}
Suppose $0 \le m \le 2$ and $P\in\mathcal{S}_m$. Then,
$$
\QQQ_{P} \in \mathsf{CohRel}_{2,14+2m}^{9+m}\, . $$
\end{cor}
Setting $P = 1$ gives Theorem~\ref{qth_rel}. Since $\mathcal{S}_1$ is $1$-dimensional and $\mathcal{S}_2$ is $2$-dimensional, Corollary \ref{cor:theta_relations} yields four relations in total.

\subsection{Proof of Proposition \ref{prop:nonsep_pullback}}\label{subsec:beads}
We now prove Proposition~\ref{prop:nonsep_pullback}, which reduces to an explicit identity in $\R^*(\Mbar_{1,2k+2})^{[1^k]\boxtimes[1^k]}$ obtained only from genus zero tautological relations.

We begin by constructing the classes $\CCC_a$. When $j^*\epsilon([\Theta_{abc}])$ is computed in terms of basic tautological classes (via the rules of \cite{GP}), there are two types of terms that occur: 
we can either cut a nonseparating edge and put $-\psi-\psi'$ on the two resulting legs, or we can split a vertex into two vertices and cut the edge between them (in such a way as to keep the graph stable and connected). In either case, exactly one of the three paths (of lengths $a,b,c$) in the $\Theta$ graph will be broken. We can therefore partition the resulting terms into three groups: let $\CCC_{a,\{b,c\}}$ be the contribution in which the edge of length $a$ is broken. Then, by construction, we have
\[
j^*\epsilon([\Theta_{abc}]) =\CCC_{a,\{b,c\}} + \CCC_{b,\{a,c\}} + \CCC_{c,\{a,b\}}\, .
\]
We will show that $\CCC_{a,\{b,c\}}$ does not depend on $b$ and $c$ (for fixed $k=a+b+c-2$), so that we can take $\CCC_a= \CCC_{a,\{b,c\}}$.
The main tool is the following result.

\begin{lem}\label{lem:cycle_identity}
Let $a$ and $b$ be positive integers, and let $k = a+b-2$. Let $\mathsf{O}_{ab}$ be the stable graph for $\Mbar_{1,2+2k}$ formed by the following three steps: 
\begin{enumerate}
\item[(i)]
take two vertices $v,w$ and connect them by two parallel paths of lengths $a,b$, 
\item[(ii)] place legs $1,2$ on the vertices $v,w$ respectively,
\item[(iii)] at each other vertex, place a pair of legs $i,i+k$ for $3 \le i \le k+2$. 
\end{enumerate}
Then, the antisymmetrized class
\[
\epsilon(\AAA_{2,k,k}([\mathsf{O}_{ab}]))\in \R^*(\Mbar_{1,2+2k})^{[1^k]\boxtimes[1^k]}
\]
only depends on $k$, and not on $a$ and $b$ individually.

Moreover, take $\mathsf{O}'_{ab}$ to be the same graph but with an extra {\em red} leg on the vertex $w$. Then, the antisymmetrized class
\[
\epsilon(\AAA_{2,k+1,k}([\mathsf{O}'_{ab}]))\in \R^*(\Mbar_{1,2+2k+1})^{[1^{k+1}]\boxtimes[1^k]}
\]
again only depends on $k$.
\end{lem}
The graphs $\mathsf{O}_{24}$ and $\mathsf{O}'_{24}$ are pictured below. Lemma \ref{lem:cycle_identity} states that we can move vertices with paired red/blue legs between the two paths freely.

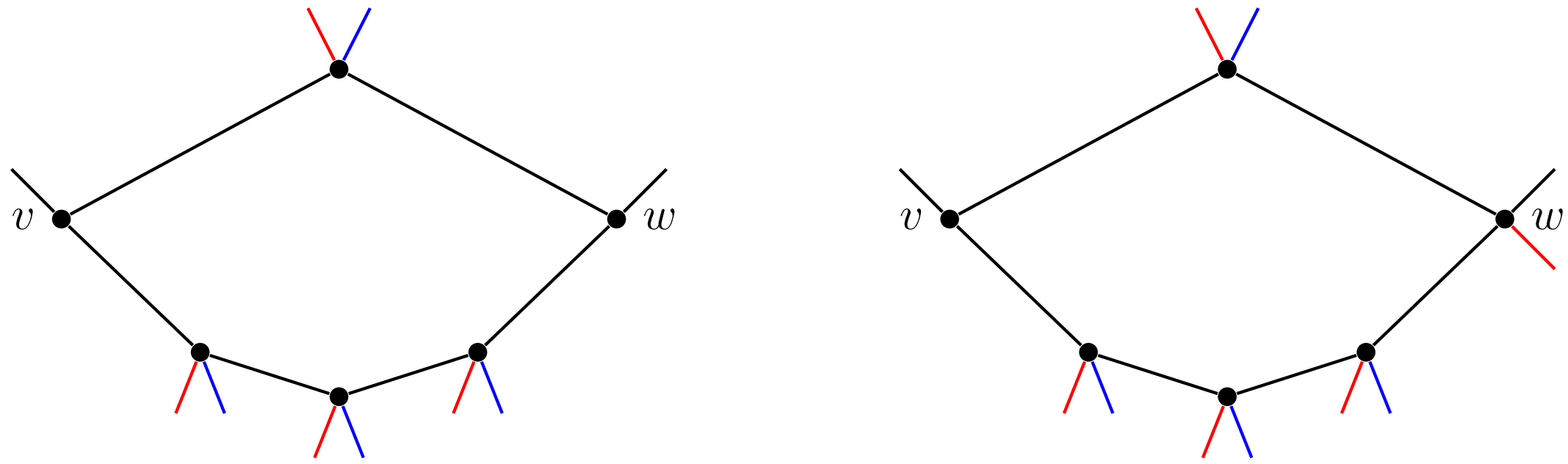
\vspace{5pt}
\begin{tikzpicture}[
    vertex/.style={circle,fill=black,inner sep=1.7pt},
    edge/.style={thick},
    leg/.style={thick},
]

\begin{scope}[xshift=0cm]

\node[vertex,label=left:$v$]  (v1) at (0,0) {};
\node[vertex,label=right:$w$] (w1) at (5,0) {};

\node[vertex] (u1) at (2.5,1.35) {};
\draw[edge] (v1) -- (u1) -- (w1);

\node[vertex] (d11) at (1.25,-1.2) {};
\node[vertex] (d12) at (2.5,-1.6) {};
\node[vertex] (d13) at (3.75,-1.2) {};
\draw[edge] (v1) -- (d11) -- (d12) -- (d13) -- (w1);

\draw[leg,red]  (u1) -- ++(-0.28,0.55);
\draw[leg,blue] (u1) -- ++( 0.28,0.55);

\foreach \x in {d11,d12,d13}{
    \draw[leg,red]  (\x) -- ++(-0.22,-0.55);
    \draw[leg,blue] (\x) -- ++( 0.22,-0.55);
}

\draw[leg] (v1) -- ++(-0.45,0.45);

\draw[leg] (w1) -- ++(0.45,0.45);

\end{scope}

\begin{scope}[xshift=8cm]

\node[vertex,label=left:$v$]  (v2) at (0,0) {};
\node[vertex,label=right:$w$] (w2) at (5,0) {};

\node[vertex] (u2) at (2.5,1.35) {};
\draw[edge] (v2) -- (u2) -- (w2);

\node[vertex] (d21) at (1.25,-1.2) {};
\node[vertex] (d22) at (2.5,-1.6) {};
\node[vertex] (d23) at (3.75,-1.2) {};
\draw[edge] (v2) -- (d21) -- (d22) -- (d23) -- (w2);

\draw[leg,red]  (u2) -- ++(-0.28,0.55);
\draw[leg,blue] (u2) -- ++( 0.28,0.55);

\foreach \x in {d21,d22,d23}{
    \draw[leg,red]  (\x) -- ++(-0.22,-0.55);
    \draw[leg,blue] (\x) -- ++( 0.22,-0.55);
}

\draw[leg] (v2) -- ++(-0.45,0.45);

\draw[leg]     (w2) -- ++(0.45, 0.45);
\draw[leg,red] (w2) -- ++(0.45,-0.45);

\end{scope}

\end{tikzpicture}
\vspace{5pt}

When proving Lemma  \ref{lem:cycle_identity} and similar identities in this paper, it will be convenient to introduce specialized notation for antisymmetrized boundary strata classes in which every vertex has genus zero. We omit all red/blue legs and instead replace each vertex by a bead of color indicating whether it has a red leg, a blue leg, both, or neither (colors red, blue, black, white respectively). Then given a bead diagram $G$ with beads of these type and $n$ labeled legs, we replace the beads with genus zero vertices with red/blue legs as indicated and choose a labeling of the red/blue legs to obtain a stable graph $\Gamma$. We take
\[
[G] = \epsilon(\AAA_{n,a,b}([\Gamma]))\, ,
\]
which is well-defined only up to sign.  To make $[G]$ truly well-defined, we must choose a convention for labeling the red/blue legs. As we have previously done, the pairs of red/blue legs on the same vertices (coming from black beads) will always have corresponding labels respecting the ordering of the markings of each color. For isolated red or blue legs, we keep track of their order by labeling them in $G$.

Lemma~\ref{lem:cycle_identity} can therefore be interpreted as saying that the classes given by the two bead diagrams 
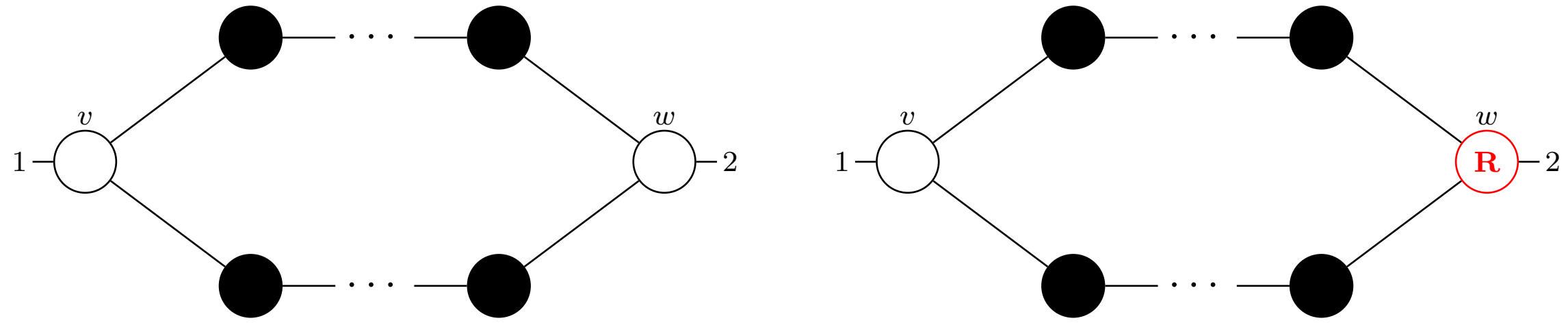
\[
  \TwoPathGraph{O}{O}{0}{0}
  \qquad
  \TwoPathGraph{O}{A}{0}{0}
\]
are unchanged by moving black beads between the two paths connecting $v$ and $w$. We will prove these invariances by direct application of the simple relations in Figure~\ref{fig:rel1}, which all follow from the WDVV relation on $\Mbar_{0,4}$ along with antisymmetrization.

\begin{figure}
\begin{tabular}{|N|l|N|l|}
  \hline
  $1$ & \BeadSum[\BeadHalfSumWidth]{
    & \BracketedBead{AB}-\BracketedBead{BA}
      } &
  $1'$ & \BeadSum[\BeadHalfSumWidth]{
    & \BracketedBead{AO_}-\BracketedBead{O_A}
      } \\
  \hline
  $2$ & \BeadSum[\BeadHalfSumWidth]{
    & \BracketedBead{AP}-\BracketedBead{PA}
      } &
  $2'$ & \BeadSum[\BeadHalfSumWidth]{
    & \BracketedBead{AA'_}-\BracketedBead{A'_A}
      } \\
  \hline
  $3$ & \multicolumn{3}{l|}{\BeadSum{
    & \BracketedBead{AB_}-\BracketedBead{B_A}+\BracketedBead{PO_}-\BracketedBead{O_P}
      }} \\
  \hline
  $4$ & \multicolumn{3}{l|}{\BeadSum{
    & \BracketedBead{AP_}-\BracketedBead{P_A}-\BracketedBead{PA_}+\BracketedBead{A_P}
      }} \\
  \hline
  $5$ & \multicolumn{3}{l|}{\BeadSum{
    & \BracketedBead{sP_}-\BracketedBead{AB_}-\BracketedBead{BA_}-\BracketedBead{PO_}
      }} \\
  \hline
\end{tabular}
\caption{Some relations in $\R^*(\Mbar_{0,n+a+b})^{[1^a]\boxtimes[1^b]}$.}
\label{fig:rel1}
\end{figure}

Relation $1$ is just exactly WDVV. Relation $1'$ is also WDVV: it is the same as Relation $1$ except that the single blue leg is uncolored (this is also the relationship between Relations $2$ and $2'$). Relation $2$ is the pullback of Relation $1$ by forgetting a red leg (since terms with two red legs on the same vertex vanish by antisymmetry). Similarly, Relation $3$ is the pullback of Relation $1'$ by forgetting a blue leg, and Relation $4$ is the pullback of Relation $3$ by forgetting a red leg (two of the four terms change sign because the red legs change order). Finally, Relation $5$ is given by pulling back the relation $\psi_1=0$ on $\Mbar_{0,3}$.

\begin{proof}[Proof of Lemma~\ref{lem:cycle_identity}]
We prove the first identity concerning $\mathsf{O}_{ab}$ first.
We begin by trying to move a single black bead (with a red leg and a blue leg) through $v$. Using Relation 3 leaves two error terms: the blue leg has moved onto $v$ and the red leg is on a vertex by itself on one of the two paths, with opposite signs. The result is that we need to show that the difference below vanishes.
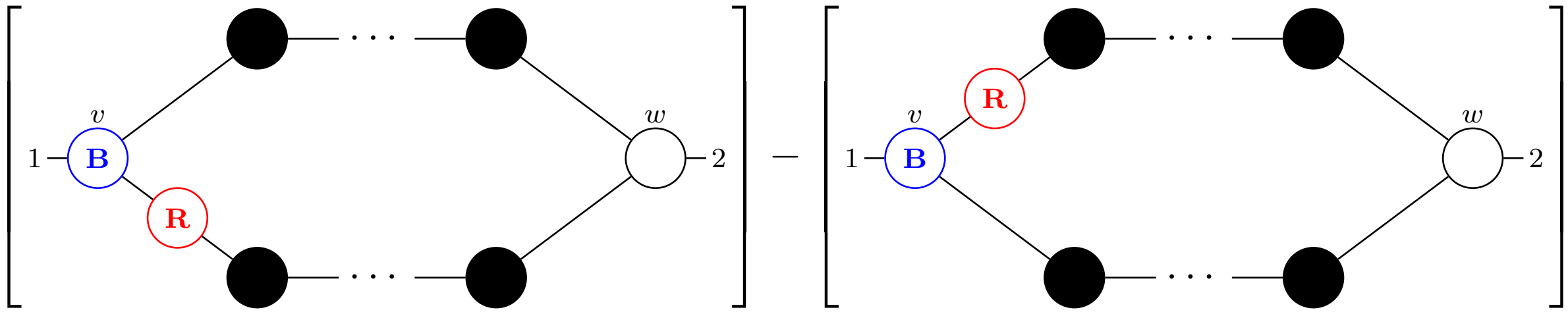
\[
  \left[\TwoPathGraph{B}{O}{0}{1}\right]
  -
  \left[\TwoPathGraph{B}{O}{1}{0}\right]
\]
We then take the red bead in one of the two terms and move it around the cycle the long way, through $w$. This can be done by using Relations $2$ and $1'$, so the two terms cancel.

The second identity (concerning $\mathsf{O}'_{ab}$) follows by an identical argument except that we use Relation $2'$ instead of $1'$ when moving the red bead past vertex $w$.
\end{proof}

\begin{figure}[htbp]
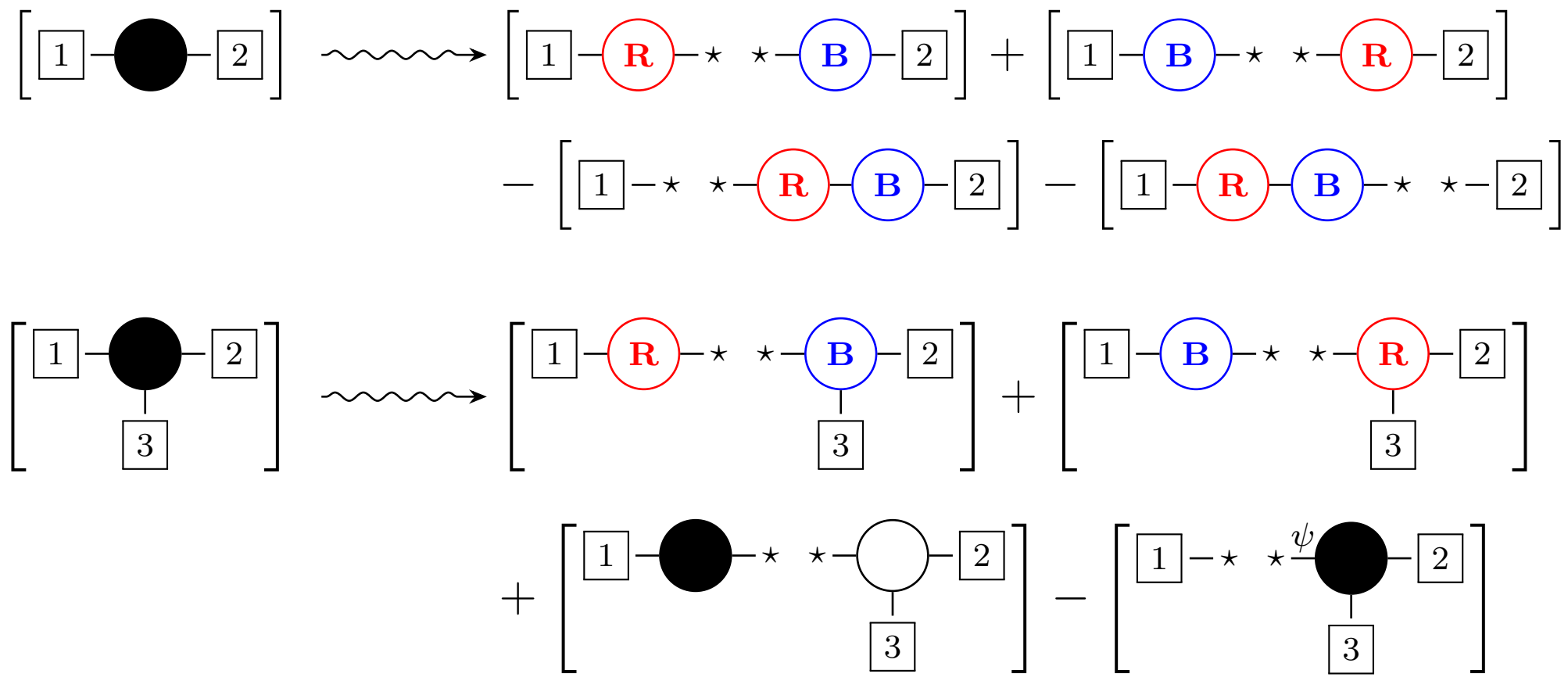

  \centering
  \(\begin{aligned}
    \BeadInput{P}\;\BeadRewriteArrow
      &\BeadPair{R}{B}+\BeadPair{B}{R}\\[2mm]
      &{}-\BeadPair{}{RB}-\BeadPair{RB}{}\\[6mm]
    \BeadInput{P_}\;\BeadRewriteArrow
      &\BeadPair{R}{B_}+\BeadPair{B}{R_}\\[2mm]
      &{}+\BeadPair{P}{O_}-\BeadPair{}{sP_}
  \end{aligned}\)
  \caption{Local operations for breaking a bead diagram when pulling back to a boundary divisor. In the first operation, psi insertions have been replaced by boundary divisors. In the second operation, only the one-third of the terms in which the path from $\boxed{2}$ to $\boxed{3}$ is unbroken are shown.}
  \label{fig:split}
\end{figure}

We now resume proving that $\CCC_{a,\{b,c\}}$ doesn't depend on $b,c$. Recall that the subdivided theta graph we start with corresponds to the triple-path bead diagram below:
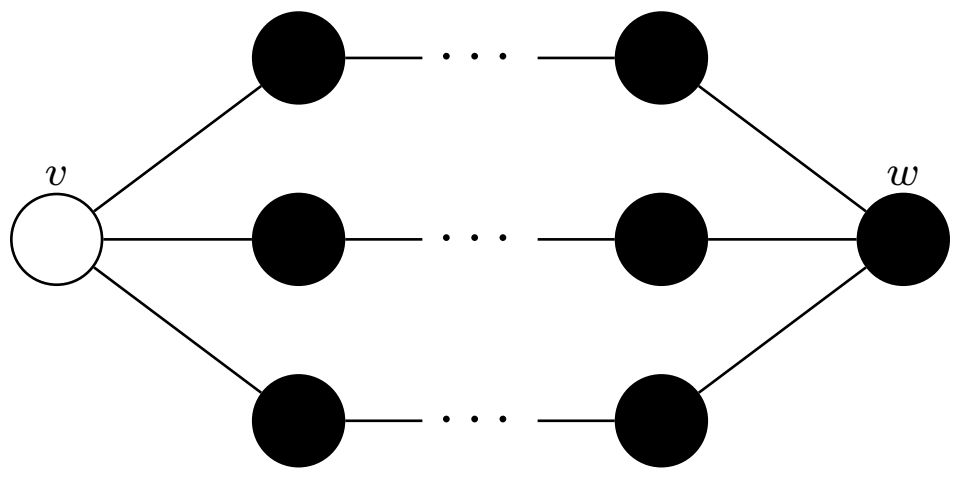
\[
  \ThreePathGraph
\]
There are three types of terms contributing to $\CCC_{a,\{b,c\}}$: we can cut an edge along the path of length $a$, we can split a vertex along the path of length $a$, or we can split vertex $w$ (breaking off the path of length $a$). Moreover, when we cut an edge, we can put $-\psi$ on either side of the cut. If we group these terms based on which vertex is split or acquires a psi insertion, then we get two types of local operations, shown in Figure~\ref{fig:split}.

We start by considering the terms coming from applying the second operation of Figure~\ref{fig:split} to vertex $w$ (so $\boxed{1}$ is the direction of the path of length $a$). The fourth term in this operation has a psi insertion, which we remove using Relation $5$ from Figure~\ref{fig:rel1}. This leaves six boundary strata, and they all have the property that the junction vertex $w$ has lost at least one colored leg. This means that Lemma~\ref{lem:cycle_identity} applies to each of them and we see that $b,c$ do not matter, as desired.

The remaining terms come from applying the first operation of Figure~\ref{fig:split} to one of the internal vertices of the path of length $a$. The idea is the same: we want to use relations to remove at least one colored leg from $w$ so that we can apply Lemma~\ref{lem:cycle_identity}.

We do this by grouping the four terms into pairs that only differ in sign and in location of the red bead. We can then use the same bead relations as in the proof of Lemma~\ref{lem:cycle_identity} to slide this red bead past black beads (and possibly through $v$) so that it is adjacent to $w$, either along a broken half of the path of length $a$ or along one of the other two paths. We then apply Relation 4 (from Figure~\ref{fig:rel1}) to replace these two terms with ones where we can apply Lemma~\ref{lem:cycle_identity}.
\qed

\vspace{5pt}
\noindent
The proofs of Theorem \ref{thm:gorenstein_kernel}, Corollary \ref{cor:theta_relations}, and Theorem \ref{qth_rel} are also now complete.

\section{Non-vanishing on the boundary}\label{sec:boundarynonvanishing}

\subsection{Boundary analysis} \label{boundan}
We will consider the following maps of moduli spaces
$$\Mbar_{1,14\cup\{ \bullet,\star\}}/\SSS_2\,  \stackrel{\epsilon}{\longrightarrow} \, \Delta_0 \, \stackrel{\iota_0}{\longrightarrow}\, 
\Delta\, \stackrel{\iota_\Delta}{\longrightarrow}\,  \Mbar_{2,14}\, ,$$
where $\epsilon$ is the normalization of the divisor 
$\Delta_0\subset \Mbar_{2,14}$ of curves with a non-separating node and $\Delta$ is the full boundary $\Mbar_{2,14}\smallsetminus \M_{2,14}$. 
The $\SSS_2$-action permutes the  markings $\{ \bullet,\star\}$ which are
glued to define $\epsilon$.
Both $\iota_0$ and $\iota_\Delta$ are closed 
inclusions of substacks.

Since the strata which appear in Definition \ref{q-def} of
$\QQQ\in \mathsf{S}_{2,14}^9$ 
include non-separating nodes,
$\QQQ$ canonically determines
homology classes
$$\QQQ_{\Delta_0} \in \HH_{16}(\Delta_0)\ \ \ 
\text{and} \ \ \ \QQQ_\Delta \in \HH_{16}(\Delta)$$
which satisfy push-forward equations:
\begin{equation}\label{pfe}
\iota_{0*}(\QQQ_{\Delta_0})=\QQQ_\Delta\, , \ \ \ \iota_{\Delta*}(\QQQ_\Delta)=0\, .
\end{equation}
The second equality in \eqref{pfe}
is a consequence of
Theorem \ref{qth_rel} and Poincar\'e duality.  We have switched to
homology for $\QQQ_{\Delta_0}$ and $\QQQ_\Delta$, since for singular spaces the cycle class map naturally takes values in Borel--Moore homology rather than cohomology.

\begin{prop} We have $\QQQ_{\Delta_0} \neq 0\in \HH_{16}(\Delta_0)$. \label{propdelq}
\end{prop}

\begin{proof}
The homology of $\Delta_0$ is obtained from
the homology of 
$\Mbar_{1,14\cup\{ \bullet,\star\}}/\SSS_2$
via standard Mayer--Vietoris sequences
obtained from the strata of the normalization morphism $\epsilon$.
More precisely, the kernel
of 
\begin{equation} \label{epspush}
\epsilon_*: \HH_{16}(\Mbar_{1,14\cup\{ \bullet,\star\}}/\SSS_2)\,  {\longrightarrow} \, \HH_{16}(\Delta_0)
\end{equation}
receives a surjection
\begin{equation}\label{mvs}
\left(\bigoplus_\Gamma \HH_{16}(\Mbar_\Gamma)\right) \, \rightarrow\,  \mathsf{Ker}(\epsilon_*)\,  \rightarrow 0\, ,
\end{equation}
where the sum on the left is over all
codimension two strata 
of $\Mbar_{2,14}$
corresponding to dual graphs in which both edges are nonseparating, and the map is given by the difference $f_*-g_*$ of the pushforwards by the two maps
\[
f_\Gamma,g_\Gamma:\Mbar_\Gamma \to \Mbar_{1,14\cup\{ \bullet,\star\}}/\SSS_2
\]
given by gluing along just one of the edges.

Let $(\mathsf{S}_{1,16}^*)^{\SSS_2}$
denote the strata algebra{\footnote{The superscript here denotes the $\SSS_2$-invariant part.}} of 
$\Mbar_{1,14\cup\{ \bullet,\star\}}$, and 
let $\widehat{\QQQ}\in 
(\mathsf{S}_{1,16}^8)^{\SSS_2}
$  be a lift{\footnote{$\widehat{\QQQ}$ is not canonical, there are many choices.}} of
$\QQQ\in 
\mathsf{S}_{2,14}^9
$
with respect to the gluing map
$$(\mathsf{S}_{1,16}^8)^{\SSS_2} \rightarrow \mathsf{S}_{2,14}^9\, .$$
We denote the corresponding homology class also by
$\widehat{\QQQ}\in \HH_{16}(
\Mbar_{1,14\cup\{ \bullet,\star\}}/\SSS_2
)$. Then,
$$\epsilon_*(\widehat{\QQQ}) = \QQQ_{\Delta_0}\, .$$
In order to prove $\QQQ_{\Delta_0}\neq 0 \in \HH_{16}(\Delta_0)$, we must
show
$\widehat{\QQQ} \notin \mathsf{Ker}(\epsilon_*)$.

Following Definition \ref{Def-QP}, let $P\in \mathcal{S}_m$ be nonzero,
and let
$\QQQ_{P}\in \mathsf{S}_{2,14+2m}^{9+m}$.
We will prove more generally that $\widehat{\QQQ}_P\notin\mathsf{Ker}(\epsilon_*)$, where $$\widehat{\QQQ}_P\in
H_{16+2m}(\Mbar_{1,14+2m\cup\{ \bullet,\star\}}/\SSS_2)
$$ is an analogous lift  with respect to
\begin{equation} \label{epspush-m}
\epsilon_*: \HH_{16+2m}(\Mbar_{1,14+2m\cup\{ \bullet,\star\}}/\SSS_2)\,  {\longrightarrow} \, \HH_{16+2m}(\Delta_0)
\end{equation}
of the class $(\QQQ_{P})_{\Delta_0}\in
\HH_{16+2m}(\Delta_0)$ associated to $\QQQ_P$.

The strategy is to find a detector class $\mathsf{D}_P\in \HH^{16+2m}(\Mbar_{1,14+2m\cup\{ \bullet,\star\}}/\SSS_2)$ with the following two properties:
\begin{itemize}
\item[(i)] For every stratum $\Gamma$ that appears in
\[
f_\Gamma,g_\Gamma:\Mbar_\Gamma \to \Mbar_{1,14+2m\cup\{ \bullet,\star\}}/\SSS_2\, ,
\]
we have
$
(f_\Gamma^* - g_\Gamma^*)\mathsf{D}_P = 0$.

\item[(ii)] The intersection number $\mathsf{D}_P\cdot\widehat{\QQQ}_P$ is nonzero.
\end{itemize}
By property (i),  property (ii)
does not depend on the choice of $\widehat{\QQQ}_P$ by the projection formula. Together, the properties show that $\widehat{\QQQ}_P\notin \mathsf{Ker}(\epsilon_*)$.
As a consequence, we obtain the following result which generalizes Proposition \ref{propdelq}.
\begin{prop} We have $(\QQQ_{P})_{\Delta_0} \neq 0\in \HH_{16+2m}(\Delta_0)$ for all nonzero $P\in\mathcal{S}_m$. \label{propdel}
\end{prop}
\noindent We will construct $\mathsf{D}_P$ and check properties (i) and (ii) in  Sections \ref{con-dp} and \ref{proof-DP}.
\end{proof}

\subsection{Construction of \texorpdfstring{$\mathsf{D}_P$}{DP}} \label{con-dp}
From now on we fix a nonzero $P\in\mathcal{S}_m$ and let $k=7+m$. Like $\QQQ_P$, $\mathsf{D}_P$ will be antisymmetric with respect to the action of $S_{k}\times S_{k}$ permuting red and blue markings. We use the same bead diagram notation as in Section~\ref{subsec:beads} to describe $\mathsf{D}_P$ and check the desired boundary pullback properties. These properties will be true in Chow, not just cohomology, and they will follow from the WDVV relation (via bead relations like the ones in Figure~\ref{fig:rel1}).

The class $\mathsf{D}_P$ will be given as a linear combination of two different types of antisymmetrized boundary strata. Consider bead diagrams $X(a,b;i_1,i_2)$ and $Y(a;i_1,i_2)$ as depicted below, with $a+b+i_1+i_2$ and $a+i_1+i_2$ edges respectively (the paths are labeled with the number of edges, not the number of black beads).  Although we've labeled legs $1$ and $2$ separately in the diagram for $X(a,b;i_1,i_2)$, in this section we will always symmetrize with respect to the $\Sigma_2$-action interchanging those legs. We assume that $a,b > 0$, but $i_1,i_2\ge 0$.

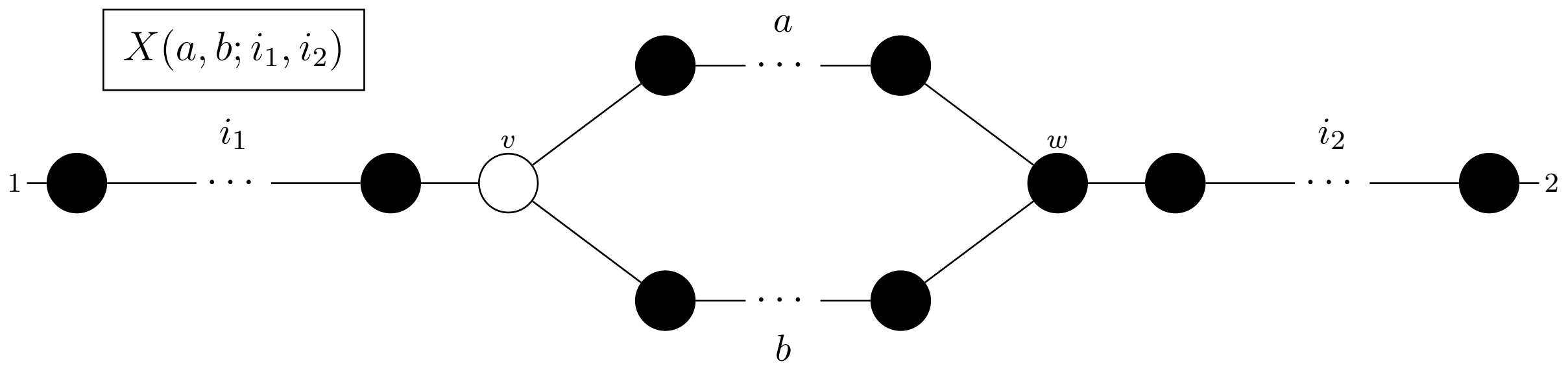
\[
  \begin{tikzpicture}[
    beads/defaults,baseline=-.5ex,
    graphs/length/.style={
      beads/label,text=black,font=\small
    }
  ]
    \GraphBead{O}{v}{0,0}
    \GraphBead{P}{w}{7,0}
    \GraphPath{upper}{1.5}{0}
    \GraphPath{lower}{-1.5}{0}
    \GraphVertexNames

    \GraphBead{P}{left-end}{-5.5,0}
    \GraphBead{P}{left-inner}{-1.5,0}
    \node[graphs/dots] (left-dots) at (-3.5,0) {$\cdots$};
    \draw[beads/edge]
      (left-end) -- (left-dots.west)
      (left-dots.east) -- (left-inner)
      (left-inner) -- (v);
    \draw[beads/edge] (left-end.west)
      -- ++(180:\pgfkeysvalueof{/tikz/beads/leg length})
      node[beads/label,anchor=east] {1};

    \GraphBead{P}{right-inner}{8.5,0}
    \GraphBead{P}{right-end}{12.5,0}
    \node[graphs/dots] (right-dots) at (10.5,0) {$\cdots$};
    \draw[beads/edge]
      (w) -- (right-inner)
      (right-inner) -- (right-dots.west)
      (right-dots.east) -- (right-end);
    \draw[beads/edge] (right-end.east)
      -- ++(0:\pgfkeysvalueof{/tikz/beads/leg length})
      node[beads/label,anchor=west] {2};

    \node[graphs/length,anchor=south,yshift=2mm]
      at (upper-dots.north) {$a$};
    \node[graphs/length,anchor=north,yshift=-2mm]
      at (lower-dots.south) {$b$};
    \node[graphs/length,anchor=south,yshift=2mm]
      at (left-dots.north) {$i_1$};
    \node[graphs/length,anchor=south,yshift=2mm]
      at (right-dots.north) {$i_2$};

    \node[rectangle,draw=black,inner sep=2mm,font=\normalsize]
      at (-3.5,1.7) {$X(a,b;i_1,i_2)$};
  \end{tikzpicture}
\]

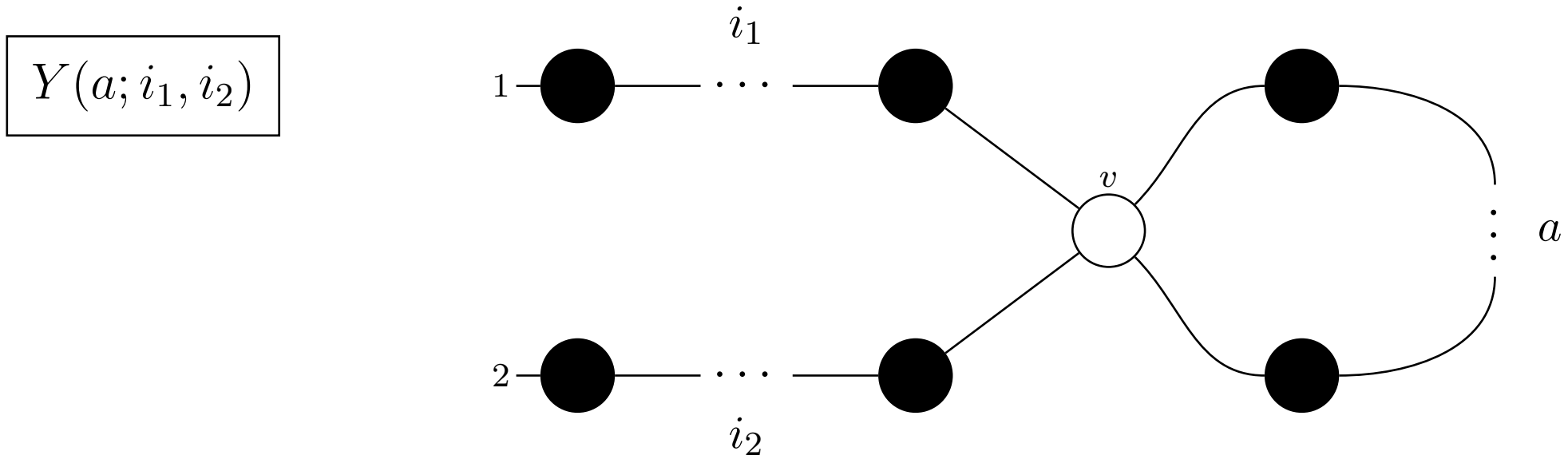
\[
  \begin{tikzpicture}[
    beads/defaults,baseline=-.5ex,
    graphs/length/.style={
      beads/label,text=black,font=\small
    }
  ]
    \GraphBead{O}{v}{0,0}
    \GraphBead{P}{loop-upper}{2,1.5}
    \GraphBead{P}{loop-lower}{2,-1.5}
    \node[graphs/dots,rotate=90] (loop-dots) at (4,0) {$\cdots$};

    \draw[beads/edge]
      (v) to[out=45,in=180] (loop-upper)
      (loop-upper) to[out=0,in=90] (loop-dots.east)
      (loop-dots.west) to[out=-90,in=0] (loop-lower)
      (loop-lower) to[out=180,in=-45] (v);
    \node[beads/label,anchor=south] at (v.north) {$v$};

    \GraphBead{P}{leg-one-end}{-5.5,1.5}
    \GraphBead{P}{leg-one-inner}{-2,1.5}
    \node[graphs/dots] (leg-one-dots) at (-3.75,1.5) {$\cdots$};
    \draw[beads/edge]
      (leg-one-end) -- (leg-one-dots.west)
      (leg-one-dots.east) -- (leg-one-inner)
      (leg-one-inner) -- (v);
    \draw[beads/edge] (leg-one-end.west)
      -- ++(180:\pgfkeysvalueof{/tikz/beads/leg length})
      node[beads/label,anchor=east] {1};

    \GraphBead{P}{leg-two-end}{-5.5,-1.5}
    \GraphBead{P}{leg-two-inner}{-2,-1.5}
    \node[graphs/dots] (leg-two-dots) at (-3.75,-1.5) {$\cdots$};
    \draw[beads/edge]
      (leg-two-end) -- (leg-two-dots.west)
      (leg-two-dots.east) -- (leg-two-inner)
      (leg-two-inner) -- (v);
    \draw[beads/edge] (leg-two-end.west)
      -- ++(180:\pgfkeysvalueof{/tikz/beads/leg length})
      node[beads/label,anchor=east] {2};

    \node[graphs/length,anchor=south,yshift=2mm]
      at (leg-one-dots.north) {$i_1$};
    \node[graphs/length,anchor=north,yshift=-2mm]
      at (leg-two-dots.south) {$i_2$};
    \node[graphs/length,anchor=west,xshift=3mm] at (4,0) {$a$};

    \node[rectangle,draw=black,inner sep=2mm,font=\normalsize]
      at (-10,1.5) {$Y(a;i_1,i_2)$};
  \end{tikzpicture}
\]

We sum over different lengths for the two tails to define the classes
\[
\mathsf{X}_{a,b;i} = \sum_{i_1 = 0}^i[X(a,b;i_1,i-i_1)],\quad \mathsf{Y}_{a;i} = \sum_{i_1 = 0}^i[Y(a;i_1,i-i_1)],
\]
both belonging to $\R^{k+1}(\Mbar_{1,2+2k})^{[1^k]\boxtimes[1^k]}$ if the sum of the subscripts is $k+1$.

Before defining $\mathsf{D}_P$, it will be convenient to define two special linear combinations of the $\mathsf{X}_{a,b;i}$ and $\mathsf{Y}_{a;i}$:
\[
\mathsf{Z}_{a;i} := \mathsf{Y}_{a;i} + \sum_{a_1=1}^{a-1}\mathsf{X}_{a_1,a-a_1;i},\quad \mathsf{Z}'_{a;i} := \mathsf{Y}_{a;i}+2\sum_{b=1}^{i}\mathsf{X}_{a,b;i-b}.
\]

\begin{definition}
    Given $P\in\mathcal{S}_m$, let (with $k=m+7$)
    \[
    \mathsf{E}_P = \sum_{a+b+i=k+1}g_{a,b;i}\mathsf{X}_{a,b;i},
    \]
    where the coefficients are determined by
\[
x^2y^2(x-z)^2(y-z)^2P = \sum_{a+b+i=k+1}g_{a,b;i}x^ay^bz^i.
\]
The class $\mathsf{D}_P$ is then defined as the unique linear combination
\[
\mathsf{D}_P = \mathsf{E}_P + \sum_{a=4}^{k+1}c_a\mathsf{Z}_{a;k+1-a}
\]
that is also (formally) a linear combination of classes of the two forms
\[
\mathsf{X}_{a,b;i}-\mathsf{X}_{a+1,b-1;i},\quad \mathsf{Z}'_{a;i}.
\]
(Such a linear combination exists and is unique because the classes $\mathsf{Z}_{a;i}-\mathsf{Z}'_{a;i}$ form a (triangular) basis for the formal span of the $\mathsf{X}_{a,b;i}$ after quotienting out by the classes $\mathsf{X}_{a,b;i}-\mathsf{X}_{a+1,b-1;i}$.)
\end{definition}

It then suffices to check several properties of these classes, which we record here before proving them in the next section. In the statement of the following proposition, let $d^*_{\text{sep}}$ and $d^*_{\text{nonsep}}$ denote (any of) the pullback difference maps $f_\Gamma^*-g_\Gamma^*$ from before with two-vertex and one-vertex $\Gamma$ respectively. In other words, in the former case we pull back along separating gluing maps
\[
f_\Gamma,g_\Gamma:\Mbar_{1,2+r}\times\Mbar_{0,2+r'}\to\Mbar_{1,2+2k}
\]
where $r+r' = 2k$ is some partition of the colored legs,
and in the latter case we pull back along nonseparating gluing maps
\[
f_\Gamma,g_\Gamma:\Mbar_{0,4+2k}\to\Mbar_{1,2+2k}.
\]
In both cases, the maps $f_\Gamma,g_\Gamma$ differ only in the choice of the markings that are glued together. For example, in the latter case we either glue markings $1,2$ or markings $3,4$.
\begin{prop}\label{prop:DP_facts}
    The following are true:
    \begin{enumerate}[(a)]
    \item $d^*_{\text{sep}}\mathsf{Z}'_{a;i} = 0$;
    \item $d^*_{\text{sep}}(\mathsf{X}_{a,b;i}-\mathsf{X}_{a+1,b-1;i}) = 0$;
    \item $f_{\text{nonsep}}^*\mathsf{Z}_{a;i} = 0$ for $f_{\text{nonsep}}:\Mbar_{0,4+2k}\to\Mbar_{1,2+2k}$ (and thus $d^*_{\text{nonsep}}\mathsf{Z}_{a;i} = 0$ as well);
        \item $d^*_{\text{nonsep}}\mathsf{E}_P = 0$ (for any $P$);

        \item $\mathsf{E}_P\cdot \widehat{\QQQ}_P \ne 0$ for $P\ne 0$, where $\widehat{\QQQ}_P$ is given by breaking an edge in each boundary stratum in the formula for $\QQQ_P$.
    \end{enumerate}
\end{prop}

Once we have verified Proposition~\ref{prop:DP_facts} in the next section, we will have completed the proof of Proposition~\ref{propdel}. To see this, note that parts (c) and (d) together imply that $d^*_{\text{nonsep}}\mathsf{D}_P = 0$, parts (a) and (b) together imply that $d^*_{\text{sep}}\mathsf{D}_P = 0$, and parts (c) and (e) together imply that $\mathsf{D}_P\cdot\widehat{\mathsf{Q}}_P\ne 0$ because the boundary strata in $\widehat{\mathsf{Q}}_P$ still contain non-separating nodes.

\subsection{Proof of Proposition \ref{prop:DP_facts}}
\label{proof-DP}
\begin{enumerate}[(a)]
    \item We observe that $\mathsf{Z}'_{a;i}$ can be obtained from a simpler class $\mathsf{W}'_i \in \R^*(\Mbar_{0,3+2i})^{[1^i]\boxtimes[1^i]}$ as follows. First, define $\mathsf{W}'_i$ as the sum of bead diagram classes coming from the diagram below, where $i_1,i_2\ge 0$ run over all options with $i_1+i_2 = i$.

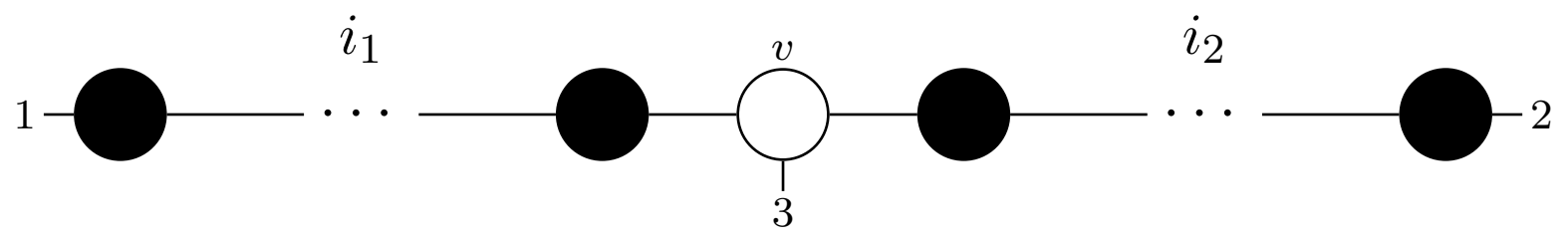
    \[
  \begin{tikzpicture}[beads/defaults,moregraphs/center]
    \GraphBead{O}{v}{0,0}
    \GraphBead{P}{left-end}{-5.5,0}
    \GraphBead{P}{left-inner}{-1.5,0}
    \GraphBead{P}{right-inner}{1.5,0}
    \GraphBead{P}{right-end}{5.5,0}
    \node[graphs/dots] (left-dots) at (-3.5,0) {$\cdots$};
    \node[graphs/dots] (right-dots) at (3.5,0) {$\cdots$};

    \draw[beads/edge]
      (left-end) -- (left-dots.west)
      (left-dots.east) -- (left-inner)
      (left-inner) -- (v)
      (v) -- (right-inner)
      (right-inner) -- (right-dots.west)
      (right-dots.east) -- (right-end);
    \MoreGraphLeg{left-end}{180}{east}{1}
    \MoreGraphLeg{right-end}{0}{west}{2}
    \MoreGraphLeg{v}{270}{north}{3}
    \node[beads/label,anchor=south] at (v.north) {$v$};
    \node[moregraphs/length,anchor=south,yshift=2mm]
      at (left-dots.north) {$i_1$};
    \node[moregraphs/length,anchor=south,yshift=2mm]
      at (right-dots.north) {$i_2$};
  \end{tikzpicture}
\]

    Then take the forgetful pullback $\pi^*\mathsf{W}'_i$ (forgetting an uncolored marking labeled $4$) and glue in a path of $a-1$ black beads between markings $3$ and $4$. The resulting class is (after re-antisymmetrizing) exactly $\mathsf{Z}'_{a;i}$.

    These operations commute with the pullback maps $f_{\Gamma}^*, g_{\Gamma}^*$ whose difference appears in $d^*_{\text{sep}}$. More precisely, consider the following commutative diagram:
    
  \[
\begin{tikzcd}
  \Mbar_{0,3+r}\times\Mbar_{0,2+r'} \arrow[r, "f"]
  & \Mbar_{0,3+2i} \\
  \Mbar_{0,4+r}\times\Mbar_{0,2+r'} \arrow[u, "\pi"] \arrow[r, "f"]
  & \Mbar_{0,4+2i} \arrow[u, "\pi"] \\
  \Mbar_{0,4+r}\times\Mbar_{0,2a}\times\Mbar_{0,2+r'}
    \arrow[u, "p_{13}"] \arrow[r, "f"] \arrow[d, "j"]
  & \Mbar_{0,4+2i}\times\Mbar_{0,2a}
    \arrow[u, "p_1"] \arrow[d, "j"] \\
  \Mbar_{1,2+r+2(a-1)}\times\Mbar_{0,2+r'} \arrow[r, "f"]
  & \Mbar_{1,2+2i+2(a-1)}
\end{tikzcd}
\]

  Here the horizontal maps $f$ glue in $\Mbar_{0,2+r'}$ along a separating node, the maps $j$ glue in $\Mbar_{0,2a}$ along two nonseparating nodes, $\pi$ forgets one of the markings glued by $j$, and $p_1,p_{13}$ project onto the indicated factors of a product. If we let $\mathsf{P}_a\in\R^*(\Mbar_{0,2a})$ be the class given by a path of $a-1$ black beads, then we have along the right column that (up to antisymmetrizing)
  \[
  \mathsf{Z}'_{a;i} = j_*((p_1^*\pi^*\mathsf{W}'_i)(p_2^*\mathsf{P}_a)).
  \]
  Since the bottom square is Cartesian, we can commute $f^*$ past $j_*$ and conclude that
  \[
  f^*\mathsf{Z}'_{a;i} = j_*((p_{13}^*\pi^*f^*\mathsf{W}'_i)(p_2^*\mathsf{P}_a)).
  \]
  In particular, it suffices to check that $(f^*-g^*)\mathsf{W}'_i = 0$, where
    \[
f,g:\Mbar_{0,3+r}\times\Mbar_{0,2+r'}\to\Mbar_{0,3+2i}
    \]
    only differ in which markings are glued (e.g.\ the first marking in each curve is glued in $f$, while the second marking in each curve is glued in $g$).

    These pullbacks are computed by breaking the graph via the first operation of Figure~\ref{fig:split}. This breaks off a path from one of the two sides of the central vertex $v$, depending on which of $f$ and $g$ is used. Since we are summing over all $i_1,i_2$, the resulting terms come in matching pairs except for the positioning of the red and blue beads created by the first operation of Figure~\ref{fig:split}. Since these colored beads can be freely commuted from one side of the graph to the other by relations 1,1',2 of Figure~\ref{fig:rel1}, the terms cancel.
    
    \item This is nearly identical to the proof of Proposition~\ref{prop:nonsep_pullback}. As in that proof, the pullbacks modify the bead diagrams $X(a,b;i_1,i_2)$ 
    by the operations in Figure~\ref{fig:split}, applied to each vertex along the two tails of lengths $i_1,i_2$. When computing a difference $(f^*-g^*)\mathsf{X}_{a,b;i}$, the only terms that don't cancel are those where at least one red or blue bead is placed on the genus $1$ component, either on the path touching $v$ or the path touching $w$, with alternating signs. But then the beads can be moved around (again using relations 1,1',2) to touch $w$ from different sides and relation 4 then converts $w$ into a red or blue bead. Then Lemma~\ref{lem:cycle_identity} applies to the resulting terms and black beads can be freely moved between the two paths in the cycle, as desired.
    
    \item This is similar to part (a): we note that $\mathsf{Z}_{a;i}$ can be obtained from a simpler class via forgetful pullback and gluing pushforward that commute with $f_{\text{nonsep}}^*$. Define $\mathsf{W}_a\in \R^*(\Mbar_{1,1+2a-2})^{[1^{a-1}]\boxtimes[1^{a-1}]}$ by taking a cycle of length $a$ with one white bead and $a-1$ black beads and attaching the remaining leg to the white bead, as in the diagram below.

    
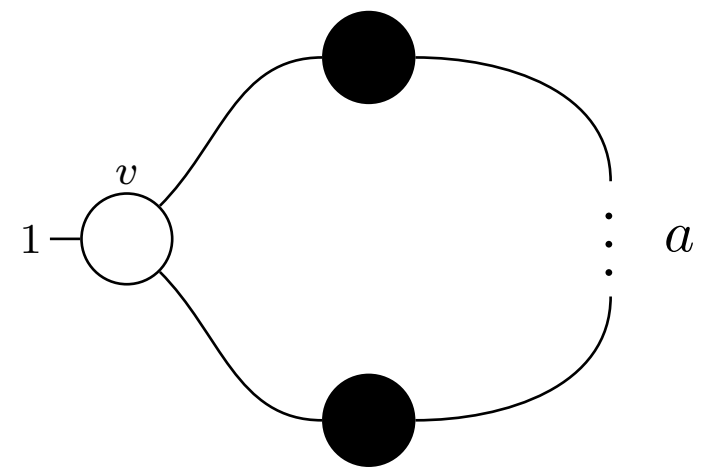
\[
  \begin{tikzpicture}[beads/defaults,moregraphs/center]
    \GraphBead{O}{v}{0,0}
    \GraphBead{P}{loop-upper}{2,1.5}
    \GraphBead{P}{loop-lower}{2,-1.5}
    \node[graphs/dots,rotate=90] (loop-dots) at (4,0) {$\cdots$};
    \draw[beads/edge]
      (v) to[out=45,in=180] (loop-upper)
      (loop-upper) to[out=0,in=90] (loop-dots.east)
      (loop-dots.west) to[out=-90,in=0] (loop-lower)
      (loop-lower) to[out=180,in=-45] (v);
    \MoreGraphLeg{v}{180}{east}{1}
    \node[beads/label,anchor=south] at (v.north) {$v$};
    \node[moregraphs/length,anchor=west,xshift=3mm]
      at (4,0) {$a$};
  \end{tikzpicture}
\]

    Then pulling back by a forgetful map (forgetting an uncolored marking labeled $2$) and gluing a path of black beads to each marking (summing over lengths $i_1,i_2$ with sum $i$) gives precisely $\mathsf{Z}_{a;i}$. These gluings across separating nodes commute with the nonseparating gluing pullback $j^*$ (similarly to part (a)), so it suffices to check that $f_{\text{nonsep}}^*\mathsf{W}_a = 0$. This pullback can be computed by the first operation in Figure~\ref{fig:split}. Since the red and blue beads can be moved all the way around the cycle in the usual way, these terms cancel.
    
    \item Let $f_{\text{nonsep}}:\Mbar_{0,4+2k}\to\Mbar_{1,2+2k}$ glue markings $1,2$. Then we need to show that $f_{\text{nonsep}}^*\mathsf{E}_P$ is invariant under the $\Sigma_2$-action simultaneously relabeling markings $1\leftrightarrow 3,2\leftrightarrow 4$. For the next part it will be convenient to compute $f_{\text{nonsep}}^*\mathsf{E}_P$ in full, so we do that here after introducing notation for the answer.

    Consider the bead diagram $T(a;i_1,i_2,j_1,j_2)$ as depicted below. Similarly to the diagrams $X(a,b;i_1,i_2)$ and $Y(a;i_1,i_2)$, we will symmetrize with respect to the $\Sigma_2$-actions interchanging legs $1,2$ or legs $3,4$.
    
    
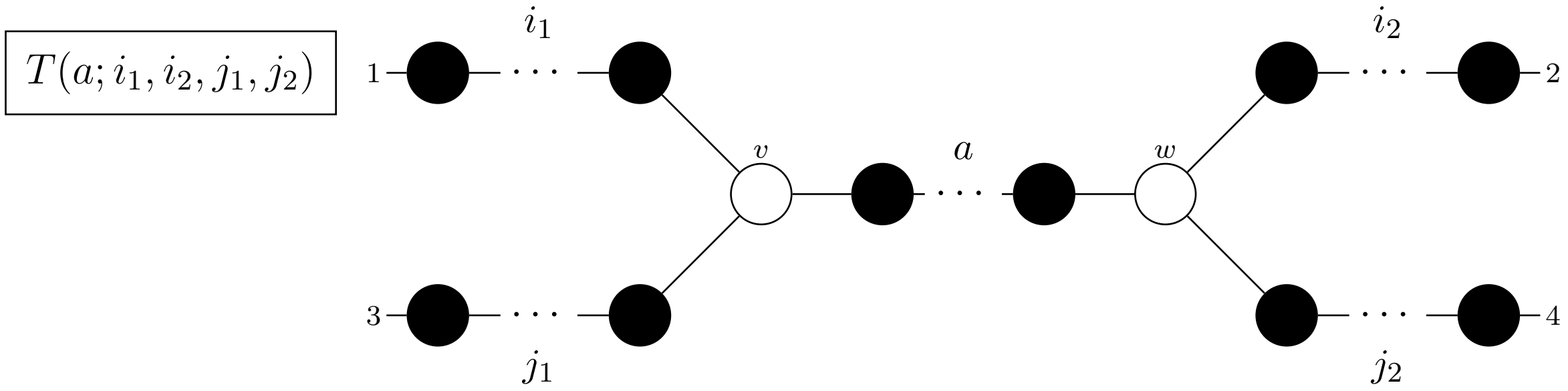
\[
  \begin{tikzpicture}[beads/defaults,moregraphs/center]
    \GraphBead{O}{v}{0,0}
    \GraphBead{O}{w}{5,0}

    \GraphBead{P}{middle-left}{1.5,0}
    \GraphBead{P}{middle-right}{3.5,0}
    \node[graphs/dots] (middle-dots) at (2.5,0) {$\cdots$};
    \draw[beads/edge]
      (v) -- (middle-left)
      (middle-left) -- (middle-dots.west)
      (middle-dots.east) -- (middle-right)
      (middle-right) -- (w);

    \foreach \side/\height/\leglabel in {upper/1.5/1,lower/-1.5/3}{
      \GraphBead{P}{left-\side-end}{-4,\height}
      \GraphBead{P}{left-\side-inner}{-1.5,\height}
      \node[graphs/dots] (left-\side-dots) at (-2.75,\height) {$\cdots$};
      \draw[beads/edge]
        (left-\side-end) -- (left-\side-dots.west)
        (left-\side-dots.east) -- (left-\side-inner)
        (left-\side-inner) -- (v);
      \MoreGraphLeg{left-\side-end}{180}{east}{\leglabel}
    }

    \foreach \side/\height/\leglabel in {upper/1.5/2,lower/-1.5/4}{
      \GraphBead{P}{right-\side-inner}{6.5,\height}
      \GraphBead{P}{right-\side-end}{9,\height}
      \node[graphs/dots] (right-\side-dots) at (7.75,\height) {$\cdots$};
      \draw[beads/edge]
        (w) -- (right-\side-inner)
        (right-\side-inner) -- (right-\side-dots.west)
        (right-\side-dots.east) -- (right-\side-end);
      \MoreGraphLeg{right-\side-end}{0}{west}{\leglabel}
    }

    \node[beads/label,anchor=south] at (v.north) {$v$};
    \node[beads/label,anchor=south] at (w.north) {$w$};
    \node[moregraphs/length,anchor=south,yshift=2mm]
      at (middle-dots.north) {$a$};
    \node[moregraphs/length,anchor=south,yshift=2mm]
      at (left-upper-dots.north) {$i_1$};
    \node[moregraphs/length,anchor=north,yshift=-2mm]
      at (left-lower-dots.south) {$j_1$};
    \node[moregraphs/length,anchor=south,yshift=2mm]
      at (right-upper-dots.north) {$i_2$};
    \node[moregraphs/length,anchor=north,yshift=-2mm]
      at (right-lower-dots.south) {$j_2$};

    \node[rectangle,draw=black,inner sep=2mm,font=\normalsize,anchor=east]
      at (-5.25,1.5) {$T(a;i_1,i_2,j_1,j_2)$};
  \end{tikzpicture}
\]
    
    We sum over lengths for the tails to define the classes
    \[
    \mathsf{T}_{a;i,j} = \sum_{i_1=0}^i\sum_{j_1=0}^j[T(a;i_1,i-i_1,j_1,j-j_1)].
    \]

     We now compute $f_{\text{nonsep}}^*\mathsf{X}_{a,b;i}$. As usual, we need to apply the two operations of Figure~\ref{fig:split}. As in Proposition~\ref{prop:nonsep_pullback}, it is convenient to group the terms based on which path between vertices $v$ and $w$ is broken. We begin by considering the terms where the path of length $b$ is broken. If we apply the first operation to an internal vertex of this path, then we obtain four terms with red and blue beads distributed across the break. Alternatively, we can apply the second operation to $w$ to get four terms, one of which still has a psi insertion.

     We use the following two relations at $w$ (moving colored beads around the cycle through $v$ when necessary) to simplify these expressions by removing all colored beads and psi insertions.

\begin{table}[H]
     \begin{tabular}{|N|l|N|l|}
  \hline
  $6$ & \multicolumn{3}{l|}{\BeadSum{
    & \BracketedBead{ABP_}
      -\BracketedBead{AP_B}
      -\BracketedBead{BP_A} \\[3mm]
    & {}+\BracketedBead{P_AB} \\[3mm]
    & {}+\BracketedBead{PPO_}
      +\BracketedBead{O_PP}
      -{\bf2}\BracketedBead{PO_P}
      }} \\
  \hline
    $7$ & \multicolumn{3}{l|}{\BeadSum{
    & \BracketedBead{sP_}
    -\BracketedBead{B_A}
    -\BracketedBead{A_B}\\[3mm]
    & {}+\BracketedBead{PO_}-{\bf2}\BracketedBead{O_P}
      }} \\
  \hline
\end{tabular}
\end{table}

     Here Relation $6$ follows from applying Relation $4$ (from Figure~\ref{fig:rel1}) twice to try to commute the blue bead past the valence $3$ black bead, and then applying Relation $3$ to try to commute the red bead past the valence $3$ blue bead thus created. Relation $7$ is simply a linear combination of Relations $5$ and $3$ (including versions with colors swapped). 

     The resulting bead diagrams are all of the form $T(a';i_1,i_2,j_1,j_2)$. Combining all the terms created in this way gives
     \[
     \mathsf{T}_{a+2;b-2,i} - 2\mathsf{T}_{a+1;b-1,i} + \mathsf{T}_{a;b,i}.
     \]
     Combining with the contributions where the path of length $a$ is broken instead, we conclude that
     \[
     f_{\text{nonsep}}^*\mathsf{X}_{a,b;i} = \mathsf{T}_{a+2;b-2,i} - 2\mathsf{T}_{a+1;b-1,i} + \mathsf{T}_{a;b,i} + \mathsf{T}_{b+2;a-2,i} - 2\mathsf{T}_{b+1;a-1,i} + \mathsf{T}_{b;a,i}
     \]

     Recalling the definition of $\mathsf{E}_P$, we see that
     \begin{equation}\label{eq:Epull}
     f_{\text{nonsep}}^*\mathsf{E}_P = 2\sum_{a+b+i = k+1}h_{a;b,i}\mathsf{T}_{a;b,i},
     \end{equation}
     where the coefficients are determined by
     \[
     x^2(x-y)^2(x-z)^2(y-z)^2P = \sum_{a+b+i = k+1}h_{a;b,i}x^ay^bz^i.
     \]
     In particular, since $P$ is symmetric this means that $h_{a;b,i} = h_{a;i,b}$, so $f_{\text{nonsep}}^*\mathsf{E}_P$ has the desired $\Sigma_2$-invariance.
    
    \item Let $T'(a;i_1,i_2,j_1,j_2)$ be the same bead diagram as in the previous part, except that the bead at vertex $w$ has been changed from white to black. We can then take
    \[
    \widehat{\mathsf{Q}}_P = (f_{\text{nonsep}})_*\sum_{a+b+c=k+2}f_{abc}[T'(a,b-1,0,c-1,0)],
    \]
    where $f_{abc}$ are the coefficients of $xyz(x-y)^2(x-z)^2(y-z)^2P$ as in Definition~\ref{Def-QP} and $f_{\text{nonsep}}:\Mbar_{0,4+2k}\to\Mbar_{1,2+2k}$ glues markings $3$ and $4$ as before.

    By the projection formula, we have the equality of intersection numbers
    \begin{equation}\label{eq:EdotQ}
    \mathsf{E}_P\cdot\widehat{\mathsf{Q}}_P = (f_{\text{nonsep}}^*\mathsf{E}_P)\cdot\sum_{a+b+c=k+2}f_{abc}[T'(a,b-1,0,c-1,0)].
    \end{equation}
    Recall that $f_{\text{nonsep}}^*\mathsf{E}_P$ was computed above in \eqref{eq:Epull} as a linear combination of boundary strata coming from the bead diagrams $T(a;i_1,i_2,j_1,j_2)$.

    Suppose $\Gamma,\Gamma'$ are stable graphs for $\Mbar_{0,4+2k}$ coming from labeling all the legs in bead diagrams of types $T,T'$ respectively, so $\Gamma$ has $k+1$ edges and $\Gamma'$ has $k$ edges. Then we claim that the intersection number $[\Gamma]\cdot[\Gamma']$ is equal to either $0$ or $(-1)^k$, and it is nonzero if and only if $\Gamma'$ is given by contracting the edge $e$ of $\Gamma$ labeled in the diagram below (the edge incident to $w$ pointing towards $v$).

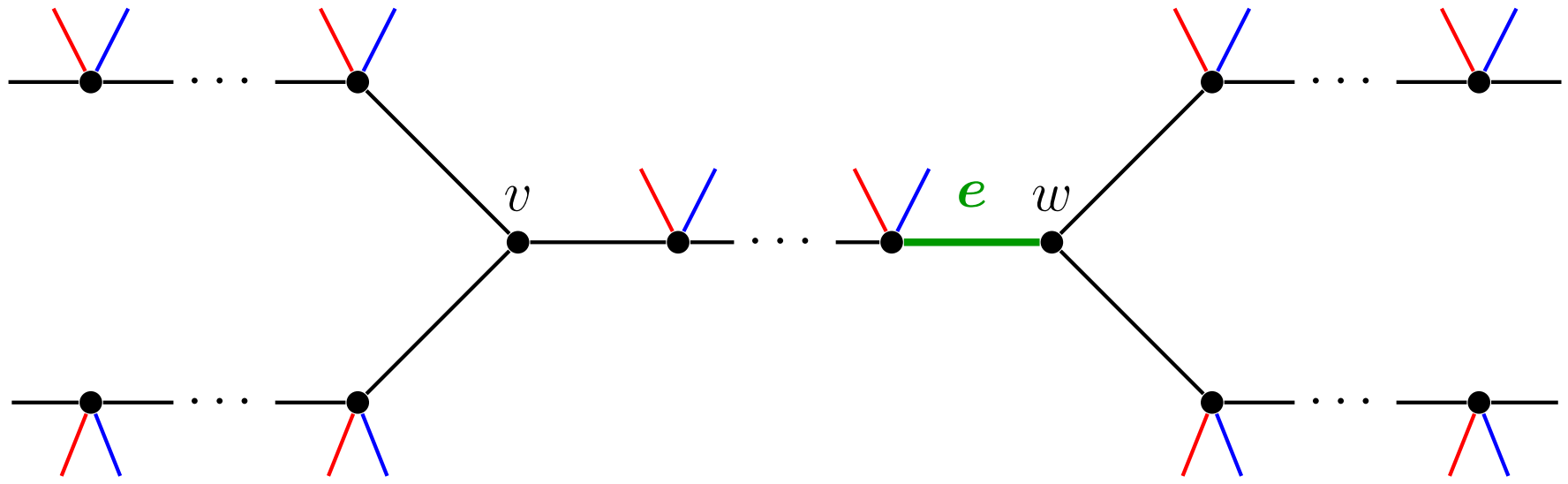
    \[
  \begin{tikzpicture}[x=8mm,y=8mm,moregraphs/center]
    \node[moregraphs/vertex,label=above:$v$] (v) at (0,0) {};
    \node[moregraphs/vertex,label=above:$w$] (w) at (5,0) {};
    \node[moregraphs/vertex] (middle-left) at (1.5,0) {};
    \node[moregraphs/vertex] (middle-right) at (3.5,0) {};
    \node[graphs/dots] (middle-dots) at (2.5,0) {$\cdots$};
    \draw[moregraphs/edge]
      (v) -- (middle-left)
      (middle-left) -- (middle-dots.west)
      (middle-dots.east) -- (middle-right);
    \draw[moregraphs/edge,green!60!black,line width=1.6pt]
      (middle-right) -- (w)
      node[midway,above=1mm] {$\boldsymbol{e}$};

    \foreach \side/\height/\legspread in {upper/1.5/2.8,lower/-1.5/2.2}{
      \pgfmathsetlengthmacro{\terminalleglength}{veclen(\legspread,5.5)*1mm}
      \node[moregraphs/vertex] (left-\side-end) at (-4,\height) {};
      \node[moregraphs/vertex] (left-\side-inner) at (-1.5,\height) {};
      \node[graphs/dots] (left-\side-dots) at (-2.75,\height) {$\cdots$};
      \draw[moregraphs/edge]
        (left-\side-end) -- (left-\side-dots.west)
        (left-\side-dots.east) -- (left-\side-inner)
        (left-\side-inner) -- (v);
      \draw[moregraphs/leg] (left-\side-end) -- ++(-\terminalleglength,0);

      \node[moregraphs/vertex] (right-\side-inner) at (6.5,\height) {};
      \node[moregraphs/vertex] (right-\side-end) at (9,\height) {};
      \node[graphs/dots] (right-\side-dots) at (7.75,\height) {$\cdots$};
      \draw[moregraphs/edge]
        (w) -- (right-\side-inner)
        (right-\side-inner) -- (right-\side-dots.west)
        (right-\side-dots.east) -- (right-\side-end);
      \draw[moregraphs/leg] (right-\side-end) -- ++(\terminalleglength,0);
    }

    \foreach \vertex in
      {middle-left,middle-right,left-upper-end,left-upper-inner,right-upper-inner,right-upper-end}{
      \draw[moregraphs/leg,red] (\vertex) -- ++(-2.8mm,5.5mm);
      \draw[moregraphs/leg,blue] (\vertex) -- ++(2.8mm,5.5mm);
    }
    \foreach \vertex in
      {left-lower-end,left-lower-inner,right-lower-inner,right-lower-end}{
      \draw[moregraphs/leg,red] (\vertex) -- ++(-2.2mm,-5.5mm);
      \draw[moregraphs/leg,blue] (\vertex) -- ++(2.2mm,-5.5mm);
    }
  \end{tikzpicture}
\]

    Indeed, suppose $G$ is a generic $(\Gamma,\Gamma')$-graph. Since every edge in both $\Gamma$ and $\Gamma'$ separates the graph into two components such that each component has at least one uncolored leg and the same number of red and blue legs, the same must be true of every edge in $G$. Splitting any vertex of $\Gamma$ will create an edge that violates this property, so we must have $G = \Gamma$. Then $\Gamma'$ must be given by contracting one edge of $\Gamma$ to form the valence $5$ vertex $w$. To give a nonzero intersection number, this edge must be along the path between $w$ and $v$ because otherwise the excess classes along the path between $v$ and $w$ will integrate to zero for dimension reasons. Finally, the resulting nonzero intersection number is $(-1)^k$ because it is obtained by integrating $-\psi$ on each of the $k$ valence $4$ vertices of $\Gamma$.

    We can now evaluate the right side of \eqref{eq:EdotQ}. Each boundary strata class $[\Gamma']$ coming from the $T'$ bead diagrams will have nonzero intersection number $(-1)^k$ with exactly one of the boundary strata classes $[\Gamma]$ coming from the $T$ bead diagrams. The coefficients with which these classes show up are equal (including possible signs coming from the antisymmetrization), so we conclude that
    \[
    \mathsf{E}_P\cdot\widehat{\mathsf{Q}}_P = C_k\sum_{a+b+c=k+2}f_{abc}^2
    \]
    for some nonzero constant $C_k$ coming from the various symmetrizations and antisymmetrizations. As a consequence, this intersection number is nonzero whenever $P$ is nonzero, since then at least one of the coefficients $f_{abc}$ is nonzero.
\end{enumerate}

\section{Relations from the boundary and the cohomology of local systems} \label{locsys}

\subsection{Ghosts}

\begin{definition}\label{ghost def}
    Let $X$ be an algebraic variety, and $Z \subset X$ a closed subvariety. A \emph{homological $Z$-ghost on $X$ of degree $k$} is a class in
    \[ \mathsf{Im}(\HH_k(Z') \to \HH_k(Z)) \cap \mathsf{Ker}(\HH_k(Z) \to \HH_k(X)),\]
    where $Z' \to Z$ is any resolution of singularities. The space of ghosts is independent of the choice of $Z'$, and we denote it $\ghost_k(X,Z)$.
\end{definition}

 According to \cref{qth_rel}, the class $\QQQ_{\Delta_0} \in \HH_{16}(\Delta_0)$ is a homological $\Delta_0$-ghost on $\Mbar_{2,14}$, and \cref{propdelq} says that it is a nontrivial $\Delta_0$-ghost. It is the first homological $\Delta_0$-ghost we have seen.{\footnote{See \cite[Section 5]{Danslecture} from 2013.}} The possibility of nontrivial $\Delta_0$-ghosts was discussed in \cite[Section 4.9.2]{P1} in the context of Chow classes.{\footnote{The terminology of \cite{P1} is not followed here. In our terminology, the discussion of \cite{P1} concerns {\em Chow $\Delta_0$-ghosts}. Examples from \cite{Liu-Pandharipande} were suggested in  \cite{P1} as candidates for $\Delta_0$-ghosts, but these have not been shown to be nonzero on
$\Delta_0$. Other variations of the definition are possible, including requiring the homology class to be algebraic or tautological.}} More generally, by Proposition \ref{propdel}, $(\QQQ_P)_{\Delta_0}$ is a $\Delta_0$-ghost for all nonzero $P\in\mathcal{S}_m$ such that $\QQQ_P$ is a cohomological relation, including all nonzero $P$ when $0\le m \le 2$ by Corollary~\ref{cor:theta_relations}.

We will only consider \cref{ghost def} in the case that $X=\Mbar_{g,n}$ and $Z=\Delta_0$ or $Z=\Delta$.

\begin{prop}\label{ghost criterion}
    Suppose that $X$ is smooth and proper and $Z \subset X$ a closed subvariety. Then there is a natural isomorphism
    \[ \ghost_k(X,Z) \cong \gr^W_{-k} \HH^\BM_{k+1}(X \setminus Z).\]
\end{prop}

\begin{proof}
The long exact sequence in Borel--Moore homology
    \[
    \dots\rightarrow \HH_{k+1}(X)\rightarrow \HH^\BM_{k+1}(X \setminus Z)\rightarrow \HH_{k}(Z)\rightarrow \HH_{k}(X)\rightarrow \dots
    \]
    respects the weights, and so we have an exact sequence
    \[
    0\rightarrow \gr^W_{-k}\HH^\BM_{k+1}(X \setminus Z)\rightarrow \gr^W_{-k}\HH_{k}(Z)\rightarrow \HH_{k}(X). 
    \]Now by Deligne \cite[Proposition (8.2.5)]{hodge3} one has \[ \gr^W_{-k}\HH_{k}(Z) = W_{-k}\HH_{k}(Z) = \operatorname{Im}(\HH_k(Z') \to \HH_k(Z)). \qedhere\] 
\end{proof}

\subsection{Main result}

The main theorem of Section \ref{locsys} concerns $\Delta_0$-ghosts and $\Delta$-ghosts on $\Mbar_{2,14}.$

\begin{thm}\label{thm section 4}The space of homological $\Delta_0$-ghosts on $\Mbar_{2,14}$ of degree $16$ is a single copy of the irreducible representation of $\SSS_{14}$ associated to the partition $(2^7)$. The natural map 
    \[ \ghost_{16}(\Mbar_{2,14},\Delta_0) \longrightarrow \ghost_{16}(\Mbar_{2,14},\Delta)\]
    is injective. In particular, since $\QQQ$ is a nontrivial $\Delta_0$-ghost, it is also a nontrivial $\Delta$-ghost. 
\end{thm}

Note that irreducibility shows that $ \ghost_{16}(\Mbar_{2,14},\Delta_0)$ must be generated as an $\SSS_{14}$-representation by $\QQQ$. 

In proving \cref{thm section 4} we will find it convenient to work with compact support cohomology instead of Borel--Moore homology. Since compact support cohomology is the linear dual of Borel--Moore homology, \cref{ghost criterion} shows that \cref{thm section 4} is equivalent to the following.

\begin{prop}\label{reduction}
 The group $\gr^W_{16} \HH^{17}_c(\M_{2,14}^{\ct})$ consists of a single copy of the irreducible representation of $\SSS_{14}$ associated to the partition $(2^7)$.  The natural map
	\[ \gr^W_{16} \HH^{17}_c(\M_{2,14}) \longrightarrow \gr^W_{16} \HH^{17}_c(\M_{2,14}^{\ct})\]
is surjective. 
\end{prop}

The {\em irreducibility} claim that $\gr^W_{16} \HH^{17}_c(\M_{2,14}^{\ct})$ is isomorphic to a single copy of $(2^7)$
follows  readily from results of Petersen \cite{Petersen2016-GenusTwo}. The {\em surjectivity} claim that the classes come from $\M_{2,14}$
is more difficult.{\footnote{The proof given here is independent of the substantially more complicated AI-generated argument, which the authors were not able to verify.}} The irreducibility claim already gives that $\ghost_{16}(\Mbar_{2,14},\Delta_0)$ is nontrivial, and generated as an $\SSS_{14}$-representation by $\QQQ$. The surjectivity claim will be needed to prove that $\QQQ \notin (\mathsf{3Spin}^9_{2,14})^{[2^7]}$ in Section \ref{S14rep}.

Proposition \ref{reduction} together with the results of Section \ref{S14rep} imply the existence of 
a class in $(\mathsf{CohRel}_{2,14}^9)^{[2^7]}$ which is not in $(\mathsf{3Spin}_{2,14}^9)^{[2^7]}$. But the
methods of Sections \ref{locsys} and \ref{S14rep} do not produce a formula. Our study of  $\QQQ$ 
depends upon an interplay between the results of Sections \ref{secvan}-\ref{S14rep}.

\subsection{Preliminaries for Proposition \ref{reduction}}
Before proving \cref{reduction} we will collect some preliminaries. 

\subsubsection{Cohomology of symplectic local systems}We will need to know something about the cohomology of local systems on $\M_2^{\ct}$ and on $\Sym^2 \M_{1,1}$, associated to algebraic representations of $\mathrm{Sp}_4$ and to $(\mathrm{SL}_2 \times \mathrm{SL}_2) \rtimes \SSS_2$, respectively; we call these \emph{symplectic local systems} for short. The cohomology of these symplectic local systems is well understood in all degrees. On $\Sym^2 \M_{1,1}$ all cohomology is determined by the classical Eichler--Shimura description of the cohomology of modular curves, as explained in \cite{d-elliptic}. The cohomology on $\M_2^\ct = \mathcal A_2$ is determined up to semisimplification in \cite{localsystemsA2}. A disclaimer is that \cite{localsystemsA2} states the result both for $\ell$-adic Galois representations  and for $\Q$-Hodge structures, but the arguments there only give the $\mathbb R$-Hodge structure. This gap has recently been filled by Gu--Lee--Prasanna \cite{GLP}. More precisely, they show that the Saito--Kurokawa lifting is compatible with the $\Q$-Hodge structure.

Propositions \ref{LS A2} and \ref{LS A11} isolate the precise properties of the cohomologies of these local systems which we shall use.

\begin{prop}\label{LS A2}
	Let $\mathbb V_{a,b}$ be a symplectic local system on $\M_2^{\ct}$. Then \[ \gr^W_{p+a+b-1} \HH^p_c(\M_2^{\ct},\mathbb V_{a,b}) \cong \begin{cases} \oplus_{f \in \mathcal S'_{a+b+4}} \qq(-b-1) & p=3, a=b \\ 0 & \text{otherwise} \end{cases}\]
	where $\mathcal S_k'$ denotes the set of normalized Hecke eigenforms for $\operatorname{SL}_2(\mathbb Z)$ of weight $k$, such that the central value $L(f,\tfrac 1 2)$ is zero. In particular, the unique symplectic local system of smallest weight with cohomology of weight exactly one less than pure weight is $\mathbb V_{7,7}$, for which 
\[ \gr^W_{16} \HH^3_c(\M_2^{\ct},\mathbb V_{7,7}) \cong \qq(-8)\,.\]\end{prop}
	
	\begin{rem}\label{nonvanishing}
		One knows that $L(f,\tfrac 1 2)=0$ for all eigenforms $f$ of weight $k \equiv 2 \pmod 4$, since the functional equation forces it to vanish of odd order. Conjecturally, eigenforms for $\operatorname{SL}_2(\mathbb Z)$ of weight $k \equiv 0 \pmod 4$ never have vanishing central value; this has been computationally verified for weights up to 14000 \cite{ghitza-mcandrew}.  The class in the cohomology of $\mathbb V_{7,7}$ is associated to the unique Hecke eigenform for $\operatorname{SL}_2(\mathbb Z)$ of weight $18$; this is the lowest weight $\equiv 2 \pmod 4$ for which a nonzero cusp form for $\operatorname{SL}_2(\mathbb Z)$ exists. 
	\end{rem}
	
	\begin{prop}\label{LS A11} If $\mathbb V$ is a symplectic local system of weight $w$ on $\Sym^2 \M_{1,1}$, then 
\[\gr^W_{p+w-1} \HH_c^p(\Sym^2 \M_{1,1},\mathbb V)=0\,.\]
\end{prop}

\begin{rem}In the Eichler--Shimura description of the cohomology of modular curves, the contributions from cusp forms give rise to pure cohomology, hence only contribute to pure weight. Classes associated to Eisenstein series of weight $k$ are of weight exactly $(k-1)$ less than pure weight; in particular, they cannot contribute to weight exactly one less than pure weight in level $1$, as there are no weight $2$ Eisenstein series for $\operatorname{SL}_2(\mathbb Z)$. 
	
\end{rem}

\subsubsection{The Decomposition Theorem and forgetting the marked points}

Results of Petersen \cite{Petersen2016-GenusTwo} allow the computation of $H^*_c(\M_{2,n}^\ct)$ for all $n$, as an $\SSS_n$-representation and a mixed Hodge structure up to semisimplification. The key input is that the Decomposition Theorem  of Beilinson--Bernstein--Deligne--Gabber \cite{BBD} takes the following simple form in the case of $\M_{2,n}^\ct \to \M_2^\ct$ \cite[Theorem 2.1]{Petersen2016-GenusTwo}.

\begin{prop} \label{bbd}
	Let $f: \M_{2,n}^{\ct} \to \M_2^{\ct}$. The complex $Rf_*\qq$ is isomorphic to its cohomology in the derived category: \[ Rf_*\qq \cong \bigoplus_q R^qf_*\qq[-q]\,.\]
	Each group $R^q f_* \qq$ decomposes as a direct sum, with each summand a weight $q$ symplectic local system on $\M_2^{\ct}$, or a symplectic local system on $\Sym^2 \M_{1,1}$ pushed forward to $\M_2^{\ct}$.
\end{prop}

The summands of $Rf_*\qq$ with full support, i.e.\ those which are local systems on $\M_2^{\ct}$, are determined by the cohomology of a single smooth fiber, i.e.\ the fibers of $\M_{2,n}^{\rt} \to \M_2$. Here, $\M_{2,n}^{\rt}$ is the moduli space of curves (of genus $2$ with $n$ markings) with rational tails. For the following see e.g. Petersen--Tommasi \cite[Proposition 4.5]{petersentommasi}, together with the discussion in the paragraph immediately before and after the proposition. 
\begin{prop}\label{which local systems occur}Let $g : \M_{2,n}^{\rt} \to \M_2$, and $h : \CC_2^n \to \M_2$. Both $Rg_*\qq$ and $Rh_*\qq$ decompose in the derived category as in \cref{bbd}. The summands are all Tate twists of local systems $\mathbb V_{a,b}$ with $a+b \leq n$. Moreover, local systems with $a+b=n$ occur only in the middle degree, with the same multiplicity in both $R^n g_*\qq$ and $R^n h_*\qq$, namely as a summand 
	\[ \bigoplus_{a+b=n} \mathbb V_{a,b} \otimes \sigma_{(2^b,1^{a-b})}\,,\]
where $\sigma_\lambda$ denotes the representation of $\SSS_n$ associated to the partition $\lambda$. This summand is mapped isomorphically under the natural inclusion $Rh_*\qq \hookrightarrow Rg_*\qq$ coming from the factorization $\M_{2,n}^{\rt} \to \CC_2^n \to \M_2$. \end{prop}

Propositions \ref{bbd} and \ref{which local systems occur} were stated here in genus $2$, since that is what we will use and what the cited references say, but they remain true with easy modifications also in higher genus. By contrast, the following Propositions \ref{pure weight} and \ref{multiplicative decomposition} are very specific to genus $2$. 

\begin{prop}\label{pure weight}
	Let $X$ be a genus $2$ curve, and let $X^n \to \CC_2^n$ be the inclusion of the corresponding fiber of $\CC_2^n \to \M_2$.  The image of the restriction map $$\HH^*(\CC_2^n) \longrightarrow \HH^*(X^n)$$ consists exactly of the monodromy invariants $\HH^*(X^n)^{\mathrm{Sp}_4}$, and it induces an isomorphism 
	\[ \RH^*(\CC_2^n) \stackrel \sim \longrightarrow \HH^*(X^n)^{\mathrm{Sp}_4}\,.\]
\end{prop}

The first part above is just Deligne's invariant cycle theorem \cite[Théorème 4.1.1(ii)]{hodge2}. The second part can be seen as a compositum of a result of Tavakol \cite{tavakol2}, that $\R^*(\CC_2^n) = \R^*(X^n)$ and both are Gorenstein, and hence agree with their image in cohomology, and Petersen--Tommasi \cite[Section 4]{petersentommasi}, who prove for $X$ of any genus that $\RH^*(X^n) \cong \HH^*(X^n)^{\mathrm{Sp}_{2g}}$. 

\begin{rem}\label{monodromy invariants}Since $\RH^{k}(\CC_2^n) = W_{2k}\HH^{2k}(\CC_2^n)$ for all $n<20$ and any $k$, by Petersen--Tommasi \cite{petersentommasi} and Petersen \cite{Petersen2016-GenusTwo}, we see in particular that $W_{2k}\HH^{2k}(\CC_2^n) \stackrel \sim \longrightarrow \HH^{2k}(X^n)^{\mathrm{Sp}_4}$ in this range.
\end{rem}

The following is explained by Petersen--Tavakol--Yin \cite[\S 5.2.3]{petersentavakolyin} on the level of Chow motives; the result here follows by Betti realization. 

\begin{prop}\label{multiplicative decomposition}
	Let $\pi : \CC \to \M_2$ be the universal curve. Then $R\pi_*\qq$ admits a {multiplicative decomposition}: there is an isomorphism
	\begin{equation}\label{derived cat iso} R\pi_*\qq\cong \bigoplus_q R^q \pi_*\qq [-q]\end{equation}
	in the derived category, compatible with cup-product, in the sense that it makes the diagram
	\[ \begin{tikzcd}
		R\pi_*\qq \otimes R\pi_*\qq \arrow[r]\arrow[d]& \arrow[d]R\pi_*\qq \\
		\bigoplus_{q_1,q_2} R^{q_1} \pi_*\qq \otimes R^{q_2} \pi_* \qq[-q_1-q_2] \arrow[r]& \bigoplus_q R^q \pi_*\qq [-q]
	\end{tikzcd}\]
commute.
\end{prop}

This result may seem abstract, but has concrete consequences. Any choice of isomorphism \eqref{derived cat iso} induces an isomorphism
\begin{equation} \label{coh iso}\HH^*_c(\CC_2^n) \cong \bigoplus_{q_1, \dots,q_n} \HH^{*-q_1-\dots-q_n}_c(\M_2,\bigotimes_{j=1}^n R^{q_j}\pi_*\qq)\end{equation}
for all $n$.
Indeed: this is because the left-hand side of \eqref{coh iso} is canonically $\mathbb H^*_c(\M_2, (R\pi_*\qq)^{\otimes n})$ and the right-hand side of \eqref{coh iso} is $\mathbb H^*_c(\M_2, (\bigoplus_q R^q \pi_*\qq [-q])^{\otimes n})$. If the isomorphism \eqref{derived cat iso} is multiplicative, then the isomorphism \eqref{coh iso} behaves well with respect to diagonal inclusions between different fibered powers of the universal curve. Specifically, consider a diagonal embedding $\Delta_{ij} : \CC_2^{n-1} \to \CC_2^n$. Multiplicativity implies that the restriction map
\[ \Delta_{ij}^* : \HH^*_c(\CC_2^n) \longrightarrow \HH^*_c(\CC_2^{n-1})\]
agrees under the isomorphism \eqref{coh iso} with the map obtained by applying $\HH^*_c(\M_2,-)$ to the map of local systems induced by the cup-product
\[ R^{q_i}\pi_*\qq \otimes R^{q_j} \pi_*\qq \longrightarrow R^{q_i+q_j} \pi_*\qq\,.\]

\subsubsection{The spectral sequence associated to the stratification of \texorpdfstring{$\CC_g^n$}{Cgn}}

The final ingredient we need is a spectral sequence computing $\HH^*_c(\M_{g,n})$ from the data of $\HH^*_c(\CC_g^m)$ for $1 \leq m \leq n$. This is a special case of the spectral sequence associated to a stratified space, introduced by Petersen \cite{spectralsequencestratification}. We consider the stratification of $\CC_g^n$ by points colliding. The poset of strata is the partition lattice, and the dense open stratum is $\M_{g,n}$. From the construction of \cite{spectralsequencestratification} we obtain 
\begin{equation} \label{SS}E_1^{pq} = \HH^q_c(\CC_g^{n-p})^{\oplus c(n,n-p)} \implies \HH^{p+q}_c(\M_{g,n})\,.\end{equation}
Here $c(n,k)$ denotes the unsigned Stirling number of the first kind. This is a spectral sequence of mixed Hodge structures, and it is functorial in maps of stratified spaces. It is a version ``in families'' of the Totaro spectral sequence \cite{totaro}, or rather its more general Poincar\'e dual version in compact support cohomology due to Getzler \cite{getzler99}. The $E_1$-differential is given by the natural maps arising from diagonal inclusions between the spaces $\CC_g^m$. For example, $c(n,n-1) = \binom n 2$, and the $E_1$-differential 
\[ E_1^{0,q} = \HH^q_c(\CC_g^n) \longrightarrow \bigoplus_{\binom n 2} \HH^q_c(\CC_g^{n-1}) = E_1^{1,q}\] 
is the sum of the restriction maps in cohomology along the $\binom n 2$ diagonals in $\CC_g^n$.

\begin{figure}[h]
	\centering
	\begin{tikzpicture}
		
		\matrix (m) [
		matrix of math nodes,
		nodes={
			minimum width=8ex,
			minimum height=5ex,
			anchor=center
		},
		column sep=1.5ex,
		row sep=1ex
		]{
			\HH_c^3(\CC_2^4)
			& \HH_c^3(\CC_2^3)^{\oplus 6}
			& \HH_c^3(\CC_2^2)^{\oplus 11}
			& \HH_c^3(\CC_2^1)^{\oplus 6} \\		
			\HH_c^2(\CC_2^4)
			& \HH_c^2(\CC_2^3)^{\oplus 6}
			& \HH_c^2(\CC_2^2)^{\oplus 11}
			& \HH_c^2(\CC_2^1)^{\oplus 6} \\
			\HH_c^1(\CC_2^4)
			& \HH_c^1(\CC_2^3)^{\oplus 6}
			& \HH_c^1(\CC_2^2)^{\oplus 11}
			& \HH_c^1(\CC_2^1)^{\oplus 6} \\
			\HH_c^0(\CC_2^4)
			& \HH_c^0(\CC_2^3)^{\oplus 6}
			& \HH_c^0(\CC_2^2)^{\oplus 11}
			& \HH_c^0(\CC_2^1)^{\oplus 6} \\
		};
		
		\coordinate (O) at (m.south west);
		\coordinate (X) at (m.south east);
		\coordinate (Y) at (m.north west);
		
		\draw[->] (O) -- (X) node[right] {$p$};
		\draw[->] (O) -- (Y) node[above] {$q$};
		\node[below=0.25em] at (m-4-1 |- O) {$0$};
		\node[below=0.25em] at (m-4-2 |- O) {$1$};
		\node[below=0.25em] at (m-4-3 |- O) {$2$};
		\node[below=0.25em] at (m-4-4 |- O) {$3$};
		
		\node[left=0.25em] at (O |- m-4-1) {$0$};
		\node[left=0.25em] at (O |- m-3-1) {$1$};
		\node[left=0.25em] at (O |- m-2-1) {$2$};
		\node[left=0.25em] at (O |- m-1-1) {$3$};
		
		\node[above=0.8em] at (m-1-1) {$\vdots$};
		\node[above=0.8em] at (m-1-2) {$\vdots$};
		\node[above=0.8em] at (m-1-3) {$\vdots$};
		\node[above=0.8em] at (m-1-4) {$\vdots$};
		
	\end{tikzpicture}
    	\captionsetup{width=0.75\textwidth}

	\caption{The $E_1$ page of the spectral sequence associated to the stratification of $\CC_g^n$ by diagonal loci, written here for $g=2$ and $n=4$.}
\end{figure}

Let $X$ be a fixed curve of genus $g$. The inclusion $X^n \hookrightarrow \CC_g^n$ is a regular embedding of codimension $3g-3$. As such, there is a Gysin homomorphism 
\[ \HH^{*-6g+6}_c(X^n)(-3g+3) \longrightarrow \HH^*_c(\CC_g^n)\,.\]
The embedding is transverse to the stratification, and therefore induces a Gysin homomorphism between the spectral sequences associated to the stratifications of $\CC_g^n$ and $X^n$, respectively: we find 
\begin{equation}\label{map of SS}\begin{tikzcd}
	E_1^{pq} = \HH^q_c(\CC_g^{n-p})^{\oplus c(n,n-p)} \arrow[r,Rightarrow]& \HH^{p+q}_c(\M_{g,n}) \\ \widetilde E_1^{pq} = \arrow[u]\HH^{q-6g+6}_c(X^{n-p})^{\oplus c(n,n-p)}(-3g+3) \arrow[r,Rightarrow]& \HH^{p+q-6g+6}_c(F(X,n))(-3g+3)\,. \arrow[u]
\end{tikzcd}\end{equation}
Here $F(X,n)$ denotes the configuration space of $n$ distinct ordered points on $X$. This is a map of spectral sequences of mixed Hodge structures. 
 
	\subsubsection{Proof of Proposition \ref{reduction}.}
	
\begin{proof}[Proof of \cref{reduction}]
	We begin with the first claim, that $\gr^W_{16} \HH^{17}_c(\M_{2,14}^{\ct})$ consists of a single copy of the irreducible representation of $\SSS_{14}$ associated to the partition $(2^7)$. By \cref{bbd} we have
	\[ \HH^*_c(\M_{2,n}^{\ct}) \cong \bigoplus_{q} \HH^{*-q}_c(\M_2^{\ct},R^q f_* \qq)\,,\]
	and $R^qf_*\qq$ decomposes as a sum of weight $q$ local systems on $\M_2^{\ct}$ and $\Sym^2 \M_{1,1}$. The latter cannot contribute to weights one below pure by \cref{LS A11}. The former summands contribute to weights one below pure for the first time when $n=14$, by \cref{LS A2} and \cref{which local systems occur}. Namely $\mathbb V_{7,7} \otimes \sigma_{(2^7)}$ is a summand of $R^{14} f_*\qq$, where $\sigma_{(2^7)}$ denotes the irreducible representation of $\SSS_{14}$ associated to the partition $(2^7)$, and this gives $$ \gr^W_{16 }\HH^{17}_c(\M_{2,14}^{\ct}) \cong \qq(-8) \otimes \sigma_{(2^7)}\,.$$

	For the second part of the proposition, we first consider the portion of the long exact sequence $$ \HH^3_c(\M_2,\mathbb V_{7,7}) \to \HH^3_c(\M_2^{\ct},\mathbb V_{7,7}) \to \HH^3_c(\Sym^2 \M_{1,1},\mathbb V_{7,7})\,.$$ The target contains nothing of weight one below pure, by \cref{LS A11}. Hence
	\[ \gr^W_{16} \HH^3_c(\M_2,\mathbb V_{7,7}) \longrightarrow \gr^W_{16} \HH^3_c(\M_2^{\ct},\mathbb V_{7,7}) \cong \qq(-8)\] is surjective (in fact an isomorphism, but we will not need this). 
We consider the diagram
\[\begin{tikzcd}
	\HH^{17}_c(\M_{2,14}) \arrow[r]\arrow[rd]& \arrow[r]\HH^{17}_c(\M_{2,14}^{\rt}) & \HH^{17}_c(\M_{2,14}^{\ct}) \\ 
	& \HH^{17}_c(\CC_2^{14})\arrow[u]&
\end{tikzcd}\]
By \cref{which local systems occur} we see that both middle terms contain $\HH^3_c(\M_2,\mathbb V_{7,7}) \otimes \sigma_{2^7}$ as a summand, mapping to $\HH^3_c(\M_2^{\ct},\mathbb V_{7,7}) \otimes \sigma_{(2^7)} \subset \HH^{17}_c(\M_{2,14}^{\ct})$. It thus suffices to show that the map $\gr^W_{16}\HH^{17}_c(\M_{2,14}) \to \gr^W_{16}\HH^{17}_c(\CC_2^{14})$ hits the summand $\gr^W_{16} \HH^3_c(\M_2^{\ct},\mathbb V_{7,7}) \otimes \sigma_{(2^7)}$. To do so, we analyze the spectral sequence \eqref{SS}. We have a summand $$ \HH^3_c(\M_2,\mathbb V_{7,7}) \otimes \sigma_{(2^7)} \subset \HH^{17}(\CC_2^{14}) = E_1^{0,17}\,.$$ 
We first argue that this summand is in the kernel of the $E_1$-differential. This follows because the $E_1$-differential is given by the sum of all diagonal embeddings $\CC_2^{13} \hookrightarrow \CC_2^{14}$, and the discussion immediately after \cref{multiplicative decomposition}. Namely, $\mathbb V_{7,7} \otimes \sigma_{(2^7)}$ is a summand of $(R^1\pi_*\qq)^{\otimes 14}$, and it lies in the kernel of all the $\binom {14} 2$ cup-product maps $R^1 \pi_*\qq \otimes R^1\pi_*\qq \to R^2\pi_*\qq$, by Weyl's construction of the irreducibles of the symplectic group. 

The $E_2$-differential is $E_2^{0,17} \longrightarrow E_2^{2,16}$. We claim that $\gr^W_{16} E_2^{2,16} = 0$, so that the $E_2$-differential out of $E_2^{0,17}$ vanishes on the weight $16$ subquotient. To see this, we consider the map of spectral sequences \eqref{map of SS}. By the Poincar\'e dual of \cref{monodromy invariants}, this map induces an isomorphism 
\[ (\widetilde E_1^{pq})^{\mathrm{Sp}_4} \longrightarrow \gr^W_q E_1^{pq} \]
between the monodromy-invariant part of the first page of the spectral sequence associated to the stratification of $X^{14}$, and the ``pure part'' of the first page of the spectral sequence associated to the stratification of $\CC_2^{14}$, compatibly with the $E_1$-differentials, when $q$ is even. Now the spectral sequence $\widetilde{E}_r^{pq}$ converges to $\HH^{p+q-6}_c(F(X,14))(-3)$, and degenerates at the second page since no differential after the first is compatible with weights. The cohomological dimension of $F(X,n)$ is $n+1$, so $\HH^k_c(F(X,14))=0$ for $k<13$. Hence $\widetilde E_2^{pq}=0$ for $p+q < 13+6=19$. Hence also $\gr^W_q E_2^{pq}=0$ for $p+q<19$, when $q$ is even. This gives the desired vanishing of the $E_2$-differential.

All further differentials emanating out of $\gr^W_{16} E_r^{0,17}$ are zero for weight reasons, since $E_r^{pq}$ is of weight $\leq q$ for all $p,q,r$.

It follows that our summand $\qq(-8)\otimes \sigma_{(2^7)} \subset \gr^W_{16} E_1^{0,17}$ is in the kernel of all differentials, and hence in the image of the edge-homomorphism from the cohomology of $\M_{2,14}$. This is what we wanted to prove. \end{proof}

\subsection{New and old ghosts}

The construction $\ghost_k(X,Z)$ has certain functorialities in $X$ and $Z$ that we now want to make explicit. 

\begin{enumerate}[(i)]
    \item If $\alpha \in \HH^p(X)$ is a cohomology class, then cap product defines a map
    \[ \alpha \cap (-) : \ghost_k(X,Z) \to \ghost_{k-p}(X,Z)\,.\]
    \item Let $f:(X,Z) \to (Y,W)$ be a map of pairs of varieties. Then there is  a pushforward map
    \[ f_* : \ghost_k(X,Z) \to \ghost_k(Y,W)\,.\]
    \item If in addition $X$ and $Y$ are smooth and proper, and $Z=f^{-1}(W)$, then there is an exceptional pullback map
    \[ f^! : \ghost_k(Y,W) \to \ghost_{k+2d-2e}(X,Z)\,,\]
    with $d=\dim(X)$, $e=\dim(Y)$.
\end{enumerate}
Point (i) is immediate from the definition and functoriality of cap product. For point (ii) we need that if $\alpha \in \HH_k(Z)$ is pushed forward from some resolution of singularities then so is $f_*\alpha$, which is \cite[Theorem 1.1(b)]{weber}. (When $Z$ and $W$ are proper, which is all we need, this is also immediate from mixed Hodge theory.) Point (iii) is \cref{ghost criterion} and Poincar\'e duality:
\[ \ghost_k(Y,W) \cong \gr^W_{-k}\HH_{k+1}^\BM(Y \setminus W) \cong \gr^W_{2e-k}\HH^{2e-k-1}(Y \setminus W)\]
which maps by pullback in cohomology to $\gr^W_{2e-k}\HH^{2e-k-1}(X \setminus Z) \cong \ghost_{k+2d-2e}(X,Z)$. 

\vspace{5pt}
\newcommand{\G}{\mathsf{G}}
\noindent The above constructions can be used to build new ghosts from existing ones. Indeed, let $$\G_{k,n}=\ghost_k(\Mbar_{2,n},\Delta_0).$$ The functorialities (i)---(iii) define maps:
\begin{equation}\label{eq:ghost-operations}
\begin{alignedat}{2}
    \pi^!                  &\colon \G_{k-2,n-1}  &\! \rightarrow \,\, & \G_{k,n} \\
    \psi_i \cap \pi_i^!(-)    &\colon \G_{k,n-1}  & \! \longrightarrow \,\, & \G_{k,n} \\
    (\xi_{ij})_*            &\colon \G_{k,n-1}  & \! \longrightarrow \,\,&  \G_{k,n},
\end{alignedat}
\end{equation}
where $\pi_i \colon \Mbar_{2,n} \to \Mbar_{2,n-1}$ forgets the $i$th marking and $\xi_{ij} \colon \Mbar_{2,n-1} \to \Mbar_{2,n}$ glues on a three-pointed rational component with markings $i$ and $j$. 
\begin{definition}
    A homological $\Delta_0$-ghost in $\G_{k,n}$ is \emph{old} if it is in the span of the image of the maps \eqref{eq:ghost-operations}. We define $\G_{k,n}^{\mathrm{new}}$ to be the quotient of $\G_{k,n}$ by the subspace of old ghosts. 
\end{definition}

We now obtain a very precise description of the space of new ghosts in terms of spaces of modular forms.

\begin{thm} \label{newghosts} The space $\G_{k,n}^{\mathrm{new}}$ is nonzero only for $n=2a$ and $k=2a+2$, and consists of $\dim \mathcal S_{2a+4}'$ copies of the representation $(2^a)$ of $\SSS_n$.
\end{thm}

We will deduce \cref{newghosts} from the following, which is  similar to \cite[Lemma 3.1]{STE}.

\begin{lem}
    Let $X$ be a smooth compact curve, and $X[n]$ the Fulton--MacPherson compactification of $n$ points on $X$. Consider the maps 
    \begin{equation}\label{eq:ghost-operations2}
\begin{alignedat}{2}
    \pi_i^*                  &\colon \HH^*(X[n-1])  &\! \longrightarrow \,\, & \HH^*(X[n]) \\
    \alpha_i \cdot \pi_i^*(-)    &\colon \HH^{*-2}(X[n-1])  & \! \to \,\, & \HH^*(X[n]) \\
    (\xi_{ij})_!            &\colon \HH^{*-2}(X[n-1])  & \! \to \,\,&  \HH^*(X[n]),
\end{alignedat}
\end{equation}analogous to \eqref{eq:ghost-operations}, and where $\alpha$ is any generator of $\HH^2(X)$. The joint cokernel of \eqref{eq:ghost-operations2} is concentrated in degree $n$, and as a representation of the mapping class group it is isomorphic to $\bigoplus_{|\lambda|=n} \mathbb V_\lambda \otimes \sigma_{\lambda^T}$. 
\end{lem}

\begin{proof}
    We use the Fulton--MacPherson presentation \cite[\S5]{FM94}, see also \cite{FMcorrection}, which
describes $\HH^*(X[n])$ as an algebra over $\HH^*(X^n)$
generated by divisor classes $D_S$, with $S \subseteq \{1,\dots,n\}$ and $|S| \geq 2$. 

Let $I$ be the span of the images of all maps \eqref{eq:ghost-operations2}. We first show that the ideal generated by all $D_S$ lies in $I$. Indeed pick an element 
\begin{equation}\label{an element}
    \gamma \cdot D_{S_1} \cdots D_{S_k}
\end{equation} 
with $\gamma \in \HH^*(X^n).$ We prove by induction on the quantity $\min_t |S_t|$ that \eqref{an element} lies in $I$. Assume $S_1$ has minimal cardinality. Pick  $i,j \in S_1$. If $S_1=\{i,j\}$, then the product \eqref{an element} lies in the image of $(\xi_{ij})_!$ by the projection formula. For the inductive step, let $r : \HH^*(X[n]) \to \HH^*(X[n])$ be the ring homomorphism $\pi_i^*\xi_{ij}^*$. 
We have 
\[ r(D_S) = \begin{cases} D_S + D_{S \setminus \{i\}} & i,j \in S \\
D_S + D_{S \cup \{i\}} & i,j \notin S \\
0 & \text{otherwise.}
     \end{cases}
\]
It follows that 
\[ r(\gamma) r(D_{S_1}) \cdots r(D_k) \equiv r(\gamma) D_{S_1} \cdots D_{S_k} \mod I\]
since if we expand the product on the left-hand side then all other terms are either in $I$ by induction, or zero by the overlap relations, or just zero. Since moreover $r(\gamma) D_{S_1} = D_{S_1}$ we see that \eqref{an element} is in $I$.  

Now the quotient of $\HH^*(X[n])$ by the ideal generated by all $D_S$ is $\HH^*(X^n)$ modulo the ideal generated by all diagonals. Taking the quotient by the first and second maps in \eqref{eq:ghost-operations2} further kills all K\"unneth components with a factor $\HH^0$ and $\HH^2$. We are left with $\HH^1(X)^{\otimes n}$ modulo all diagonals, which by symplectic Schur--Weyl duality is precisely $\bigoplus_{|\lambda|=n} \mathbb V_\lambda \otimes \sigma_{\lambda^T}$.
\end{proof}

We now give the proof of \cref{newghosts}. Since we will not directly use this result in the rest of the paper, we allow ourselves to be brief.

\begin{proof}[Proof of \cref{newghosts}]
    By \cref{ghost criterion} and Poincar\'e duality, we are equivalently considering 
    \[ \G_{k,n} \cong \gr^W_{6+2n-k} \HH^{5+2n-k}(\M_{2,n}^\ct)\]
    modulo the image of the maps defined as in \eqref{eq:ghost-operations2}. As in \cref{bbd} we decompose $\HH^*(\M_{2,n}^\ct)$ into summands by decomposing the complex $Rf_*\Q$, where $f: \M_{2,n}^\ct \to \M_2^\ct$. By \cref{LS A11} and \cref{LS A2} the only contributions to $\gr^W_{*+1} \HH^*(\M_{2,n}^\ct)$ come from the local systems of the form $\mathbb V_{a,a}$ on $\M_2^\ct$. The summands of full support in $R^q f_*\Q$ are determined by their restriction to $\M_2$, and hence given by the action of the mapping class group on $\HH^q(X[n])$, where $X$ is a smooth genus two curve. Working modulo the image of the maps \eqref{eq:ghost-operations2} now leaves only the space
    \[ \G_{k,n}^{\mathrm{new}} \cong \bigoplus_{|\lambda|=n} \gr^W_{6+2n-k} \HH^{5+2n-k-n}(\M_2^\ct,\mathbb V_\lambda) \otimes \sigma_{\lambda^T}. \]
    By \cref{LS A2}, we only have a nonzero contribution when $\lambda = (a,a)$, and $5+2n-k-n=3$, whence $n=2a$ and $k=2a+2$, and the precise contribution from $\mathbb V_{a,a}$ in \cref{LS A2} gives exactly the desired formula for $\G_{k,n}^{\mathrm{new}}$. 
\end{proof}
\newcommand{\sgn}{\mathrm{sgn}}
\begin{rem}\cref{newghosts} is interesting in  light of \cref{propdel}. For example, \cref{propdel} shows the existence of a nontrivial $\Delta_0$-ghost on $\Mbar_{2,16}$, but \cref{newghosts} shows that $\G_{18,16}^{\mathrm{new}} = 0$. So the ghost must be old. Similarly on $\Mbar_{2,18}$, \cref{propdel} constructs two $\Delta_0$-ghosts. But  both transform according to the representation $\mathsf{sgn} \boxtimes \mathsf{sgn}$ of $\SSS_9 \times \SSS_9$, and the $\mathsf{sgn} \boxtimes\mathsf{sgn}$-isotypical component of $\G_{20,18}^{\mathrm{new}} \cong \sigma_{(2^9)}$ is one-dimensional. So a linear combination of the two classes is old. It would be interesting to understand better which ghosts constructed from \cref{propdel} are old or new. 
\end{rem}

\begin{question}\label{q:moreghosts-non3spin}
    Do the $\Delta_0$-ghosts constructed in  \cref{newghosts} lie outside of $\mathsf{3Spin}$?
\end{question}
\begin{rem}\cref{newghosts} produces nontrivial new ghosts for all $n \equiv 2 \pmod 4$, $n \geq 14$. 
    When $n=14$, the answer to \cref{q:moreghosts-non3spin} is yes (Theorem \ref{qth_not3spin}). When $n=18$, the ghosts from \cref{newghosts} also lie outside of $\mathsf{3Spin}^{11}_{2,18}$, because the proof of  \cref{thm section 4} goes through, and so does the subsequent analysis in Section \ref{S14rep}. We do not know which of the $\Delta_0$-ghosts constructed in \cref{newghosts} survive under the map
    \[ \ghost_{n+2}(\Mbar_{2,n},\Delta_0) \longrightarrow \ghost_{n+2}(\Mbar_{2,n},\Delta)\]
	for $n\geq 22$, because the proof of \cref{thm section 4} uses Remark \ref{monodromy invariants}, which requires $n<20$.
\end{rem}

\section{Representation theory of \texorpdfstring{$3$}{3}-spin relations} \label{S14rep}

\subsection{\texorpdfstring{$3$}{3}-spin relations}
We review the standard set of $\mathbb{Q}$-linear generators for the ideal of $3$-spin relations $\mathsf{3Spin}^*_{g,n}$, following the notation of \cite{Pixton2012-ConjecturalRelations}. Suppose $g,n,r$ are nonnegative integers, $\sigma$ is a partition with no parts of size $2$ mod $3$, and $a_1,\ldots,a_n$ are nonnegative integers that are not $2$ mod $3$. If the inequality
\begin{equation}\label{eq:3spin_ineq}
3r\ge g+1+|\sigma|+\sum_{i=1}^na_i\, ,
\end{equation}
holds, then  a relation
\[
\mathcal{R}(g,n,r;\sigma,a_1,\ldots,a_n)\in\mathsf{S}^r_{g,n}\, 
\]
is described in \cite{Pixton2012-ConjecturalRelations}.
We will call such relations {\it interior relations}. The action of the symmetric group $\SSS_n$ by relabeling points acts on these relations by permuting the $a_i$.

The ideal of $3$-spin relations $\mathsf{3Spin}^*_{g,n}$ is then $\mathbb{Q}$-linearly spanned by these interior relations along with boundary pushforwards
\begin{equation}\label{eq:boundary_rels}
(j_{\Gamma})_*\left(x_1\boxtimes x_2\boxtimes\cdots\boxtimes x_m\right)\, ,
\end{equation}
where $\Gamma$ is a nontrivial stable graph (with $m$ vertices) and one of the $x_i$ insertions is equal to some interior relation $\mathcal{R}(g_v,n_v,r_v;\sigma,a_1,\ldots,a_{n_v})$. 

By the following result,
interior relations do not
contribute to $(\mathsf{3Spin}_{2,14}^9)^{[2^7]}$.
\begin{lem} \label{27rep}
    Suppose $\mathcal{R}(2,14,9;\sigma,a_1,\ldots,a_{14})\in\mathsf{S}^9_{2,14}$ is an interior $3$-spin relation. Then its $(2^7)$-isotypical projection is the zero relation.
\end{lem}
\begin{proof}
The inequality \eqref{eq:3spin_ineq} implies that $a_1+\cdots+a_{14}\le 3\cdot9-2-1 = 24$. Since $$0+0+1+1+3+3+4+4+6+6+7+7+9+9 = 60 > 24\, ,$$ this implies that some value occurs at least three times among the $a_i$. Therefore this $3$-spin relation is $\SSS_3$-invariant with respect to some subgroup $\SSS_3\subset \SSS_{14}$ given by permuting three of the $14$ markings, which implies that it has zero $(2^7)$-isotypical projection.
\end{proof}

\subsection{\texorpdfstring{$\QQQ$}{Q} is not a \texorpdfstring{$3$}{3}-spin relation} \label{prthm}

Let $$\mathsf{b3Spin}^*_{g,n}\subset \mathsf{3Spin}^*_{g,n}$$ be the subspace
obtained from boundary pushforwards \eqref{eq:boundary_rels}. 
Then \cref{27rep} says
$$(\mathsf{b3Spin}^9_{2,14})^{[2^7]} = (\mathsf{3Spin}^9_{2,14})^{[2^7]}\, .$$

\begin{proof}[Proof of \cref{qth_not3spin}]
    By construction, 
$\QQQ\in (\mathsf{S}^9_{2,14})^{[2^7]}$. Hence
\[ \QQQ\in \mathsf{3Spin}^9_{2,14} \implies \QQQ\in (\mathsf{3Spin}^9_{2,14})^{[2^7]} 
\xRightarrow{\text{\cref{27rep}}}
\QQQ\in 
(\mathsf{b3Spin}^9_{2,14})^{[2^7]} \implies \QQQ_\Delta=0 \in \HH_{16}(\Delta).\]
But 
by \cref{thm section 4}, we know $\QQQ_\Delta \neq 0 \in \HH_{16}(\Delta)$. So $\QQQ$ is not in $\mathsf{3Spin}^9_{2,14}$. \end{proof}
That $\QQQ_{\Delta_0}$ is a nontrivial $\Delta_0$-ghost is {\em not} sufficient for the proof of Theorem \ref{qth_not3spin}. The crucial property used is that $\QQQ_\Delta$ is a nontrivial $\Delta$-ghost.

\section{Does \texorpdfstring{$\QQQ$}{Q} hold in Chow?} \label{chow-quest}
\subsection{Beilinson's conjectures}
References for the following subsection are \cite{nekovar,jannsen} and the original \cite{beilinson}.
Let $X$ be an algebraic variety over a field $F$. There is a cycle map $\CH_k(X) \to \HH_{2k}^\BM(X,\Q(-k))$ and we denote its kernel by $\CH_k(X)_0$, the group of homologically trivial cycles. The \emph{Abel--Jacobi map}
\[ \CH_k(X)_0 \longrightarrow \mathrm{Ext}^1(\Q,\HH_{2k+1}^\BM(X,\Q(-k))) \]
may be constructed as follows. Let $\alpha \in \mathsf Z_k(X)_0$ be a homologically trivial cycle, and $Z \subset X$ the support of the cycle. Consider the following piece of the long exact sequence of Borel--Moore homology
\[ 0 \to \HH_{2k+1}^\BM(X,\Q(-k))\to \HH_{2k+1}^\BM(X \setminus Z,\Q(-k)) \to \HH_{2k}^\BM(Z,\Q(-k)) \to \HH_{2k}^\BM(X,\Q(-k)),\]
the first entry being $\HH_{2k+1}^\BM(Z,\Q(-k))=0$. The fundamental class of the cycle $\alpha$ is a homomorphism $\Q \longrightarrow \HH_{2k}^\BM(Z,\Q(-k))$. Since its composite into $\HH_{2k}^\BM(X,\Q(-k))$ is zero by hypothesis, we may form the pullback to obtain an extension:
\begin{equation}
    \label{AJ diagram}\begin{tikzcd}[column sep=small]
	0 \arrow[r]& \HH_{2k+1}^\BM(X,\Q(-k))\arrow[r]& \HH_{2k+1}^\BM(X \setminus Z,\Q(-k)) \arrow[r]& \HH_{2k}^\BM(Z,\Q(-k)) & \\
	0 \arrow[r]&\arrow[u,equal] \HH_{2k+1}^\BM(X,\Q(-k))\arrow[r]& \arrow[u]E_\alpha \arrow[r]& \Q \arrow[r]\arrow[u]& 0
\end{tikzcd}
\end{equation} 
The class of the extension $[E_\alpha]$ in $\mathrm{Ext}^1(\Q,\HH_{2k+1}^\BM(X,\Q(-k)))$ is the Abel--Jacobi image of $\alpha$; it depends only on the rational equivalence class of $\alpha$. We have been deliberately vague about in which category the $\mathrm{Ext}$ group is taken: if $F=\C$ we may take mixed Hodge structures; in general we may take cohomology to be \'etale cohomology over $F^{\mathrm{sep}}$, replace $\Q$ with $\Q_\ell$, and take the Ext in continuous $\ell$-adic $\mathrm{Gal}(F^{\mathrm{sep}}/F)$-representations, with $\ell \neq \operatorname{char} F$. 

Finally, what is most relevant here is that the construction is expected to work in a conjectural category $\mathsf{MMot}_F$ of \emph{mixed motives} over $F$. This makes the target group significantly more refined: in the motivic setting, one expects \emph{every} extension to arise from an algebraic cycle. In formulating the following conjecture we assume the existence of an abelian category $\mathsf{MMot}_F$ and motivic (Borel--Moore) homology groups $h_*^\BM(X)$ for any variety $X/F$, and the conjecture may be taken as a desideratum for the expected behavior of the category $\mathsf{MMot}_F$.

\begin{conjecture}[Beilinson]\label{BB1}
	Suppose $X$ is an algebraic variety over a number field $F$. Then
	\[ \CH_k(X)_0 \to \mathrm{Ext}^1_{\mathsf{MMot}_{F}}(\Q,h_{2k+1}^\BM(X,\Q(-k)))\]
	is an isomorphism. 
\end{conjecture}
If $X$ is smooth and proper, then $h_{2k+1}^\BM(X,\Q(-k))$ is a pure motive of weight $-1$. In this case, the above $\mathrm{Ext}$-group is expected to be deeply related to special values of $L$-functions.

\begin{conjecture}[Beilinson]\label{BB2}
	Let $M \in \mathsf{MMot}_F$ be a pure motive of weight $-1$ over a number field $F$. Then the dimension of $\mathrm{Ext}^1_{\mathsf{MMot}_{F}}(\Q,M)$ equals the order of vanishing of the $L$-function $L(M,s)$ of $M$ at the central point. 
\end{conjecture}

\begin{rem}
	There is also the more refined Bloch--Kato conjecture \cite{BK}, which gives a prediction for the leading coefficient of the $L$-function at the central point. It will not play a role in this paper.
\end{rem}

\begin{rem}
	The ``central point'' of the $L$-function depends on the normalizations chosen. Number theorists often impose that the functional equation of every $L$-function relates the value at $s$ and $1-s$, making the central point $s=1/2$. Arithmetic geometers often normalize $L(M,s)$ to have a functional equation relating $s$ and $w+1-s$, when $M$ is a motive of weight $w$; this would give the central point as $s=0$ in our weight $-1$ setting, and the familiar $s=1$ in the Birch--Swinnerton-Dyer conjecture where $M$ is $H^1$ of an elliptic curve.
\end{rem}

\begin{rem}
    \label{beilinson conjectures}
The formulation of \cref{BB1} is not quite the conventional one, which assumes $X$ smooth and projective. These hypotheses should however not be necessary with the above formulation. The general picture of Beilinson is that for any variety $X$ over any field there is a motivic spectral sequence calculating $\CH(X)$ and inducing a motivic filtration on Chow groups. When $X$ is smooth and projective the spectral sequence should degenerate, expressing the graded pieces of the filtration in terms of $\mathrm{Ext}$ groups of mixed motives. When $F$ is a number field, degeneration should hold without smoothness or properness, since we should have $\mathrm{Ext}^p=0$ for $p>1$ in this case.
\end{rem}

\subsection{Consequences in genus \texorpdfstring{$1$}{1}}

Beilinson's conjectures make interesting and highly nontrivial predictions about the Chow groups of $\Mbar_{1,n}$, since the cohomology in this case is quite well understood. Let us now explain this. We refer to \cite{semiclassical,getzler99,gorinov1,Petersen2014-GenusOne} for more background on how one computes the cohomology of $\Mbar_{1,n}$. The upshot is that all even cohomology is of Tate type, and the odd cohomology can all be described in terms of cusp forms for $\mathrm{SL}_2(\mathbb Z)$. More specifically, to every Hecke eigenform for $\mathrm{SL}_2(\mathbb Z)$ of weight $N$, there is attached a two-dimensional ``motive'' \cite{deligne69,scholl90}, which appears in the cohomology of $\Mbar_{1,n}$ for all $n \geq N-1$, and when it appears for the first time it occurs only in the middle cohomology group $\HH^{N-1}(\Mbar_{1,N-1})$. In particular, odd cohomology arises for the first time when $n=11$, arising from the Ramanujan cusp form $\Delta$. We write $S[N]$ for the motive attached to all cusp forms for $\mathrm{SL}_2(\mathbb Z)$ of weight $N$. Its Hodge numbers are $(N-1,0)$ and $(0,N-1)$. 

\vspace{5pt}
\emph{The case of at most 10 points.} When $n \leq 10$, we have just explained that the cohomology is entirely in even degree and of Tate type. The Hodge/Tate conjectures thus predict that the cycle map $\CH(\Mbar_{1,n}) \to \HH(\Mbar_{1,n})$ is surjective; moreover, Beilinson's conjectures on filtrations on Chow groups (\cref{beilinson conjectures}) predict that in fact $\CH(\Mbar_{1,n}) \to \HH(\Mbar_{1,n})$ should be an isomorphism over any base field. It is a theorem of Jannsen \cite[Theorem 3.6]{jannsen} that if the cycle map $\CH(X) \to \HH(X)$ is injective, for some smooth proper complex variety $X$, then it is also surjective; thus Beilinson's conjectures predict the converse. 
In this case, all these results are known unconditionally. Indeed, Belorousski \cite{belorousski} proved that $\R(\Mbar_{1,n})=\CH(\Mbar_{1,n})$ for $n \leq 10$, and Petersen \cite{Petersen2014-GenusOne} proved that $\R(\Mbar_{1,n}) \to \HH^{\text{even}}(\Mbar_{1,n})$ is an isomorphism for all $n$.

\vspace{5pt}
\emph{The case of 11--16 points.}
When $n=11$ one finds the first occurrence of odd cohomology and the motive $S[12]$ in $H^{11}(\Mbar_{1,11})$. Hence Jannsen's theorem implies the existence of nonzero homologically trivial cycles. Moreover, the Mumford--Ro\u{\i}tman theorem \cite{mumford,roitman} implies that the Chow group of zero-cycles is enormous.

\begin{thm}[Mumford--Ro\u{\i}tman] \label{mumfordroitman} Let $X$ be a smooth proper variety over an uncountable algebraically closed field of characteristic zero. Suppose that $\HH^0(X,\Omega^p) \neq 0$ for some $p>1$. Then for all $m$, the difference map
\[ \Sym^m(X) \times \Sym^m(X) \longrightarrow \CH_0(X)_0\]
is not surjective.     
\end{thm}

Indeed the theorem applies as in fact $\HH^{11,0}(\Mbar_{1,n}) \neq 0$ for all $n \geq 11$. So there must be a large supply of nontautological algebraic cycle classes over $\C$. Beilinson's conjectures predict that the situation over $\overline{\mathbb Q}$ is completely different. Indeed:
\begin{align*}
    \text{\cref{mumfordroitman}} & \implies && \CH_0((\Mbar_{1,11})_{\C})_0  &&\text{infinite dimensional in Mumford's sense} \\
    \text{\cref{BB1}} & \implies &&\CH_k((\Mbar_{1,11})_{\overline{ \mathbb Q}})_0 && \text{vanishes for $k \neq 5$.}
\end{align*}Moreover, for the group of $5$-cycles, the prediction is a more subtle statement: if $F$ is a number field, then according to \cref{BB2} the dimension of $\CH_5((\Mbar_{1,11})_{F})_0$ equals the order of vanishing at the center point, of the $L$-function of the motive $S[12]$, base changed to $F$. This may or may not vanish, depending on $F$. But in particular when $F=\Q$ we get the $L$-function of the cusp form $\Delta$ itself, which does not vanish at the center point, as discussed in \cref{nonvanishing}. When $n \leq 16$, the only further new motive appearing in the cohomology is $S[16]$, whose $L$-function likewise does not vanish at the center point. Beilinson's conjectures thus imply the following concrete prediction, which is widely open.

\begin{prediction}
    If $n \leq 16$, then $\R(\Mbar_{1,n})$ equals the subring of $\CH(\Mbar_{1,n})$ given by cycles defined over $\Q$.
\end{prediction}

Another intriguing prediction arising from the Beilinson conjectures (\cref{beilinson conjectures}) is the following: 
\begin{prediction}
    Let $F$ be any field. Then $\CH_k((\Mbar_{1,11})_F)_0$ is the summand of $\CH_k((\Mbar_{1,11})_F)$ which transforms according to the sign representation of $\SSS_{11}.$
\end{prediction}

Indeed, the non-Tate cohomology of $\Mbar_{1,11}$ is exactly the part that transforms by the sign representation \cite{faberconsani}.

\vspace{5pt}
\emph{The case of at least 17 points}. In the cohomology of $\Mbar_{1,17}$, we find the motive $S[18]$ attached to the unique cusp form for $\mathrm{SL}_2(\mathbb Z)$ of weight $18$. In this case, the $L$-function has a simple zero, and \cref{BB2} predicts instead the existence of homologically trivial cycles over $\Q$. More precisely we have the following.

\begin{prediction}
The space $\CH_8((\Mbar_{1,17})_\Q)_0$ is one-dimensional and transforms according to the sign representation of $\SSS_{17}$. There is no nonzero homologically trivial algebraic cycle defined over $\Q$ in any other degree. 
\end{prediction}

We note that such a class is necessarily nontautological, since the tautological ring maps injectively under the cycle map to cohomology in genus $1$ \cite{Petersen2014-GenusOne}. One part of this prediction --- the nontriviality of $\CH_8((\Mbar_{1,17})_\Q)_0$ --- is in fact known unconditionally, after work of S.-W. Zhang \cite{heegner}. Let us explain this.

Let $E$ be an elliptic curve with complex multiplication by the imaginary quadratic field $K$. The curve $E$ is defined over the Hilbert class field $H$. Zhang considers a certain cycle $\mathfrak z_E \in \CH_k(E^{2k})$ (he denotes it $S_{k+1}(E)$), and proves the following properties:
\begin{enumerate}
    \item Let $\mathcal E^n$ be the $n$-fold fibered power of the universal elliptic curve over $\M_{1,1}$, and let $\overline{\mathcal E^n}$ be the smooth compactification of it constructed by Deligne \cite{deligne69}. Then $\mathfrak z_E$, or rather its pushforward along $E^{2k} \hookrightarrow \mathcal E^{2k}$, is a homologically trivial cycle on $\overline{\mathcal E^{2k}}$. By construction $\mathfrak z_E$ is supported in a single smooth fiber over $\M_{1,1}$, and $\mathfrak z_E$ is defined over $H$. Let $\mathfrak z_K$ denote its trace from $H$ to $K$.
    \item Consider the submodule of $\CH_k(\overline{\mathcal E^{2k}})_0$ generated by $\mathfrak z_K$ under the action of the Hecke algebra. For a weight $2k+2$ eigenform $f$ for $\mathrm{SL}_2(\mathbb Z)$, let $\mathfrak z_f$ be the $f$-isotypic component in this module. For a smooth proper variety of dimension $2d+1$, let $\langle -,-\rangle$ denote the arithmetic height pairing of homologically trivial cycles $\CH_d(X)_0 \times \CH_d(X)_0 \to \mathbb R$. Then $ \langle \mathfrak z_f,\mathfrak z_f\rangle $ is an explicit nonzero constant times the derivative of $L(f/K,s)$ at the central point. Here $L(f/K,s) = L(f,s)L(f \otimes \chi_K,s)$.
\end{enumerate}
Zhang's construction is more general: he considers rather the submodule generated by $\mathfrak z_E$ under the Hecke algebra and $\mathrm{Gal}(H/K)$, which allows to replace $f$ by a twist of any character of $\mathrm{Gal}(H/K)$. Moreover, we have stated the result in level $1$, whereas Zhang works in arbitrary level, which he in fact takes to be at least $3$ to avoid stacky issues. We shall gloss over this point. 

Now Zhang's construction furnishes a large supply of nonzero homologically trivial cycles over quadratic fields. For example we may take $k=8$, and $f$ the unique weight $18$ and level $1$ cusp form, $K=H=\Q(i)$, and $E$ the curve $y^2=x^3+x$ of $j$-invariant $1728$. Then $L(f/K,s)$ has nonzero derivative at the center point, as one can look up in the \emph{$L$-functions and modular forms database}\footnote{Specifically, \href{https://www.lmfdb.org/ModularForm/GL2/Q/holomorphic/1/18/a/a/}{1.18.a.a} gives the analytic rank of $f$ as $1$, meaning that the zero at the central point is simple. It also lists \href{https://www.lmfdb.org/ModularForm/GL2/Q/holomorphic/16/18/a/b/}{16.18.a.b} as the quadratic twist of $f$ by the Dirichlet character $\chi_K$ (called 4.b in loc.~cit.), and records the analytic rank of the latter as zero.} \cite{lmfdb}. Hence $\mathfrak z_K$ is nonzero in this case. Moreover, Zhang shows \cite[2.4.1(b)]{heegner} that changing the embedding of $K$ to its complex conjugate changes the sign of $\mathfrak z_K$ by $(-1)^k$, which in our case is $+1$. So $\mathfrak z_K$ descends to a nonzero class defined over $\Q$.  

We should argue that Zhang's construction, which as stated only gives a nonzero class in $\CH_8(\overline{\mathcal E^{16}})_0$, produces an element of $\CH_8(\Mbar_{1,17})_0$. The argument is in fact completely general. Let $U$ be a smooth variety, and consider an algebraic cycle $\mathfrak z$ on $U$ all of whose irreducible components are proper. Let $X$ and $X'$ be smooth compactifications of $U$. Then $\mathfrak z$ is nontrivial in $\CH(X)$ if and only if it is nontrivial in $\CH(X')$. The reason is that the closure of $U$ in $X \times X'$ gives a correspondence between $X$ and $X'$, and pulling and pushing $\mathfrak z$ along this correspondence gives $\mathfrak z$ back since its support is contained in $U$. In our case, the support of $\mathfrak z_K$ is contained in $\mathcal E^{16}$, which we may instead compactify using Hassett's space of weighted pointed curves $\Mbar_{1,\epsilon^{17}}$ \cite{hassettweighted}. The latter receives a proper blowdown map $\Mbar_{1,17} \longrightarrow \Mbar_{1,\epsilon^{17}}$ under which pullback on Chow groups is injective. 

\begin{rem}
    Pulling this class back gives a large supply of nontautological classes over $\Q$ for any $n>17$, and higher weight modular forms produce even more. It seems an extremely hard problem to prove that all homologically trivial cycles on $\Mbar_{1,n}$ over $\Q$ arise in this way from modular forms.
\end{rem}

\begin{rem}
    Zhang's construction as explained above works verbatim already in weight $12$ to produce nontrivial elements of $\CH_5((\Mbar_{1,11})_{F})_0$ for imaginary quadratic fields $F$. In this case the cycle will however not descend to $\Q$ as complex conjugation instead acts by $-1$. This is of course what Beilinson's conjectures predict, too.
\end{rem}

\subsection{Consequences in genus \texorpdfstring{$2$}{2}}

The cohomology of $\Mbar_{2,n}$ is not as simple as that of $\Mbar_{1,n}$, but can be calculated by similar methods given the results on local systems from \cite{localsystemsA2}. Results have been conveniently tabulated by Bergstr\"om \cite{jonas-github}. Some of the arguments in the following are a natural extension of what was already discussed when $g=1$, and we permit ourselves to be brief in those cases.

In the following we shall say that a certain ``motive'' occurs in the cohomology of $\Mbar_{2,n}$. What we really mean is that our description gives the $\ell$-adic Galois representation in cohomology and the $\Q$-Hodge structure, using results of Gu--Lee--Prasanna \cite{GLP} for the latter.

\vspace{5pt}

\emph{The case of at most 9 points.} All cohomology is pure Tate, and one therefore expects $\CH(\Mbar_{2,n}) \to \HH(\Mbar_{2,n})$ to be an isomorphism in this range. This is in fact known; it is a theorem of Canning--Larson \cite{CL-jems}.

\vspace{5pt}

\emph{The case of 10--13 points.} One finds odd cohomology. The only non-Tate motive occurring in this range is $S[12]$, and it always appears with a positive Tate twist. The positive Tate twist means that the Mumford--Ro\u{\i}tman theorem on zero-cycles no longer applies. However, Jannsen's theorem implies that $\CH(\Mbar_{2,n}) \to \HH(\Mbar_{2,n})$ is not injective over $\C$, so \emph{some} nonzero homologically trivial cycle must exist in this range. Since the $L$-function associated to $S[12]$ is nonzero at the central point, Beilinson's conjectures predict the following:

    \begin{prediction}
        Let $n \leq 13.$ Then $\CH((\Mbar_{2,n})_\Q)_0=0$.
    \end{prediction}

Beilinson's conjectures make further intriguing predictions. For example, consider the correspondence 
\[ \Mbar_{1,11} \stackrel p\longleftarrow \Mbar_{1,12} \stackrel q \longrightarrow \Mbar_{2,10} \] given by forgetting the last point and gluing together the last two points, respectively. One can show that it induces an isomorphism $$q_*p^* :H^{11}(\Mbar_{1,11})(-1) \stackrel\sim\longrightarrow H^{13}(\Mbar_{2,10})$$
between the unique non-Tate cohomology groups of $\Mbar_{1,11}$ and $\Mbar_{2,10}$, up to Tate twist. One shows this by using the localization sequence for $\M_{2,10}^\ct \subset \Mbar_{2,10}$ and the computations of \cite{Petersen2016-GenusTwo} to see that $q_*$ is surjective on $H^{13}$. Compare also Canning--Larson--Payne \cite[Theorem 1.6]{STE}. The Beilinson conjectures (in the version that holds over any field, \cref{beilinson conjectures}) therefore predict the following.

\begin{prediction}Let $F$ be any field. 
    The map $\alpha \mapsto q_*p^*( \alpha )$ induces an isomorphism
    \[ \CH_k((\Mbar_{1,11})_F)_0 \longrightarrow \CH_{k+1}((\Mbar_{2,10})_F)_0.\]
\end{prediction}

\vspace{5pt}\emph{The case of at least 14 points.} In the cohomology of $\Mbar_{2,14}$, one finds two new motives, both in the middle degree $H^{17}$. Namely, they are $S[16](-1)$, occurring with multiplicity $2$ and transforming according to the representation $(1^{14})$ of $\SSS_{14}$, and $S[18]$, occurring with a single copy of the representation $(2^7)$. Mumford--Ro\u{\i}tman applies to the motive $S[18]$ to produce an enormous supply of homologically trivial zero-cycles. Assuming Beilinson's conjectures, one can be more precise---both summands give rise to nonzero homologically trivial cycles, but of rather different nature:
\begin{itemize}
    \item Working over $\Q$, the motive $S[16](-1)$ should not produce any nonzero homologically trivial cycle; the motive $S[18]$ should produce a single copy of the representation $(2^7)$ in $\CH_8((\Mbar_{2,14})_\Q)_0$.
    \item Working over $\C$, both motives should produce an uncountable-dimensional space of nonzero homologically trivial cycles, transforming entirely according to the representations $(1^{14})$ and $(2^7)$, respectively. In the case of $S[16](-1)$ all these cycles should be pushed forward from $\Mbar_{1,15} \times \Mbar_{1,1}$ and $\Mbar_{1,16}/\SSS_2$; in the case of $S[18]$ this should be far from true.
\end{itemize}

The last sentence is perhaps a bit vague. The point is that the motive $S[18]$ survives to the cohomology of $\M_{2,14}$, which is why it should contribute to $\CH(\M_{2,14})$ in the motivic spectral sequence; however, the motivic spectral sequence might be arbitrarily non-degenerate in this case, which complicates the analysis. But this is a digression.

In any case, returning to the matter at hand, the prediction is that $\CH_8((\Mbar_{2,14})_\Q)_0$ is a single copy of the representation $(2^7)$ of $\SSS_{14}$. Our expectation is that this group is generated by the homologically trivial cycle class $\QQQ$; in particular, we expect $\QQQ$ to be nonzero. Hence the motive $S[18]$ appears in the cohomology of both $\Mbar_{1,17}$ and $\Mbar_{2,14}$, and should produce nonzero homologically trivial $8$-cycles over $\Q$ on both spaces, but our expectation is that on $\Mbar_{1,17}$ we obtain a non-tautological class and on $\Mbar_{2,14}$ we obtain a tautological class. That the class on $\Mbar_{1,17}$ should be nontautological could have been seen a priori from the length bounds on tautological classes of Faber--Pandharipande \cite{FP3}; no class in $\R^9(\Mbar_{1,16})$ could transform according to the sign representation. But it is consistent with the length bound for a class in $\R^9(\Mbar_{2,14})$ to transform according to the representation $(2^7)$.

The discussion above only indicates that $\QQQ$ \emph{could} be nonzero in Chow: Beilinson's conjectures say that there should be a nontrivial homologically trivial cycle in $\CH_8(\Mbar_{2,14})$ over $\Q$, and it is consistent with everything we know that that cycle could be a tautological class. But in fact a conjecture of Harder implies that $\QQQ$ \emph{should} be nonzero in Chow. We now want to explain this implication. 

\subsection{A conjecture of Harder}

Let $V$ be the $\SSS_{14}$-representation generated by $\QQQ$ inside the strata algebra. By our construction of $\QQQ$, the space $V$ is isomorphic to a single copy of the representation $\sigma_{(2^7)}$. Since every cycle in $V$ is homologically trivial, there is the Abel--Jacobi map
\begin{equation} \label{equivariant AJ}
    V \longrightarrow \mathrm{Ext}^1(\Q,\HH_{17}(\Mbar_{2,14},\Q(-8))).
\end{equation} 
The map \eqref{equivariant AJ} is $\SSS_{14}$-equivariant. From the tables of Bergstr\"om \cite{jonas-github}, we may read off the decomposition of $\HH^{17}(\Mbar_{2,14})$ into irreducible representations of $\SSS_{14}$: one finds
\[ \HH^{17}(\Mbar_{2,14}) \cong S[18] \otimes \sigma_{(2^7)} \oplus S[16](-1) \otimes \sigma_{(1^{14})}^{\oplus 2} \oplus S[12](-3) \otimes W,\]
for a certain reducible representation $W$ whose formula we omit. The key point is that the $(2^7)$-isotypical component is only the first term $S[18] \otimes \sigma_{(2^7)}$. Now Schur's lemma gives
\begin{align*}
    \mathrm{Hom}_{\SSS_{14}}(V,\mathrm{Ext}^1(\Q,\HH_{17}(\Mbar_{2,14},\Q(-8))))  & \cong \mathrm{Hom}_{\SSS_{14}}(V,\mathrm{Ext}^1(\Q,S[18](9)) \otimes \sigma_{(2^7)}) \\
    & \cong \mathrm{Ext}^1(\Q,S[18](9)).
\end{align*} We will call the resulting element of $\mathrm{Ext}^1(\Q,S[18](9))$ the \emph{Abel--Jacobi class of $\QQQ$}. The class is only defined up to multiplication by a nonzero rational number, since the above isomorphisms are not canonical. 

An alternative description of the same class is as follows. Let $G= \SSS_2 \wr \SSS_7 = (\SSS_2)^7 \rtimes \SSS_7$. The cohomology of $\Mbar_{2,14}/G$ simplifies significantly compared to $\Mbar_{2,14}$: all its cohomology is in even degree and tautological, except $\HH^{17}(\Mbar_{2,14}/G)$, which consists of a single copy of the motive $S[18]$. Indeed, this is a direct computation from \cite{jonas-github}, using that the $\SSS_{14}$-representation $\sigma_{(2^7)}$ has a $1$-dimensional $G$-invariant subspace, and no other $\SSS_{14}$-representation occurring among non-Tate terms has any $G$-invariants. Moreover, the one-dimensional space $V^G$ has a preferred generator, namely the $G$-average $\Qinv := \smash{{|G|^{-1}}\sum_{g \in G} g\QQQ}$. It is nonzero (as an element of the strata algebra) by the following \cref{Qinv not zero} and descends to a cycle on $\Mbar_{2,14}/G$, which we shall continue to denote $\Qinv$. The Abel--Jacobi class of $\QQQ$ is then the same as the image of $\Qinv$ under the Abel--Jacobi map
    \begin{equation*}\label{AJ-ext}
        \CH_8(\Mbar_{2,14}/G)_0 \longrightarrow \mathrm{Ext}^1(\Q,\HH_{17}(\Mbar_{2,14}/G,\Q(-8))) \cong \mathrm{Ext}^1(\Q,S[18](9)).
    \end{equation*} 

\begin{lem}\label{Qinv not zero}
    The class $\Qinv$ is not zero.
\end{lem}

\begin{proof}
    Let $a,b\in\Q[\SSS_{14}]$ be the idempotents obtained by  symmetrizing over $(\SSS_2)^7$, and antisymmetrizing over $\SSS_7\times\SSS_7$, respectively. Since the diagonal $\SSS_7$ already acts trivially on $\QQQ$, we have $a\QQQ = \Qinv$. By the theory of Young symmetrizers, $ba$ is idempotent up to a nonzero scalar. 
Since $\QQQ$ is a nonzero element in the image of $ba$, we obtain
\[
0 \neq ba\QQQ = b\Qinv. \qedhere
\]
\end{proof}

On the other hand, let us consider the cohomology of symplectic local systems on $\mathcal A_2$. These cohomology groups were calculated up to semisimplification in \cite{localsystemsA2}, building on earlier work of several people, in particular Harder \cite{harderA2}. Knowing the answer only up to semisimplification leaves the possibility of nontrivial extensions, and Harder has conjectured that extensions should indeed occur between the part of cohomology arising from Saito--Kurokawa lifting of elliptic modular forms, and a specific piece of the Eisenstein cohomology. Let us consider the first nontrivial case of the Saito--Kurokawa lift $\chi_{10}$ and the local system $\mathbb V_{7,7}$. From \cite{localsystemsA2} we see that \[ \gr^W_k \HH^3_c(\mathcal A_2,\mathbb V_{7,7}) \cong \begin{cases}
S[18] & k=17 \\
\Q(-8) & k= 16 \\
\Q & k=0 \\
0 & k \notin \{0,16,17\}.
    \end{cases} \]
    Dualizing, we obtain an extension
    \begin{equation}\label{LS-ext}
        0 \longrightarrow S[18] \longrightarrow W_{18} \HH^{3}(\mathcal A_2,\mathbb V_{7,7}) \longrightarrow \Q(-9) \longrightarrow 0.
    \end{equation} 
    Harder has conjectured that this extension is nontrivial \cite[pp.~81--82]{harderbook}. We have been somewhat vague about in what category we are taking this $\mathrm{Ext}$-group. Certainly Harder was thinking about a conjectural category of mixed motives over $\Q$ or $\mathbb Z$. But we may also consider the extension in mixed Hodge structures or Galois representations. Either way, we have the following theorem.

\begin{thm}Under the isomorphism $\mathrm{Ext}^1(\Q,S[18](9)) = \mathrm{Ext}^1(\Q(-9),S[18])$, the Abel--Jacobi class of $\QQQ$ equals the class of the extension \eqref{LS-ext}, up to multiplication by a nonzero rational number. 
\end{thm}

\begin{proof}
Consider the diagram \eqref{AJ diagram} defining the Abel--Jacobi class of $\Qinv$, i.e.~
\[ \begin{tikzcd}[column sep=small]
	0 \arrow[r] & \HH_{17}(X,\Q(-8))\arrow[r]& \HH_{17}^\BM(X \setminus Z,\Q(-8)) \arrow[r]& \HH_{16}(Z,\Q(-8)) & \\
    0 \arrow[r]&\arrow[u,equal] \HH_{17}(X,\Q(-8))\arrow[r]& \arrow[u]E_{\Qinv} \arrow[r]& \Q \arrow[r]\arrow[u]& 0
\end{tikzcd}\]
with $X=\Mbar_{2,14}/G$ and $Z=\operatorname{supp} \Qinv$. Let $Z' = \Delta_0/G \subset X$ , so that $X \setminus Z' = \M_{2,14}^\ct/G$.  We may then augment this diagram with another exact row above it, as follows:
\[ \begin{tikzcd}[column sep=small]
	\HH_{17}(Z',\Q(-8)) \arrow[r]& \HH_{17}(X,\Q(-8))\arrow[r]& \HH_{17}^\BM(\M_{2,14}^\ct/G,\Q(-8)) \arrow[r]& \HH_{16}(Z',\Q(-8)) & \\
	0 \arrow[u]\arrow[r] & \arrow[u,equal]\HH_{17}(X,\Q(-8))\arrow[r]& \arrow[u]\HH_{17}^\BM(X \setminus Z,\Q(-8)) \arrow[r]& \arrow[u]\HH_{16}(Z,\Q(-8)) \end{tikzcd}\]
    Now the map $\HH_{17}(Z',\Q(-8)) \longrightarrow \HH_{17}(X,\Q(-8))$ is zero. Indeed the codomain is given by the motive $S[18](9)$ which does not occur in $\HH_{17}(Z',\Q(-8))$. One sees this by writing down a Mayer--Vietoris spectral sequence calculating the homology of $Z'$ (see the earlier discussion in connection with \cref{propdelq}) and noting that all terms on the $E_1$-page involve spaces $\Mbar_{1,n}$ with $n \leq 16$ and $\Mbar_{0,n}$ with $n \leq 18$; the motive $S[18]$ does not appear in the homology of $\Mbar_{1,n}$ until $n=17$. Moreover, the composite
    \[ \Q \longrightarrow \HH_{16}(Z,\Q(-8)) \longrightarrow \HH_{16}(Z',\Q(-8))\]
    is nonzero, precisely because $\QQQ$ is homologically nontrivial on $\Delta_0$ (\cref{propdelq}), and hence also $\Qinv$. This means that $E_{\Qinv}$ is pulled back from an extension occurring in the homology group $$\HH_{17}^\BM(\M_{2,14}^\ct/G,\Q(-8)) \stackrel{\mathrm{PD}}\cong \HH^{17}(\M_{2,14}^\ct/G,\Q(9)).$$
    Furthermore, as already discussed, the cohomology groups $\HH^*(\M_{2,n}^\ct)$ may be computed $\SSS_n$-equivariantly and up to semisimplification as in \cite{Petersen2016-GenusTwo}. We have a direct summand
    \[ \HH^{*+14}(\M_{2,14}^\ct,\Q) \supset H^*(\mathcal A_2,\mathbb V_{7,7}) \otimes \sigma_{(2^7)}\]
    arising from the Decomposition Theorem, and hence also
    \[ \HH^{*+14}(\M_{2,14}^\ct/G,\Q) \supset H^*(\mathcal A_2,\mathbb V_{7,7}).\]
    We have already noted that $\gr^W_{18}\HH^{17}(\M_{2,14}^\ct,\Q) \cong \Q(-9) \otimes \sigma_{(2^7)}$, and the class lies in precisely this summand. Thus  
    \[ \gr^W_{0} \HH^{17}(\M_{2,14}^\ct/G,\Q(9)) \cong \Q,\]
    and since this group is one-dimensional, the image of a generator of it in $\HH_{16}(Z,\Q(-8))$ must be proportional to the class of $\Qinv$ (using that $\Qinv$ is homologically trivial on $X$ and exactness of the rows above). Hence the pullback extension $E_{\Qinv}$ must be specifically the pullback of the extension occurring in $\HH^3(\mathcal A_2,\mathbb V_{7,7})$.
\end{proof}

\begin{cor}
    If Harder's extension \eqref{LS-ext} is nontrivial, then $\QQQ$ is nonzero in Chow.
\end{cor}

\begin{rem}
    The \emph{Manin--Drinfeld theorem} \cite{manin,drinfeld} says that the difference between two cusps is always a torsion point of the Jacobian on any modular curve. Deligne \cite[Remarque 7.5]{D79} (see also Elkik \cite{elkik}) has explained how to derive the theorem by ``pure thought''. Namely, if $Y$ is an open modular curve, then the theorem is equivalent to saying that the mixed Hodge structure on $H^1(Y,\Q)$ is split (i.e.~the direct sum of its weight $1$ and weight $2$ subquotients); this follows from a standard description of extensions of mixed Hodge structures. But the Hecke algebra acts on this group, with eigenvalues of different absolute value on the weight $1$ and $2$ subquotients, and this splits the extension. Similarly the mixed Hodge structure on cohomology with local system coefficients $H^1(Y,\mathbb V_k)$ is likewise split. In Harder's work, he says that the ``Manin--Drinfeld principle is satisfied'' whenever motivic extensions are forced to split in this way by the Hecke algebra, in the cohomology of more general Shimura varieties. What typically happens is that the Manin--Drinfeld principle is indeed satisfied, but the cohomology group $\HH^3_c(\mathcal A_2,\mathbb V_{7,7})$ is a striking exception. The subquotient $\Q(0)$ is split off by the Hecke algebra, but the Hecke algebra acts in the same way on the weight $16$ and $17$ subquotients --- heuristically, this is because both subquotients $\Q(-8)$ and $S[18]$ arise from the same elliptic Hecke eigenform, the former through Eisenstein cohomology and the latter through the Saito--Kurokawa lifting. Hence a nontrivial extension is possible. The fact that the local system $\mathbb V_{7,7}$ has singular highest weight is important, since this permits the existence of CAP automorphic representations (such as the Saito--Kurokawa lifting) contributing to cohomology. This offers some motivation for Harder's conjecture: we are in an exceptional case where both Beilinson's conjectures and the action of the Hecke algebra  permit the possibility of a nontrivial extension, and then this is a sign from the universe that the extension in question is actually there in the geometry. 
\end{rem}

\begin{rem}
    As stated, Harder conjectured the nontriviality of the extension \eqref{LS-ext} already in his 1993 book \cite{harderbook}. The problem is still open. Harder returned to this question several times, in particular in his final paper \cite[pp.~161--162]{harderdenominators}. There he considers a certain biextension $E$ which in our case would have 
    \[ \gr^W_k E = \begin{cases}
        \Q(-8) & k = 16 \\
        S[18] & k = 17 \\
        \Q(-9) & k = 18
    \end{cases}\]
    and no other nontrivial graded pieces. The simple extensions in this biextension are \eqref{LS-ext} and its dual. Harder sketches a computation of the global height pairing of this biextension but declines to state it as a theorem, because ``some computations still need to be checked''. The formula is in terms of the derivative of the $L$-function at the center point, which in particular gives a nonzero height when the central zero is simple. Nonvanishing of the global height pairing would imply nontriviality of \eqref{LS-ext} and nonvanishing of $\QQQ$. A natural goal is therefore to make these computations rigorous in our setting and to identify the resulting pairing explicitly with the arithmetic height pairing of $\QQQ$ with itself. Note that more refined conjectures of Beilinson \cite[Conjecture 6.5(3)]{nekovar} predict that $\langle \QQQ,\QQQ\rangle$ should be a $\Q$-multiple of the derivative of the $L$-function at the center point, divided by a suitable period. Zhang's formula for $\langle \mathfrak z_f,\mathfrak z_f\rangle$ provides an analogue in the genus $1$ setting. We leave these questions for future work.
\end{rem}

\section{Higher genus}\label{sec:highergenus}


\subsection{Properties of \texorpdfstring{$\ThreeSpin$}{3-Spin}}
As before, let $\mathsf{S}^*_{g,n}$ be the strata algebra of $\Mb_{g,n}$. For a stable graph $\Gamma$ of genus $g$ with $n$ legs, there is an associated strata algebra gluing pushforward
\[\xi_{\Gamma*}:\mathsf{S}^*_{\Gamma}=\bigotimes_{v\in \Gamma} \mathsf{S}^*_{g(v),n(v)}\rightarrow \mathsf{S}^*_{g,n} \]
and pullback
\[
\xi_{\Gamma}^*:\mathsf{S}^*_{g,n}\rightarrow \bigotimes_{v\in \Gamma} \mathsf{S}^*_{g(v),n(v)}\, .
\]
Similarly, there are pushforward and pullback morphisms associated to the map 
$$\pi:\Mb_{g,n+1}\rightarrow \Mb_{g,n}$$ forgetting a marked point:
\[
\pi_*: \mathsf{S}^*_{g,n+1}\rightarrow \mathsf{S}^*_{g,n} \qquad \text{ and }\qquad  \pi^*:\mathsf{S}^*_{g,n}\rightarrow \mathsf{S}^*_{g,n+1}\, .
\]
There are natural quotient maps
\[
\mathsf{S}^*_{g,n}\rightarrow \mathsf{R}^*(\Mb_{g,n})\rightarrow \mathsf{RH}^*(\Mb_{g,n})\, .
\]
We have an integration map
\[
\int:\mathsf{S}^*_{g,n}\rightarrow \qq
\]
defined by mapping first to $\RH^*(\Mb_{g,n})$ and then integrating over $\Mb_{g,n}$.
For each vertex $v$ in a stable graph $\Gamma$, there is a projection map
\[
\mathsf{pr}_{v*}: \mathsf{S}^*_{\Gamma}\rightarrow \mathsf{S}^*_{g(v),n(v)}
\]
defined by integration of the tensor factors corresponding to vertices $w\neq v$.

Pixton proved that the collection of ideals $\{\ThreeSpin_{g,n}^*\}$ is closed under the operations $\pi_*$, $\pi^*$, $\xi_{\Gamma*}$, and $\xi_{\Gamma}^*$, see \cite[Section 5.3]{Pix}. For $\xi_{\Gamma}^*$, closure means $\xi_{\Gamma}^*(\ThreeSpin^*_{g,n})\subset \ThreeSpin^*_{\Gamma}$, where $\ThreeSpin^*_{\Gamma}$ is the ideal generated by the ideals $\ThreeSpin^*_{g(v),n(v)}$ for vertices $v\in \Gamma$.

\begin{lem}\label{lem:onefactor}
   Let $\Gamma$ be a stable graph of genus $g$ with $n$ legs. If $\alpha \in \ThreeSpin_{g,n}^*$, then $$\mathsf{pr}_{v*} \xi_\Gamma^*(\alpha)\in \ThreeSpin_{g(v),n(v)}^*$$ for every $v\in \Gamma$.
\end{lem}
\begin{proof}
    We have \[\xi_{\Gamma}^*(\alpha)\in \sum_{u\in \Gamma} \big(\ThreeSpin_{g(u),n(u)}^*\otimes\bigotimes_{w\neq u} \mathsf{S}^*_{g(w),n(w)}\big)\,.\]
    For any
    \[
    \beta\in \ThreeSpin_{g(u),n(u)}^*\otimes\bigotimes_{w\neq u} \mathsf{S}^*_{g(w),n(w)}
    \]
    with $u\neq v$, $\mathsf{pr}_{v*}(\beta)=0$ because the $3$-spin relations hold in cohomology and integration factors through cohomology. For any 
    \[
    \gamma \in \ThreeSpin_{g(v),n(v)}^*\otimes\bigotimes_{w\neq v} \mathsf{S}^*_{g(w),n(w)}\,,
    \]
    $\mathsf{pr}_{v*}(\gamma)\in \ThreeSpin_{g(v),n(v)}$ by definition.
\end{proof}


\begin{lem}\label{lem:addmarkings}
    Let $\gamma\in \mathsf{S}^*_{g,n}$. If $\gamma\notin \ThreeSpin^*_{g,n}$, then $\pi^*\gamma \notin \ThreeSpin^*_{g,n+1}$.
\end{lem}
\begin{proof}
    If $\pi^*\gamma \in \ThreeSpin^*_{g,n+1}$, then so is $\psi_{n+1}\cdot \pi^*\gamma$. After pushing forward and applying the projection formula, we obtain $\gamma\in \ThreeSpin^*_{g,n}$.
\end{proof}

Let $g\geq 2$. Let $\mathsf{E}_{g}$ be the stable graph of genus $g+1$ (with no markings) corresponding to the locus of elliptic tails: $\mathsf{E}_{g}$ has two vertices $\{v,v'\}$ connected by a single edge, the vertex $v$ is of genus $g$ and the vertex $v'$ is of genus $1$. 
Let $$\mathsf{pr}_v:\Mbar_{\mathsf{E}_g}\rightarrow \Mbar_{g,1}$$
be projection onto the factor corresponding to the vertex $v$.
Let $$\xi_{\mathsf{E}_{g}}:\Mb_{\mathsf{E}_g}= \Mb_{g,1}\times \Mb_{1,1}\rightarrow \Mb_{g+1}$$ be the associated gluing map. 
\vspace{8pt}

\noindent $\bullet$ For $\gamma\in \mathsf{S}^k_g$, we define the genus raising operator
\[
\R_g(\gamma)=\xi_{\mathsf{E}_{g*}}(\pi^*\gamma\otimes 1)\in \mathsf{S}^{k+1}_{g+1}\,.
\]

\vspace{8pt}

\noindent $\bullet$ For $\gamma\in \mathsf{S}^k_{g}$, we define the genus 
lowering operator
\[
\L_g(\gamma)=24\pi_*\mathsf{pr}_{v*}\big((1\otimes \psi)\cdot \xi_{\mathsf{E}_{g-1}}^*\gamma\big)\in \mathsf{S}^{k-1}_{g-1}\,.
\]
\noindent In the formulas above, $\pi$ denotes the forgetful map
$$\pi:\Mbar_{h,1}\rightarrow \Mbar_{h}$$ and
$\psi\in\mathsf{S}^*_{1,1}$.


\begin{lem}\label{lem:commutator}    
    For $\gamma \in \mathsf{S}_g^*$,
    we have the commutation relation
    \[
    \L_{g+1}(\R_{g}(\gamma)) = \R_{g-1}(\L_{g}(\gamma))+(2-2g)\gamma \, .
    \]
\end{lem}
\begin{proof}
    We have
    \[
    \L_{g+1}(\R_{g}(\gamma)) = 24\pi_*\pr_{v*}\big((1\otimes \psi)\cdot \xi_{\mathsf{E}_g}^*\xi_{\mathsf{E}_{g*}}(\pi^*\gamma\otimes 1)\big)\,.
    \]
    By definition of the product in the strata algebra via excess intersection theory, we have
    \begin{equation}\label{eq:excess0}
    \xi_{\mathsf{E}_g}^*\xi_{\mathsf{E}_{g*}}(\pi^*\gamma\otimes 1) = \sum_{(\Gamma,f,g)\in \mathfrak{G}_{\mathsf{E},\mathsf{E}}} \xi_{f*}\bigg (\xi_g^*(\pi^*\gamma\otimes 1)\cdot \prod_{(h,h')\in \mathrm{Im}(f)\cap \mathrm{Im}(g)} (-\psi_h - \psi_{h'}) \bigg )\,.
    \end{equation}
Here, $\mathfrak{G}_{\mathsf{E}_g,\mathsf{E}_g}$ is the set of isomorphism classes of generic $(\mathsf{E}_g,\mathsf{E}_g)$ structures $(\Gamma,f,g)$. There are only two isomorphism classes of generic $(\mathsf{E}_g,\mathsf{E}_g)$ structures. 

The first isomorphism class is $\Gamma=\mathsf{E}_g$ with both $f$ and $g$ the identity morphisms. 
The associated contribution to \eqref{eq:excess0} is 
\[
\mathsf{Cont}_{\mathsf{E}_g} = -\psi\cdot\pi^*\gamma\otimes 1 - \pi^*\gamma\otimes \psi\,,
\]
where $\psi$ on each side corresponds to the half-edge constituents of the edge of $\mathsf{E}_g$. Therefore, we have 
\[
24\pi_*\pr_{v*}\big((1\otimes\psi)\cdot\mathsf{Cont}_{\mathsf{E}_g}\big)=(2-2g)\gamma\,.
\]

The second isomorphism class has underlying graph $\Gamma$ with three vertices, two of genus $1$ and one of genus $g-1$. We find
$$\mathsf{Cont}_\Gamma = \R_{g-1}(\L_g(\gamma))\, .$$
We leave the straightforward check for the reader. \end{proof}

\begin{lem}\label{lem:addgenus}
    Let $\gamma\in \mathsf{CohRel}^*_{g}$. If $\gamma \notin \ThreeSpin^*_{g}$, then $\R_{g}(\gamma)\notin \ThreeSpin^*_{g+1}$.
\end{lem}
\begin{proof}
Suppose $\R_g(\gamma)\in \ThreeSpin^*_{g+1}$. The lowering and raising operators satisfy
\[
\mathsf{L}_g(\ThreeSpin_{g}^*)\subset \ThreeSpin_{g-1}^* \qquad \text{ and } \qquad \mathsf{R}_{g}(\ThreeSpin_g^*)\subset \ThreeSpin_{g+1}^*
\]
by Lemma \ref{lem:onefactor} and the basic properties of $\ThreeSpin$. 
Iteratively applying Lemma \ref{lem:commutator}, we obtain
\begin{equation}\label{eq:iteration}
\mathsf{R}_{g-1}\circ\dots\circ\R_{2}\circ \mathsf{L}_3\circ\dots\circ \mathsf{L}_g(\gamma) = \gamma\prod_{k=0}^{g-3} (k+1)(2g-k-2)\in \mathsf{S}^*_{g}/\ThreeSpin^*_{g}\,.
\end{equation}
Because $\gamma\in \mathsf{CohRel}^*_{g}$, and hence $\mathsf{L}_3\circ\dots\circ \mathsf{L}_g(\gamma)\in \mathsf{CohRel}^*_2= \ThreeSpin^*_2$, we have
\[
\mathsf{L}_3\circ\dots\circ \mathsf{L}_g(\gamma) = 0 \in \mathsf{S}_2^*/\ThreeSpin^*_{2}\,,
\]
and so $\gamma\in \ThreeSpin^*_g$ by \eqref{eq:iteration}. \end{proof}

\subsection{Proof of Theorem \ref{thm:highergenus}}

Let $0\leq k\leq 7$. Let $\mathsf{A}_{2+2k,14-2k}$ be the stable graph of genus $2+2k$ with $14-2k$ legs such that
\begin{enumerate}
    \item[(i)] $\mathsf{A}_{2+2k,14-2k}$ has one vertex of genus $2$ and $k$ vertices of genus $1$,
    \item[(ii)] the genus $2$ vertex has $14-2k$ legs,
    \item[(iii)] each genus $1$ vertex is attached to the genus $2$ vertex by $2$ edges, and there are no other edges,
\end{enumerate}
as depicted in Figure \ref{fig:sunflower}. 
Recall that by Lemma \ref{Qinv not zero}, $\Qinv$ is a nonzero $(\SSS_2\wr \SSS_7)$-invariant element of the $\SSS_{14}$-representation $\mathcal{Q}$ of cohomological relations outside of $\ThreeSpin^9_{2,14}$. For $0\leq k \leq 7$, we set \[
\alpha_k=\xi_{\mathsf{A}_{2+2k,14-2k}*}(\mathsf{Q}_{\mathsf{inv}}\otimes1^{\otimes k})\in \mathsf{S}^{*}(\Mb_{2+2k,14-2k})\,.
\]
\begin{figure}
\centering
\begin{tikzpicture}[
  moregraphs/center,draw=black,text=black,
  sunflower vertex/.style={circle,draw=black,fill=white,line width=0.85pt,
    minimum size=11mm,inner sep=0pt,align=center,font=\small},
  sunflower edge/.style={draw=black,line width=0.85pt},
  sunflower leg/.style={draw=black!50,line width=0.55pt},
  sunflower label/.style={text=black,font=\small,inner sep=3pt}
]
  \node[sunflower vertex,minimum size=18.5mm,font=\normalsize]
    (sunflower-center) at (0,0) {$2$};
  \node[sunflower vertex] (sunflower-top) at (0,3.05) {$1$};
  \node[sunflower vertex] (sunflower-upper-right) at (2.39,1.9) {$1$};
  \node[sunflower vertex] (sunflower-lower-right) at (2.98,-0.69) {$1$};
  \foreach \petal in {top,upper-right,lower-right}{
    \draw[sunflower edge] (sunflower-center) to[bend left=18] (sunflower-\petal);
    \draw[sunflower edge] (sunflower-center) to[bend right=18] (sunflower-\petal);
  }
  \foreach \mark/\x/\y/\anchor in {
    14/-2/1.96/south east,7/-2.35/1.52/south east,
    13/-2.78/-0.35/east,6/-2.66/-0.9/north east,
    12/-1.46/-2.4/north east,5/-0.95/-2.64/north,
    11/0.95/-2.64/north,4/1.46/-2.4/north west}{
    \draw[sunflower leg] (sunflower-center) -- (\x,\y);
  }
\end{tikzpicture}
\caption{The graph $\mathsf{A}_{8,8}$}\label{fig:sunflower}
\end{figure}
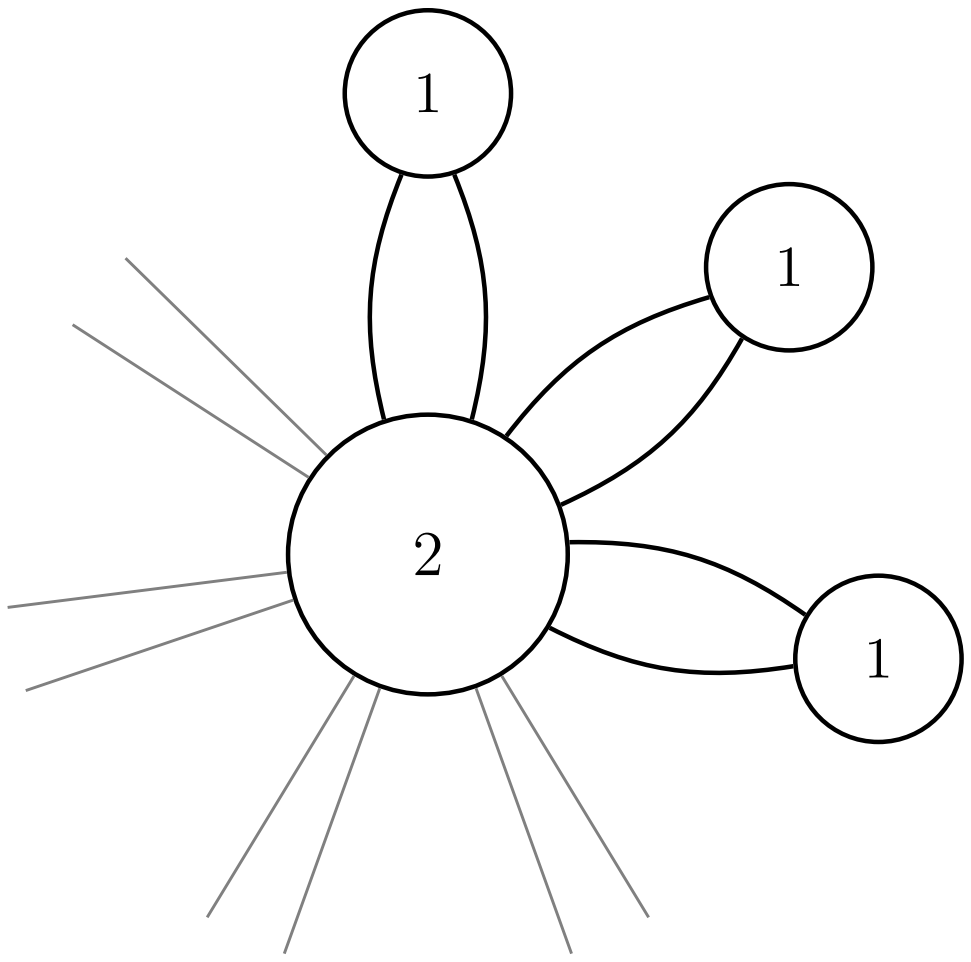
\begin{lem}\label{lem:sunflower}
    We have $\alpha_k \notin \ThreeSpin_{2+2k, 14-2k}^*$.
\end{lem}
\begin{proof}
     It suffices to prove the statement for $\alpha_7$. Set $\mathsf{A} = \mathsf{A}_{16,0}$.     As in the proof of Lemma \ref{lem:commutator}, we have 
    \begin{equation}\label{eq:excess}
    \xi_{\mathsf{A}}^*\alpha_k = \sum_{(\Gamma,f,g)\in \mathfrak{G}_{\mathsf{A},\mathsf{A}}} \xi_{f*}\bigg (\xi_g^*(\mathsf{Q}_{\mathsf{inv}}\otimes 1^{\otimes 7})\cdot \prod_{(h,h')\in \mathrm{Im}(f)\cap \mathrm{Im}(g)} (-\psi_h - \psi_{h'}) \bigg )\,.
    \end{equation}

    First, we consider generic $(\mathsf{A},\mathsf{A})$ structures $(\Gamma,f,g)$ with $\Gamma = \mathsf{A}$. Up to isomorphism, we may assume $f$ is the identity and $g$ is any automorphism of $\mathsf{A}$. Therefore, there are $|\Aut(\mathsf{A})|$ 
    generic $(\mathsf{A},\mathsf{A})$ structures on $\mathsf{A}$. We have $\xi_g^*(\mathsf{Q}_{\mathsf{inv}}\otimes 1^{\otimes 7})=\mathsf{Q}_{\mathsf{inv}}\otimes 1^{\otimes 7}$ for all $g$ because $\mathsf{Q}_{\mathsf{inv}}$ is $\Aut(\mathsf{A})$ invariant. The contribution to equation \eqref{eq:excess} from the generic $(\mathsf{A},\mathsf{A})$ structures  of the form $(\mathsf{A},f,g)$ is
    \[
    |\Aut(\mathsf{A})| (\mathsf{Q}_{\mathsf{inv}}\otimes 1^{\otimes7}) \prod_{(h,h')\in E(\mathsf{A})} (-\psi_h-\psi_{h'})\,.
    \]
    Let $w$ be the genus $2$ vertex of $\mathsf{A}$, and $\mathsf{pr}_w:\Mbar_{\mathsf{A}}\rightarrow \Mbar_{2,14}$ the projection map.
    We have
    \[
    |\Aut(\mathsf{A})|\pr_{w*}\big((\mathsf{Q}_{\mathsf{inv}}\otimes 1^{\otimes7}) \prod_{(h,h')\in E(\mathsf{A})} (-\psi_h-\psi_{h'})\big)=
    |\Aut(\mathsf{A})|\mathsf{pr}_{w*}(\mathsf{Q}_{\mathsf{inv}}\otimes (\psi_1\psi_2)^{\otimes 7})= c\cdot\mathsf{Q}_{\mathsf{inv}}
    \]
    for some nonzero constant $c$.

    Next, we consider generic $(\mathsf{A},\mathsf{A})$ structures whose underlying stable graph $\Gamma$ is not isomorphic to $\mathsf{A}$. For these, every vertex of $\Gamma$ has genus $0$ or $1$ by genericity. We have $\xi_g^*(\mathsf{Q}_{\mathsf{inv}}\otimes 1^{\otimes 7})\in \ThreeSpin_{\Gamma}^*$ because the $3$-spin relations are complete in genus $0$ and $1$ and $\mathsf{Q}_{\mathsf{inv}}\in \mathsf{CohRel}^*_{2,14}$. 
    Let $\Gamma'$ be the stable graph of genus $2$ with $14$ markings that is the preimage  $f^{-1}(v)$ in $\Gamma$. Consider the projection map
    \[
    \mathsf{pr}_{\Gamma'}: \Mbar_{\Gamma}\rightarrow \Mbar_{\Gamma'}\,.
    \]
    Then for the same reason as before,
     \[
     \mathsf{pr}_{\Gamma'_*}\big(\xi_g^*(\mathsf{Q}_{\mathsf{inv}}\otimes 1^{\otimes 7})\cdot \prod_{(h,h')\in \mathrm{Im}(f)\cap \mathrm{Im}(g)} (-\psi_h - \psi_{h'}) \big )\in \ThreeSpin_{\Gamma'}^*\,,\]
     and hence
     \[
     \xi_{\Gamma'*}\mathsf{pr}_{\Gamma'_*}\big(\xi_g^*(\mathsf{Q}_{\mathsf{inv}}\otimes 1^{\otimes 7})\cdot \prod_{(h,h')\in \mathrm{Im}(f)\cap \mathrm{Im}(g)} (-\psi_h - \psi_{h'}) \big )\in \ThreeSpin_{2,14}^*\,.
     \]
     Now note that 
     \[
     \xi_{\Gamma'*}\mathsf{pr}_{\Gamma'_*}= \mathsf{pr}_{w*}\xi_{f*}.
     \]
     Combining the results from the two types of $(\mathsf{A},\mathsf{A})$ structures, we have 
     \[
     \mathsf{pr}_{w*}\xi_{\mathsf{A}}^*\alpha_k\notin \ThreeSpin_{2,14}^*\,,
     \]
     and the result follows from Lemma \ref{lem:onefactor}.
\end{proof}

For $0\leq k\leq 6$, let $\mathsf{E}=\mathsf{E}_{3+2k,13-2k}$ be the stable graph of genus $3+2k$ with $13-2k$ legs corresponding to the locus of elliptic tails. We consider the gluing morphism 
\[
\xi_{\mathsf{E}}:\Mbar_{2+2k,14-2k}\times \Mbar_{1,1}\rightarrow \Mbar_{3+2k,13-2k}\,.
\]
\begin{lem}\label{lem:odd-sunflower}
We have
$\xi_{\mathsf{E}*}(\alpha_k\otimes1)\notin\ThreeSpin^*_{3+2k,13-2k}$.
\end{lem}

\begin{proof}

Let $v$ be the genus $2+2k$ vertex of $\mathsf{E}$. We claim that
\begin{equation}\label{eq:odd-sunflower-recovery}
\mathsf{pr}_{v*}\xi_{\mathsf{E}}^*\xi_{\mathsf{E}*}(\alpha_k\otimes1)=-\frac{1}{24}\alpha_k\,.
\end{equation}
By Lemma~\ref{lem:sunflower}, $\alpha_k$ is not a $3$-spin relation. Lemma~\ref{lem:onefactor} and \eqref{eq:odd-sunflower-recovery} therefore imply the same for $\xi_{\mathsf{E}*}(\alpha_k\otimes1)$.

As in the proof of Lemma~\ref{lem:commutator}, there are two isomorphism classes of generic $(\mathsf{E},\mathsf{E})$ structures. The first isomorphism class is $\mathsf{E}$ itself with identity morphisms. We have
\[
\mathsf{pr}_{v*}(\mathsf{Cont}_{\mathsf{E}})=\mathsf{pr}_{v*}(-\psi\cdot\alpha_k\otimes1-\alpha_k\otimes\psi)=-\frac{1}{24}\alpha_k\,.
\]

The second isomorphism class has underlying graph $\Gamma$ with three vertices, two of genus $1$ and one of genus $1+2k$. There are two edges, one connecting each genus $1$ vertex to the genus $1+2k$ vertex. All legs are on the genus $1+2k$ vertex. The two edge contractions $f,g:\Gamma\rightarrow \mathsf{E}$ contract different edges. The gluing map
\[
\xi_f:\Mbar_{1,1}\times \Mbar_{1+2k,15-2k}\times \Mbar_{1,1}\rightarrow \Mbar_{2+2k,14-2k}\times \Mbar_{1,1}
\]
can be written as $\xi_f=\xi_{\mathsf{B}}\times \id$, where $\xi_{\mathsf{B}}$ is the gluing map
\[
\Mbar_{\mathsf{B}}=\Mbar_{1,1}\times \Mbar_{1+2k,15-2k}\rightarrow \Mbar_{2+2k,14-2k}\,.
\]

It suffices to check that $\mathsf{pr}_{v*}\xi_\mathsf{B}^*\alpha_k = 0$, where as usual $\mathsf{pr}_{v*}$ integrates the $\Mbar_{1,1}$ factor. Let $\mathsf{A}=\mathsf{A}_{2+2k,14-2k}$, so $\alpha_k = \xi_{\mathsf{A}*}(\mathsf{Q}_{\mathsf{inv}}\otimes1^{\otimes k})$. The usual formula (valid in the strata algebra) for the operator composition $\xi_{\mathsf{B}}^*\xi_{\mathsf{A}*}$ expresses the result as a sum over generic $(\mathsf{A},\mathsf{B})$-structures $\mathsf{G}$. Since $\mathsf{A}$ has no separating edges and the only edge of $\mathsf{B}$ is separating, no edge of $\mathsf{G}$ can correspond to an edge of both $\mathsf{A}$ and $\mathsf{B}$. Thus there are no excess intersections and we simply have
\[
\xi_{\mathsf{B}}^*\xi_{\mathsf{A}*} = \sum_{\mathsf{G}} (\xi_{\mathsf{G}\to\mathsf{B}})_*\xi_{\mathsf{G}\to\mathsf{A}}^*,
\]
where $\xi_{\mathsf{G}\to\mathsf{A}},\xi_{\mathsf{G}\to\mathsf{B}}$ are the gluing maps corresponding to the contractions $\mathsf{G}\to\mathsf{A},\mathsf{G}\to\mathsf{B}$.

Therefore the class in the $\Mbar_{1,1}$ factor of $\xi_{\mathsf{B}}^*\xi_{\mathsf{A}*}(\mathsf{Q}_{\mathsf{inv}}\otimes1^{\otimes k})$ comes from pulling back either $\mathsf{Q}_{\mathsf{inv}}$ or $1$ via a simple elliptic tail gluing map. In the former case the pullback is zero because the boundary strata used to define $\mathsf{Q}$ do not have any separating edges or vertices of positive genus (this argument was already used in the proof of Lemma~\ref{lem:sep_pullbacks}), while the latter case gives a copy of $1$ that integrates to zero when $\mathsf{pr}_{v*}$ is applied. This completes the proof that $\mathsf{pr}_{v*}(\mathsf{Cont}_\Gamma) = 0$, and we are left with \eqref{eq:odd-sunflower-recovery}.
\end{proof}

\begin{proof}[Proof of Theorem \ref{thm:highergenus}]
For part (i), the result follows from combining Theorems \ref{qth_rel}, \ref{qth_not3spin}, and Lemmas \ref{lem:addmarkings}, \ref{lem:addgenus}, \ref{lem:sunflower}, and \ref{lem:odd-sunflower}.

The argument for part (ii) is essentially the same. First, note that if $\gamma\neq 0\in \R^*(\Mb_{g,n})$, then $\pi^*\gamma \neq 0 \in \R^*(\Mb_{g,n+1})$ by the same argument as in Lemma \ref{lem:addmarkings}. Similarly, if $\gamma \neq 0 \in \R^*(\Mb_{g})$ but $\gamma=0\in \RH^*(\Mb_g)$, then $\xi_{\mathsf{E}*}(\pi^*\gamma\otimes 1)\neq 0\in \R^*(\Mb_{g+1})$ by the same argument as in Lemma \ref{lem:addgenus}, using that $\mathsf{ChowRel}_2^*=\mathsf{CohRel}_2^*$. The calculation in Lemma \ref{lem:sunflower} implies that for $1\leq k\leq 7$, $\alpha_k\neq 0\in \R^*(\Mb_{2+2k,14-2k})$ if $\mathsf{Q}_{\mathsf{inv}}\neq 0\in \R^*(\Mb_{2,14})$. Similarly, the calculation in Lemma \ref{lem:odd-sunflower} gives $\xi_{\mathsf{E}*}(\alpha_k\otimes1)\neq 0\in \R^*(\Mb_{3+2k,13-2k})$ for $0\leq k\leq 6$ if $\mathsf{Q}_{\mathsf{inv}}\neq 0\in \R^*(\Mb_{2,14})$. The result now follows from Theorem \ref{qth_rel}.
\end{proof}

\begin{rem}
    The proofs of the statements in Section \ref{sec:highergenus} 
    are quite similar to the proofs in \cite[Section 8]{CLPWholomorphic}, which is not a coincidence. The holomorphic forms studied therein come from  $\gr^W_{17}\HH^3(\M_2^{\ct},\mathbb{V}_{7,7})$, whereas $\gr^W_{18}\HH^3(\M_2^{\ct},\mathbb{V}_{7,7})$ gives rise to $\QQQ$. The $\SSS_{14}$-representation theory is thus the same in both settings.
\end{rem}

\bibliographystyle{plain}
\bibliography{main}

\end{document}